\documentclass{aims}
\usepackage{lipsum}
\usepackage{amssymb,amsthm}
\usepackage[foot]{amsaddr}
\usepackage{amsmath}
\usepackage{tikz}
\usetikzlibrary{quotes,angles}
\usepackage{graphicx}
\usepackage{subcaption}
\usepackage[shortlabels]{enumitem}

\usetikzlibrary{arrows}
\usepackage{ulem}
\usepackage{float}
\usepackage{graphicx}
\usepackage[section]{placeins}
\graphicspath{{AnalBragg/}}
\usetikzlibrary{decorations.pathmorphing}
\tikzset{snake it/.style={decorate, decoration=snake}}

\usepackage{sansmath}
\usetikzlibrary{shadings,intersections}

\makeatletter
\newcommand*{\rom}[1]{\expandafter\@slowromancap\romannumeral #1@}
 \usepackage[colorlinks=true]{hyperref}
\hypersetup{urlcolor=blue, citecolor=red}

\def\currentvolume{X}
 \def\currentissue{X}
  \def\currentyear{200X}
   \def\currentmonth{XX}
    \def\ppages{X--XX}

\numberwithin{equation}{section}

\theoremstyle{plain}
\newtheorem{theorem}{Theorem}
\newtheorem{corollary}[theorem]{Corollary}
\numberwithin{theorem}{section}

\theoremstyle{definition}
\newtheorem{definition}[theorem]{Definition}
\theoremstyle{remark}
\newtheorem{remark}[theorem]{Remark}
\theoremstyle{remark}
\newtheorem{discussion}[theorem]{Discussion}

\DeclareMathOperator*{\argmin}{arg\,min}

\newcommand{\abs}[1]{\left|#1\right|}                 
\newcommand{\paren}[1]{\left(#1\right)}               

\newcommand{\vp}{\varphi}

\newcommand{\vb}{\mathbf{b}}
\newcommand{\vc}{\mathbf{c}}

\newcommand{\vv}{\mathbf{v}}

\newcommand{\vx}{\mathbf{x}}
\newcommand{\vy}{\mathbf{y}}

\newcommand{\vd}{\mathbf{d}}
\newcommand{\vs}{\mathbf{s}}
\newcommand{\vr}{\mathbf{r}}

\newcommand{\be}{\begin{equation}}
\newcommand{\ee}{\end{equation}}
\newcommand{\bea}{\begin{eqnarray}}
\newcommand{\eea}{\end{eqnarray}}
\newcommand{\bean}{\begin{eqnarray*}}
\newcommand{\eean}{\end{eqnarray*}}

\newcommand{\bel}[1]{\begin{equation}\label{#1}}

\newcommand{\eel}[1]{{\label{#1}\end{equation}}}

\title[Bragg scatter imaging] 
      {Bragg scattering tomography}

\author[James W. Webber and Eric L. Miller]{}


 \email{james.webber@tufts.edu}
 \email{eric.miller@tufts.edu}

\thanks{This material is supported by the U.S.\ Department of
Homeland Security, Science and Technology Directorate, Office of
University Programs, under Grant Award 2013-ST-061-ED0001.}

\thanks{$^*$ Corresponding author: James W. Webber}

\begin{document}
\maketitle

\centerline{\scshape James W. Webber$^*$}
\medskip
{\footnotesize
 \centerline{161 College Avenue}
   \centerline{Halligan Hall}
   \centerline{Medford, MA 02155, USA}
} 

\medskip

\centerline{\scshape Eric L. Miller}
\medskip
{\footnotesize
 \centerline{161 College Avenue}
   \centerline{Halligan Hall}
   \centerline{Medford, MA 02155, USA}
}

\bigskip

 \centerline{(Communicated by Mikko Salo)}

\begin{abstract}
Here we introduce a new forward model and imaging modality for Bragg Scattering Tomography (BST). The model we propose is based on an X-ray portal scanner with linear detector collimation, currently being developed for use in airport baggage screening. The geometry under consideration leads us to a novel two-dimensional inverse problem, where we aim to reconstruct the Bragg scattering differential cross section function from its integrals over a set of symmetric $C^2$ curves in the plane. The integral transform which describes the forward problem in BST is a new type of Radon transform, which we introduce and denote as the Bragg transform. We provide new {injectivity results} for the Bragg transform here, and describe how the conditions of our theorems can be applied to assist in the machine design of the portal scanner. Further we provide an extension of our results to $n$-dimensions, where a generalization of the Bragg transform is introduced. Here we aim to reconstruct a real valued function on $\mathbb{R}^{n+1}$ from its integrals over $n$-dimensional surfaces of revolution of $C^2$ curves embedded in $\mathbb{R}^{n+1}$. {Injectivity proofs} are provided also for the generalized Bragg transform. 
\end{abstract}

\section{Introduction} 
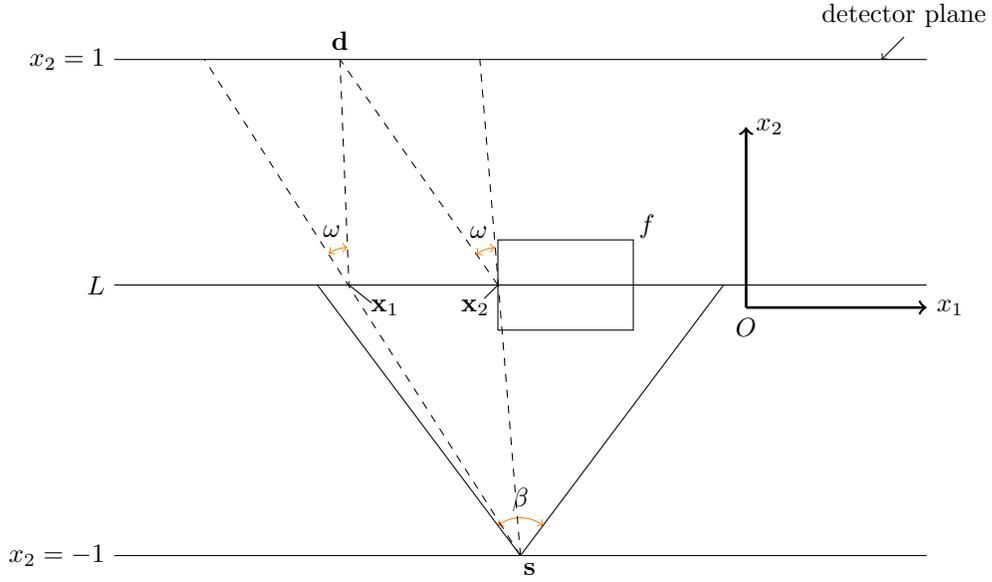
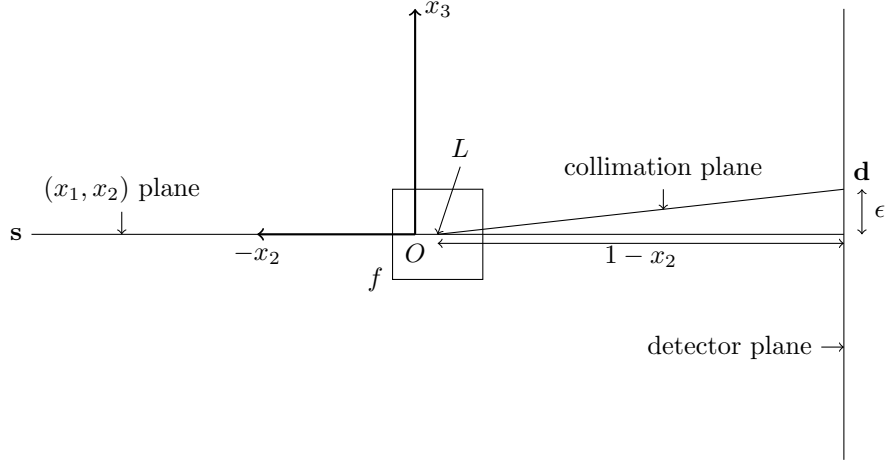
\begin{figure}[!h]
\begin{subfigure}{1\linewidth}
\centering
\begin{tikzpicture}[scale=6]
\draw  (-0.9,0.5)node[left] {$L$}--(0.9,0.5);
\draw  (-0.9,1)node[left] {$x_2=1$}--(0.9,1);
\draw  [->,line width=1pt] (0.5,0.45)--(0.9,0.45)node[right] {$x_1$};
\draw [->,line width=1pt]  (0.5,0.45)node[below] {$O$}--(0.5,0.85)node[right] {$x_2$};
\coordinate (origo) at (0,-0.1);
\coordinate (pivot) at (-0.45,0.5);
\coordinate (bob) at (0.45,0.5);
\draw pic[draw=orange, <->,"$\beta$", angle eccentricity=1.5] {angle = bob--origo--pivot};
\draw [dashed] (0,-0.1)--(-0.7,1);
\draw [dashed] (-0.4,1)--(-0.38,0.5);
\coordinate (origo2) at (-0.38,0.5);
\coordinate (pivot2) at (-0.7,1);
\coordinate (bob2) at (-0.4,1);
\draw pic[draw=orange, <->,"$\omega$", angle eccentricity=1.5] {angle = bob2--origo2--pivot2};
\node at (-0.3,0.45) {$\vx_1$};
\draw (-0.38,0.5)--(-0.33,0.46);
\draw (0,-0.1)--(-0.45,0.5);
\draw [dashed] (0,-0.1)--(-0.09,1); 
\draw (0,-0.1)--(0.45,0.5);
\draw (-0.9,-0.1) node[left] {$x_2=-1$} --(0.9,-0.1);
\draw  (-0.05,0.4) rectangle (0.25,0.6);
\node at (0.28,0.63) {$f$};
\node at (-0.1,0.45) {$\vx_2$};
\draw  (-0.08,0.47)->(-0.05,0.5);
\draw [dashed] (-0.05,0.5)--(-0.4,1)node[above] {$\vd$};
\draw [->] (0.85,1.05)node[above] {detector plane}--(0.8,1);
\node at (0.02,-0.13) {$\vs$};
\coordinate (origo1) at (-0.05,0.5);
\coordinate (pivot1) at (-0.4,1);
\coordinate (bob1) at (-0.09,1);
\draw pic[draw=orange, <->,"$\omega$", angle eccentricity=1.5] {angle = bob1--origo1--pivot1};
\end{tikzpicture}
\caption{$(x_1,x_2)$ (source fan-beam) plane cross-section. The source ($\vs$) opening angle is $\beta$ and we have shown two scattering locations at $\vx_1, \vx_2\in L$ with scattering angle $\omega$.}
\label{fig1}
\end{subfigure}
\begin{subfigure}{1\linewidth}
\centering
\hspace*{0.6cm}
\begin{tikzpicture}[scale=6]
\begin{scope}
\draw  [->,line width=0.8pt] (-0.15,0)--(-0.5,0)node[below] {$-x_2$};
\draw  [->,line width=0.8pt] (-0.15,0)node[below] {$O$}--(-0.15,0.5)node[right] {$x_3$};
\draw  (-1,0) node[left] {$\vs$}--(0.8,0);
\draw  (0.8,-0.5)--(0.8,0.5);
\draw [<->] (0.84,0)--(0.84,0.1)node[above] {$\vd$};
\node at (0.88,0.05) {$\epsilon$};
\draw  (0.8,0.1)--(-0.1,0);
\draw (-0.2,-0.1) node[left] {$f$} rectangle (0,0.1);
\draw [->] (0.4,0.1) node[above] {collimation plane}--(0.4,0.055);
\draw [->] (-0.8,0.05) node[above] {$(x_1,x_2)$ plane}--(-0.8,0);
\draw [->] (0.75,-0.25) node[left] {detector plane}--(0.8,-0.25);
\draw [<->] (-0.1,-0.02)--(0.8,-0.02);
\node at (0.35,-0.05) {$1-x_2$};
\draw [->](-0.05,0.15)node[above]{$L$}--(-0.1,0);
\end{scope}
\end{tikzpicture}
\caption{$(x_2,x_3)$ plane cross-section. Note that $L$ is now orthogonal to the page (parallel to $x_1$).}
\label{figyz}
\end{subfigure}
\caption{The portal scanner geometry. The scanned object is labelled as $f$. The detectors are collimated to planes, and the scattering events occur along lines $L=\{x_2=a,x_3=0\}$, for some $-1<a<1$. The scatter from $L$ is measured by detectors $\vd\in\{x_2=1,x_3=\epsilon\}$, for some $\epsilon>0$.}
\label{figmain}
\end{figure}

In this paper we introduce new Radon transforms in $\mathbb{R}^{n+1}$ which describe the integrals of $L^2$ functions of compact support over the $n$-dimensional surfaces of revolution of a class of $C^2$ curves. A special, motivating case of interest describes the BST problem for an X-ray scanning geometry in airport baggage screening with linear detector collimators, which we will refer to as ``Venetian blind" type collimation. Specifically we focus on the scanning geometry depicted in figure \ref{figmain}. The scanner sources (with coordinate $\vs$) are fixed and switched along the linear array $\{x_2=-1,x_3=0\}$, and are assumed to be polychromatic 2-D fan-beam (in the $(x_1,x_2)$ plane) with opening angle $\beta$. The detectors (with coordinate $\vd$) are assumed to be energy-resolved and lie on the $\{x_2=1\}$ plane, with small (relative to the scanning tunnel size) offset $\epsilon$ in the $x_3$ direction. The detectors are collimated to record photons which scatter on planes in $\mathbb{R}^3$, and the planes of collimation are orientated to intersect the source $(x_1,x_2)$ plane along horizontal lines (parallel to $x_1$). Hence the photon arrivals measured by the portal scanner detectors are scattered from horizontal lines embedded in the $(x_1,x_2)$ plane. An example line of intersection is illustrated by $L$ in figure \ref{figmain}. 

Here we introduce new physical models for the linear detector collimation design which estimate the Bragg scattered signal from line-samples of randomly orientated crystallites (powder scattering). When the effects due to attenuation are ignored, the physical models presented lead us to a new, linear Radon transform, which we denote as the Bragg transform. The removal of attenuation from the modelling is a common assumption made in the scattering tomography literature, for example in Compton Scattering Tomography (CST) \cite{2D1,2D2,2D7,3D1,3D3}. While neglecting the attenuative effects introduces a systematic error in the modelling, the linearization of the model allows us to apply the theory of linear integral equations and Radon transforms to obtain {a solution}. Further the analysis conducted here will likely shed light on the inversion and stability properties in the non-linear case (with attenuation included) and provides the theoretical groundwork required to move forward with such problems. 

The Bragg transform maps the Bragg differential cross section function (the reconstruction target) associated with the crystalline material to its integrals over a set of bounded, symmetric $C^2$ curves in the plane. By exploiting the translational invariance of the Bragg integral data, and using the established theory on linear Volterra integral equations \cite{Tric} and analytic continuation ideas, we prove {the injectivity} of the Bragg transform. {Similar inversion ideas using Volterra equation theory are applied to the circular Radon transform in \cite{ambartsoumian2010inversion}, which has applications in ultrasound imaging.} This work lays the foundation for new imaging techniques in Bragg spectroscopy, a decades old idea, dating back to the experiments of Debye-Scherrer \cite{DS} and the works of \cite{Xray,BE,diff}, which image the crystal structure of point-samples using monochromatic, pencil-beam sources. Our analysis considers the BST problem with polychromatic, fan-beam sources, where the (3-D) crystal sample is observed along lines in the $(x_1,x_2)$ plane. The restriction of the scatter to lines occurs due to the linear collimation technology. The energy-resolved capabilities of the detectors (a technology not available to \cite{DS,Xray,BE,diff}) allow us also to distinguish the energies of the spectrum, which adds a new dimensionality to the data, and will play a crucial role in the inversion ideas presented in our main theorems (e.g. in Theorem \ref{main1}). 

As a natural continuation to the inversion results presented on the Bragg transform, we introduce a generalization to the Bragg transform in section \ref{ndim}. The {injectivity results} presented prior to section \ref{ndim} consider a specific set of symmetric $q_1\in C^2([0,\infty))$ embedded in $\mathbb{R}^2$ with relevance in BST. The generalized Bragg inversion broadens this result to a larger class of $q_1\in C^2([0,\infty))$, and to $(n+1)$-dimensions (for $n\geq 1$), where the integrals are taken over the $n$-dimensional surfaces of revolution of $q_1$. The idea follows similar intuition to that of CST. In 2-D CST the Radon transforms of interest take integrals of the electron density over toric sections \cite{2D1,2D5,2D8}, which is generalized in 3-D to the integrals over spindle tori \cite{3D3,3D6,ctrans}, namely the surfaces of revolution of toric sections about their axis of reflection. See also \cite{2D2,2D3,2D4,3D1,3D4} for more work on generalized Radon transforms in CST. The $n=2$ case considered here follows similar ideas to that of \cite[Theorem 5.1]{ctrans}, which proves the injectivity of a generalized apple transform. The results presented in this paper however consider a different integral kernel to that of \cite{ctrans} {and we prove injectivity for $n\geq 1$} (\cite{ctrans} considers $n=2$) {on the domain of compactly supported $L^2$ functions}. We also consider both monotone decreasing and increasing curves in the $n=1$ case (\cite{ctrans} covers increasing curves). 

{The main purpose of the theorems presented here is to prove injectivity of the Bragg transforms presented (i.e. to show that the solution is unique) and to provide some insight into the problem stability and the portal scanner dimensions required for unique inversion. Proof of injectivity is important since it removes concern for image artifacts due to null space, and, with sufficient data samples and noise reduction, it provides mathematical guarantee that the reconstructed cross section is a sensible representation of the ground truth. 

Injectivity is proven by showing that the Bragg differential cross section function can be expressed explicitly as an infinite sum, using Neumann series and Volterra equation theory. The proofs presented thus lay out novel inversion methods for BST. The inversion methods detailed in the main theorems are not recommended for practical use however, due to unstable steps discovered in the inversion process. These are discussed in detail in Remark \ref{rem1}. Also, to prove injectivity, we consider a lower dimensional subset of the full data (the complete data is 4-D and the data used in our theorems is 3-D). 
So as not to discard any data, and to avoid the instabilities discussed in Remark \ref{rem1}, we choose not to apply the inversion methods presented in the main theorems using Volterra equation theory. Instead, for practical reconstruction, we introduce an algebraic reconstruction algorithm using Total Variation (TV) regularization and non-negativity constraints in section \ref{results}. To test our algorithm, we present image reconstructions of two imaging phantoms which consist of multiple spheres of crystalline powder arranged in a row. The spheres are made up of materials such as graphite, diamond, salt and Aluminum. Our algorithm is shown to offer high quality results with $F_1$ scores $F_1>0.8$ in most cases, and we see no evidence of any artifacts due to null space, which provides further validation of our injectivity theorems in a practical setting.}

The literature {includes} the inversion of generalized Radon transforms \cite{gen1,gen2,gen3,gen4,Co1981,Co1987,GU1989}. In \cite{gen1} the authors present a uniform reconstruction formula for the generalized Radon transform which describes the integrals of continuous functions over hypersurfaces with ``regular" generating function. See \cite[page 2]{gen1} for the definition of a regular generating function. The formulae of \cite{gen1} do not apply to our problem however since the hypersurfaces we consider fail to satisfy the regularity condition. We give examples as to why this is this case later in section \ref{ndim}, when the notation of the paper is settled upon. The works of \cite{gen2,gen3} present inversion formulae{, up to modulo smoothing (i.e. a microlocal inversion), for a class of generalized Radon transforms.} A microlocal inversion, dubbed an ``almost" inversion by Belykin \cite[page 584]{gen3}, recovers the locations and directions of the image singularities. The smooth parts  however are undetermined. We provide {exact inversion results} for the solution here, which recover both the smooth and non-smooth parts (including singularities) of the target function. In \cite{Co1981}, the integrals over $\alpha$ and $\beta$ type curves are considered, under rotational invariance, and a solution is obtained through expansion into the Fourier series. The Radon transforms we consider satisfy translational invariance, and we reach a solution after applying the Fourier transform. In \cite{gen4} the authors consider the inversion of Radon transforms which define the integrals over weighted hyperplanes with general defining measure. The hypersurfaces we consider are of a different class to that of hyperplanes and hence the theory we present does not fall into the framework of \cite{gen4}. {In \cite{webber2020microlocal} the authors consider the stability properties of some generalized Radon transforms from scattering tomography from a microlocal perspective, which includes an analysis of the Bragg transform presented here. No injectivity results are provided however. We aim to address this here. Later in section \ref{results} we provide some numerical evidence for the theory of \cite{webber2020microlocal} (in particular Theorem (3.10)), and our results are in line with the authors stability estimates.}

The rest of this paper is organized as follows. In section \ref{pre} we state our notation, and some preliminary results on the Fourier transform and Volterra integral equations. In section \ref{3.1} we define the Bragg differential cross section function $f$ (the target for reconstruction) and state the equations from scattering physics which will be used in the physical modeling. We then show, in section \ref{sec3.1}, how the physical models translate to the geometry of figure \ref{figmain}. 

In section \ref{bgsec1} we introduce the Bragg transform $\mathfrak{B}$ and show how, under certain physical assumptions, the Bragg signal can be approximated by $\mathfrak{B}f$. We then prove that $\mathfrak{B}$ is invertible and bounded. The physical assumptions made are to neglect attenuation (as discussed earlier), and to set $\epsilon=0$ (see figure \ref{figyz}) when calculating the Bragg angle. In section \ref{ndim} we introduce the generalized Bragg transform $\mathfrak{B}_n$ to $\mathbb{R}^{n+1}$ (for $n\geq 1$), as a natural extension of the results of section \ref{bgsec1}. The injectivity of $\mathfrak{B}_n$ is proven thereafter. 

In section \ref{6.2} we show how, under certain restrictions to the machine design, we can lift the $\epsilon=0$ assumption of section \ref{bgsec1}. Here we introduce the offset Bragg transform $\mathfrak{B}_{\epsilon}$ which models the Bragg intensity for $\epsilon\geq 0$, excluding only the attenuative effects from the exact model.  We then show (in Theorem \ref{main4}) that $\mathfrak{B}_{\epsilon}$ is invertible when the machine design conditions are satisfied. Roughly speaking, the design conditions specify that the detector offset $\epsilon$ not be too large relative to the source fan width $w$. To finish the paper, in section \ref{MD} we provide example machine parameters which satisfy the conditions of Theorem \ref{main4}. Here we consider the ranges of $w$, $\epsilon$ and energy which allow for inversion and illustrate example detector configurations. {Finally, in section \ref{results}, we provide image reconstructions from simulated Bragg scatter data measured in the acquisition geometry of figure \ref{figmain}.}

\section{Preliminary results and definitions}
\label{pre}
Here we state our notational conventions, and some definitions and preliminary results which will be used in our theorems. Throughout this paper we denote:
\begin{enumerate}[label=(\roman*)]
\item $L^2_0(\Omega)$ as the set of $L^2$ functions with compact support on $\Omega\subset\mathbb{R}^n$.
\item $C^k(\Omega)$ as the set of $k$-continuously differentiable functions on $\Omega\subset\mathbb{R}^n$.
\item $\mathbb{R}_+=(0,\infty)$ as the set of positive real numbers not including zero.
\item $\mathcal{I}=(-1,1)$ as the interval length two, center zero with the end-points removed.
\end{enumerate}
\noindent Now we have the definition of the Fourier transform in terms of angular frequency.
\begin{definition}
Let $f\in L^2(\mathbb{R}^n)$. Then we define the Fourier transform $\hat{f}$ of $f$ in terms of angular frequency
\begin{equation}
\hat{f}(\boldsymbol{\eta})=(2\pi)^{-n/2}\int_{\mathbb{R}^n}f(\vx)e^{-i\vx\cdot \boldsymbol{\eta}}\mathrm{d}\vx.
\end{equation}
\end{definition}
\noindent We state the Plancherel theorem \cite{plan}.
\begin{theorem}[Plancherel theorem]
\label{plan}
Let $f\in L^2(\mathbb{R}^n)$. Then $\hat{f}\in L^2(\mathbb{R}^n)$ and
\begin{equation}
\|f\|_{L^2(\mathbb{R}^n)}=\|\hat{f}\|_{L^2(\mathbb{R}^n)}.
\end{equation}
\end{theorem}
\noindent We now state some results on Volterra type integral equations from \cite[page 10]{Tric}.
\begin{definition}
We define a Volterra equation of the second kind to be an equation of the form
\begin{equation}
g(x)=\lambda\int_0^xK(x,y)f(y)\mathrm{d}y+f(x)
\end{equation}
with real valued kernel $K$ on a triangle $T'=\{0<x<x', 0<y<x\}$. $K$ is said to be an $L^2$ kernel if
\begin{equation}
\|K\|^2_{L^2(T')}=\int_0^{x'}\int_0^x K^2(x,y)\mathrm{d}y\mathrm{d}x\leq N^2
\end{equation}
for some $N>0$.
\end{definition}
\begin{theorem}
\label{volt}
Let $g\in L^2([0,x'])$ and let $K(x,y)$ be an $L^2$ kernel on $T'=\{0<x<x', 0<y<x\}$ for some $x'>0$. Then the Volterra integral equation of the second kind
\begin{equation}
g(x)=\lambda\int_0^xK(x,y)f(y)\mathrm{d}y+f(x)
\end{equation}
has one and only one solution in $L^2([0,x'])$, and the solution is given by the formula
\begin{equation}
f(x)=\lambda\int_0^xH(x,y;\lambda)g(y)\mathrm{d}y+g(x),
\end{equation}
where
\begin{equation}
H(x,y;\lambda)=\sum_{l=0}^{\infty}{(-1)^l}\lambda^{l}K_{l+1}(x,y)
\end{equation}
and the iterated kernels $K_{l}$ are defined by $K_1(x,y)=K(x,y)$ and
$$K_{l+1}(x,y)=\int_0^xK(x,z)K_{l}(z,y)\mathrm{d}z$$
for $l\geq 1$.
\end{theorem}
\section{The forward Bragg scattering model}
\label{sec3}
Here we present the equations from scattering theory used for the modeling of the Bragg intensity in the portal scanner geometry (of figure \ref{figmain}). We assume that only single scatter effects occurs and neglect the effects due to multiple scattering. We further assume a polychromatic fan-beam source with finite width $w$ and detectors with energy-resolving capabilities, so we can distinguish between the photon energies of the spectrum. 
\subsection{The physical model}
\label{3.1}
\begin{figure}[!htb]
\centering
\begin{tikzpicture}[scale=5]
\path (1.5,0,0.1) coordinate (S);
\path (0,0.5,0.5) coordinate (D);
\path (0.5,0,0.75) coordinate (w);
\path (0.25,0,0.9125) coordinate (a);
\draw [thin,dashed] (w)--(a);
\node at (0.05,0.3,0.5) {$E_s$};
\node at (1.1,0,0.6) {$E$};
\draw [fill] (S) circle (0.01);
\draw pic[draw=orange, <->,"$\omega$", angle eccentricity=1.5] {angle = D--w--a};
\draw[thick,->] (0,0,0) -- (1,0,0) node[above]{$x_2$};
\draw[thick,->] (0,0,0) -- (0,1,0) node[anchor=north west]{$x_3$};
\draw[thick,->] (0,0,0) -- (0,0,1) node[anchor=south]{$x_1$};
\draw (0,0.45,0.45)--(0,0.55,0.45);
\draw (0,0.45,0.45)--(0,0.45,0.55);
\draw (0,0.55,0.55)--(0,0.45,0.55);
\draw (0,0.55,0.55)--(0,0.55,0.45);
\draw[thick,->] (0,0.5,0.5)--(0.1,0.5,0.5)node[above]{$\vv$};
\draw [snake it, ->,blue](S)node[right]{\textcolor{black}{$\vs$ (point source)}}--(w)node[below]{\textcolor{black}{$\vx$}};
\draw [snake it, ->,red](w)--(D);
\draw [->] (0,0.75,0.75)node[above]{$\vd$ (area $D_A$)}--(0,0.55,0.55);
\end{tikzpicture}
\caption{A scattering event occurs at a scattering site $\vx$, for photons emitted from a source $\vs$ and recorded at a detector $\vd$. The initial photon energy is $E$ and the scattered energy is $E_s$. Here $\vv$ is the direction normal to the detector surface, displayed as a square in the $x_1x_3$ plane. The scattering angle is $\omega=2\theta$. where $\theta$ is the Bragg angle.}
\label{fig3}
\end{figure}
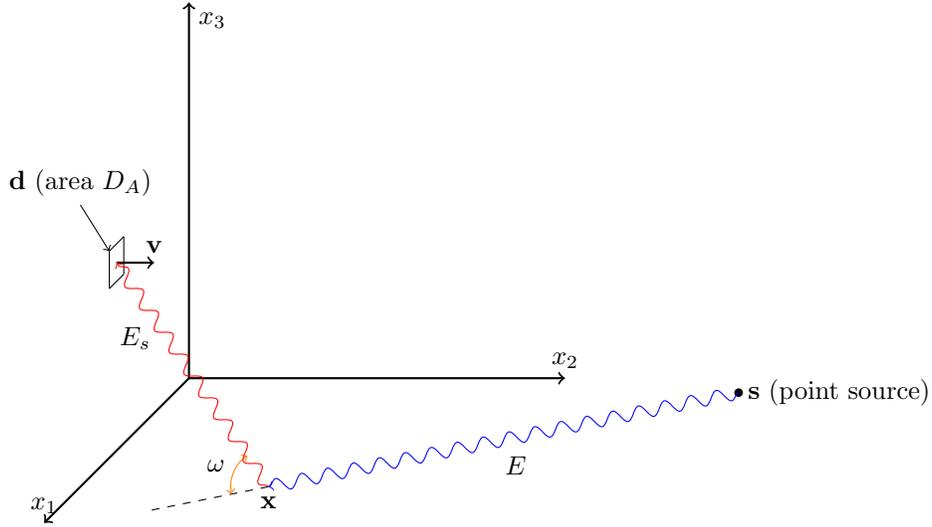

Consider the scattering event pictured in figure \ref{fig3}. The scattering contribution from the scattering site $\vx$ is \cite[page 73]{model} {(see also \cite{greenberg2013snapshot,maccabe2012pencil})}
\begin{equation}
\label{equ17}
\begin{split}
I\left(E_s,\vs,\vd\right)=I_0\left(E\right)e^{-\int_{L_{\vs\vx}}\mu(E,Z)}
n_c\left(\vx\right)\mathrm{d}V\times \frac{\mathrm{d}\sigma_{\textbf{e}}}{\mathrm{d}\Omega}\left(E,\omega,Z\right)e^{-\int_{L_{\vx\vd}}\mu(E_s,Z)}
\mathrm{d}\Omega_{\vx,\vd},
\end{split}
\end{equation}
where $Z=Z(\vx)$ denotes the effective atomic number of the material (as a function of $\vx$), $I_0$ is the initial intensity, $\omega$ is the scattering angle and $L_{\vs\vx}$ and $L_{\vx\vd}$ are the line segments connecting $\vs$ to $\vx$ and $\vx$ to $\vd$ respectively. 
The line integrals in the exponents are taken over the linear attenuation coefficient $\mu$, which depends on $Z$ and the photon energy $E$. Since energy and wavelength $\lambda$ are inversely proportional by $E=\frac{hc}{\lambda}$, where $h$ is Planck's constant and $c$ is the speed of light in a vacuum, we use the convention in spectroscopy and set $1\AA^{-1}\approx 12.4$keV, where $\AA=10^{-8}$cm denotes one Angstrom. That is we will use the terms ``energy" and ``inverse length" interchangeably with a conversion factor of 12.4. The number density (number of cells or electrons per unit volume)  is denoted by $n_c$ and $\mathrm{d}V$ is the volume measure. The solid angle $\mathrm{d}\Omega_{\vx,\vd}$ subtended by $\vx$ and $\vd$ is given by
\begin{equation}
\mathrm{d}\Omega_{\vx,\vd}=D_A\times\frac{\left|\vr\cdot \vv\right|}{|\vr|^3},
\end{equation}
where $D_A$ is the detector area, $\vr=\vx-\vd$ and $\vv$ is the direction normal to the detector surface (see figure \ref{fig3}). The differential cross section $\frac{\mathrm{d}\sigma_c}{\mathrm{d}\Omega}\left(E,\omega,Z\right)$ describes the angular distribution (in $\omega$) of scattering events from a material $Z$ with incident photon energy $E$.

The Bragg--Laue (total) cross section is \cite{ITC,paper}
\begin{equation}
\label{equ1}
\sigma_{\textbf{e}}(E,Z)=\frac{r^2_0}{2E^2a^3_0(Z)}\sum_{H\in \mathcal{H}}P(\theta)d_{H}j_{H}\left|F_{H}\left(\frac{1}{2d_H}\right)\right|^2e^{-2M}.
\end{equation}
Here the sum is over all Miller indices
$$\mathcal{H}=\left\{H=(h,k,l)\in\mathbb{N}^3\backslash\{(0,0,0)\} : \frac{1}{2d_HE}<1\right\},$$
which each correspond to a reflection at an angle $\omega=2\theta$ determined by the Bragg equation
\begin{equation}
\label{braggequ}
 \frac{1}{E}=2d_{H}\sin\theta,
\end{equation}
where for cubic structures
\begin{equation}
\label{dH}
d_{H}=\frac{a_0(Z)}{\sqrt{h^2+k^2+l^2}},
\end{equation}
is the spacing between the reflection planes within the crystal and $\theta$ is the Bragg angle. For non-cubic structures (e.g. hexagonal, tetragonal) there also exist explicit formulae for $d_H$. In equation \eqref{dH} $a_0$ denotes the uniform lattice spacing of the crystal. The higher order reflections occur at integer multiples of the wavelength $n\lambda=n/E$ and correspond to reflections with $1/n$ times the $d_H$-spacing. As Bragg diffraction is a coherent scattering event there is no decrease in energy and $E=E_s$. The scattering factor $F_H$ is defined as
\begin{equation}
\label{FH}
F_H\left(q\right)=\sum_{i=1}^{n_a}F_i\left(q\right) e^{-2\pi i \vx_i\cdot H},
\end{equation}
where $n_a$ is the number of atoms in a single cell of the crystal, the $\vx_i\in [0,1]^3$ are the coordinates of the atoms within the cell and $F_i$ is the atomic form factor \cite{hubbell,hubbell1} of atom $i$. The momentum transfer $q$ is defined by 
\begin{equation}
\boxed{q=E\sin\theta,}
\end{equation}
which has units of $\AA^{-1}$.

The polarisation factor $P(\theta)$ is given by
\begin{equation}
P(\theta)=\frac{1+\cos^22\theta}{2}.
\end{equation}
The remaining terms in equation (\ref{equ1}), $j_H$ and $e^{-2M}$, account for the multiplicity factors of the powder and the temperature factor respectively. However we neglect effects due to temperature and plane multiplicity and set $j_H=1$ for all $H\in\mathcal{H}$ and $M=0$. We can now write the Bragg--Laue scattering differential cross section as
\begin{equation}
\label{Bgm}
\begin{split}
\frac{\mathrm{d}\sigma_{\textbf{e}}}{\mathrm{d}\Omega}\left(q,E,Z\right)&=P\left(\sin^{-1}\frac{q}{E}\right)F(q,Z)
\end{split}
\end{equation}
for $0<q<E$ where,
\begin{equation}
\label{Bgm1}
\begin{split}
F(q,Z)&=\frac{1}{\pi q}\paren{\frac{r^2_0}{16 a^3_0(Z)}}\sum_{H\in \mathcal{H}}\delta\left(\frac{1}{2d_H}-q\right)d_H\left|F_{H}\left(q\right)\right|^2
\end{split}
\end{equation}
and $\delta$ is the Dirac--delta function. 
Let us parametrize the scattering direction in spherical coordinates, so that the scattering angle $\omega$ is the polar angle and $\varphi$ is the azimuth angle. Then the solid angle element is $\mathrm{d}\Omega=\sin\omega\mathrm{d}\omega\mathrm{d}\varphi=\frac{4q}{E^2}\mathrm{d}q\mathrm{d}\varphi$ and we verify that
\begin{equation}
\int\frac{\mathrm{d}\sigma_{\textbf{e}}}{\mathrm{d}\Omega}\mathrm{d}\Omega=\frac{8\pi}{E^2}\int_{0}^{E}\frac{\mathrm{d}\sigma_{\textbf{e}}}{\mathrm{d}\Omega}\left(q,E,Z\right)q\mathrm{d}q=\sigma_{\textbf{e}}(E,Z).
\end{equation}
\subsection{Model transition to the portal scanner geometry}
\label{sec3.1}
In the acquisition geometry of figure \ref{figmain}, the scattering is restricted to lines parallel to $x_1$ given the Venetian blind type collimation of the detector array. The relation between the horizontal line profile scanned (i.e. $\{x_2=a\}$ for some $a\in\mathcal{I}$) and the detector offset $\epsilon$ (see figure \ref{figyz}) has not yet been discussed. To keep the discussion general, in this paper we consider diffeomorphic mappings from $x_2$ to $\epsilon$. To this end let
\begin{equation}
\label{phi}
\boxed{\Phi : \mathcal{I}\to \Phi(\mathcal{I}),\ \ \text{be an $x_2\to\epsilon$ diffeomorphism},}
\end{equation}
which describes the relation between the detector offset $\epsilon$ and the scattering line profile $x_2$. Throughout this paper we assume that $0\leq \Phi(x_2)\leq M$ for all $x_2\in\mathcal{I}$, for some $M>0$. Then we can model the intensity of Bragg scattered photons measured by the portal scanner detectors as integrals of the point scatterer model \eqref{equ17} over lines parallel to the $x_1$ axis
\begin{equation}
\label{equBG}
\begin{split}
\mathfrak{B}_af(E,s_1,d_1,\Phi(x_2))=C\int_{\mathbb{R}}&\chi_{\left[-w(x_2),w(x_2)\right]}(x_1-s_1)I_0(E,x_1)P(\theta(\vd,\vs,\vx))\mathrm{d}\Omega_{\vx,\vd}\\
&\times f\left(E\sin\theta(\vd,\vs,\vx),\vx\right)A_1(E,\vs,\vx)A_2(E,\vx,\vd)\mathrm{d}x_1,
\end{split}
\end{equation}
where $\vs=(s_1,-1,0)$, $\vd=(d_1,1,\Phi(x_2))$, $\chi_{S}$ denotes the characteristic function on a set $S$ and $f(q,\vx)=n_c(\vx)F(q,\vx)$. Here the spatial variable $\vx=(x_1,x_2)$ determines the effective atomic number $Z$ and the function $F(q,Z)=F(q,\vx)$ of \eqref{Bgm} is written as a function of $\vx$. The source width $w$ is determined by the source opening angle $\beta$ (see figure \ref{fig1})
\begin{equation}
\label{wdef}
\boxed{w(x_2)=(1+x_2)\tan\frac{\beta}{2}.}
\end{equation}
{{Throughout this paper (barring section \ref{ndim}) we use the shorthand notation $w=w(x_2)$, as the dependence of $w$ on $x_2$ is not important to the proofs presented.}} $\beta\in (0,\pi)$ will remain fixed until section \ref{MD} where we consider varying machine configurations.

The solid angle is
\begin{equation}
\label{sangle}
\mathrm{d}\Omega_{\vx,\vd}=D_A\times \frac{((\vx,0)-\vd)\cdot (0,-1,0)^T}{|(\vx,0)-\vd|^3},
\end{equation}
and the Bragg angle ($\theta$) is determined by
\begin{equation}
\label{bangle}
\cos\omega=\cos 2\theta(\vd,\vs,\vx)=\frac{((\vx,0)-\vs)\cdot(\vd-(\vx,0))}{|((\vx,0)-\vs)||(\vd-(\vx,0))|}.
\end{equation}
Here $A_1(E,\vs,\vx)=e^{-\int_{L_{\vs(\vx,0)}}\mu(E,Z)}$ and $A_2(E,\vx,\vd)=e^{-\int_{L_{(\vx,0)\vd}}\mu(E,Z)}$ account for the attenuation of the incoming and scattered rays respectively, $n_c$ is the number of crystal cells per unit volume and $C$ is a small cross sectional area (the area of an image pixel). 

The target for reconstruction is $f(q,\vx)$, which is a function of the momentum transfer $q\in\mathbb{R}_+$ and $\vx\in\mathbb{R}^2$, where $x_2\in\mathcal{I}$ is fixed for each linear detector array (with offset $\epsilon$) considered. So we consider a one-dimensional set of 2-D inverse problems to recover the full 3-D $f$. That is we recover $f(\cdot,\cdot,x_2)$ for every $x_2\in\mathcal{I}$.

By equation \eqref{dH} the maximum value of $d_H$ (for cubic structures) over all Miller indices $H$ is $d_H=a_0$. Hence the minimum value of $q$ for which the Bragg equation is satisfied is $1/2d_H=1/2a_0$. Let $a_M$ be the maximum lattice spacing among the crystals of interest. Then from equation \eqref{Bgm} it follows that $f(q,\vx)=0$ for $q<E_m=1/2a_M$ for the cubic crystals in consideration. One would see a similar bound away from zero on $\text{supp}(f)$ also for non-cubic structures. Hence we aim to recover $f$ in the momentum transfer range $q\in\mathfrak{E}$ for some maximum energy $E_M$ and minimum energy $ E_m$ of interest. Letting $\mathfrak{E}=[E_m,E_M]$ denote the energy range of interest, we model $f\in L^2_0(\mathfrak{E}\times\mathbb{R}\times\mathcal{I})$ as an $L^2$ function of compact support. 
\section{The Bragg transform}
\label{bgsec1}
In this section we consider the injectivity and inversion properties of a new Radon transform which describes the idealized case for the physical model \eqref{equBG}. We assume a negligible detector offset for the theoretical results presented in this section, so the photon transport occurs exclusively within the $(x_1,x_2)$ plane. That is we assume a maximum offset $\epsilon_M=\max_{x_2\in\mathcal{I}}\Phi(x_2)$ with negligible size relative to the scanning tunnel length (i.e. 2 or the length of $\mathcal{I}$), and we approximate $\epsilon=0$ in the calculation of the Bragg angle (equation \eqref{bangle}). In this case there will be systematic errors introduced in the modeling due to the finite detector offset. In section \ref{sdoff}, we show how the $\epsilon=0$ assumption may be lifted to suit $\epsilon>0$. We choose to consider the idealized mathematical model first for the Bragg problem as it leads more naturally into the generalizations to $\mathbb{R}^{n+1}$ presented in section \ref{ndim}. Further the mathematical derivations in the $\epsilon=0$ example are more elegant than those of the $\epsilon>0$ case, and hence provide a better gateway for the reader into the inversion ideas of the paper.

The Bragg data $\mathfrak{B}_a$ is formally overdetermined (the data is 4-D and the target function is 3-D). To prove the invertibility of $\mathfrak{B}_a$ we consider the three-dimensional subset of the data when $s_1=d_1$. We also set the detector offset equal to zero in the calculation of the Bragg angle (as discussed in the last paragraph) and neglect the attenuative effects (setting $A_1=A_2=1$) as discussed in the introduction. With this in mind we define the Bragg transform $\mathfrak{B} : L^2_0(\mathfrak{E}\times\mathbb{R}\times\mathcal{I})\to L^2(\mathbb{R}_+\times\mathbb{R}\times\Phi(\mathcal{I}))$ as
\begin{equation}
\label{bg11}
\begin{split}
\mathfrak{B}f(E,s_1,\Phi(x_2))=\int_{\mathbb{R}}{\chi_{\left[-w(x_2),w(x_2)\right]}}(x_1-s_1)&W(E,x_1-s_1,x_2)\\
&\ \ \ \ \  \times f(Eq_1(x_1-s_1),x_1,x_2)\mathrm{d}x_1,
\end{split}
\end{equation}
where
\begin{equation}
\begin{split}
\label{bcur}
q_1(x_1)&=\sin\theta\paren{(0,-1,0),(0,1,0),\vx}\\
&=\frac{1}{\sqrt{2}}\sqrt{1+\frac{x_1^2-(1-x_2^2)}{\sqrt{x_1^2+(x_2+1)^2}\sqrt{x_1^2+(1-x_2)^2}}}
\end{split}
\end{equation}
after setting $\epsilon=\Phi(x_2)=0$ in the calculation of the Bragg angle $\theta$ (see equation \eqref{bangle}) under the negligible $\epsilon$ assumption. Here $W : \mathbb{R}_+\times\mathbb{R}\times \mathcal{I}\to \mathbb{R}$ is a weighting which accounts for the physical modeling. We will consider general $W$ for now. Later in Corollary \ref{cor1} we consider the specific $W$ which describes the physical modeling terms of \eqref{equBG} (barring attenuation). 

With $W$ appropriately chosen (see Corollary \ref{cor1}), equation \eqref{bg11} approximates the Bragg signal in the following sense
$$\mathfrak{B}f(E,s_1,\Phi(x_2))\approx\mathfrak{B}_af(E,s_1,s_1,\Phi(x_2)),$$
where the approximation error arises from neglecting attenuation and the negligible $\epsilon=\Phi(x_2)$ assumption. 

\begin{figure}[!h]
\centering
\begin{subfigure}{0.49\textwidth}
\includegraphics[width=0.9\linewidth, height=5.2cm, keepaspectratio]{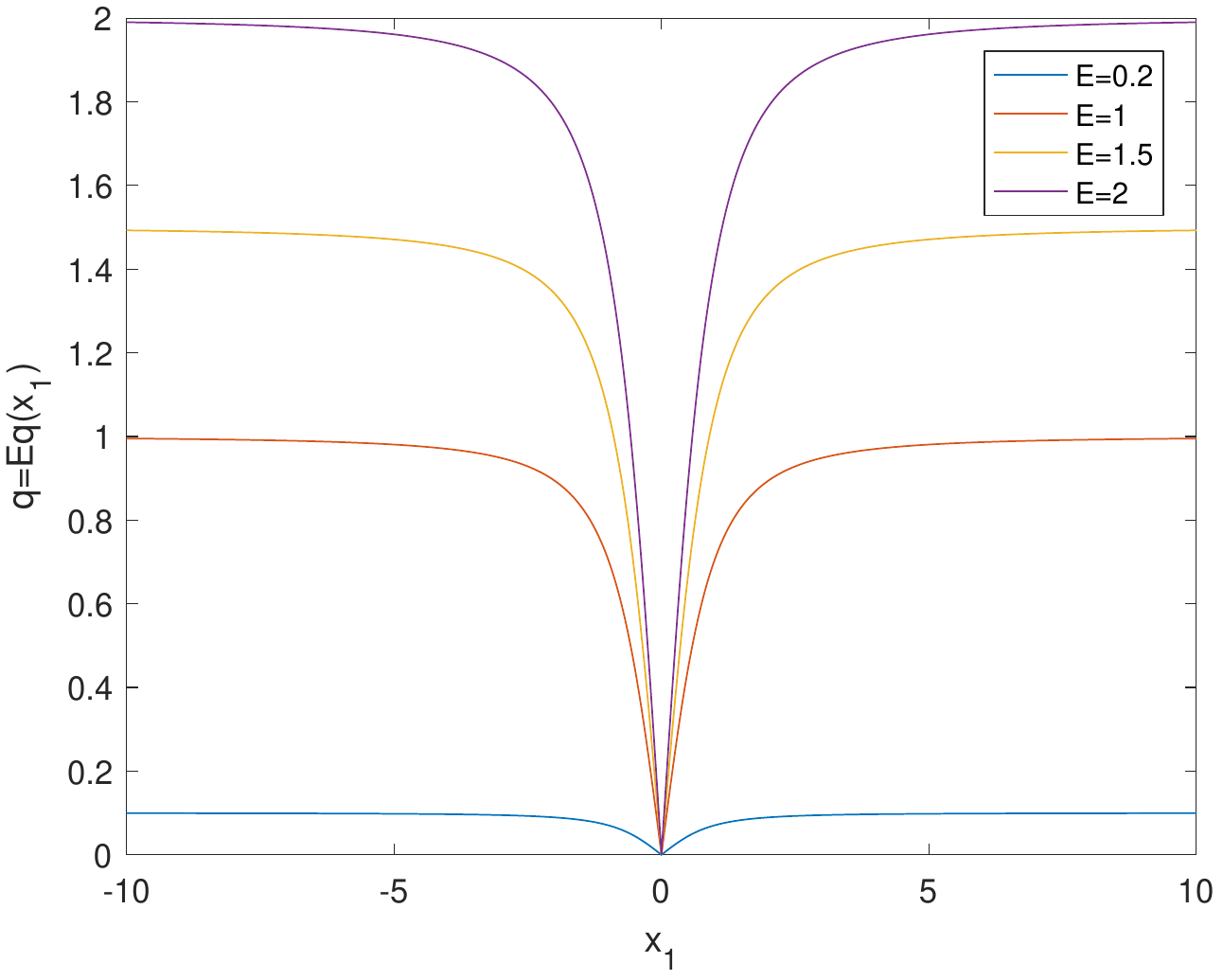}
\subcaption{$x_2=0$}
\end{subfigure}
\begin{subfigure}{0.49\textwidth}
\includegraphics[width=0.9\linewidth, height=5.2cm, keepaspectratio]{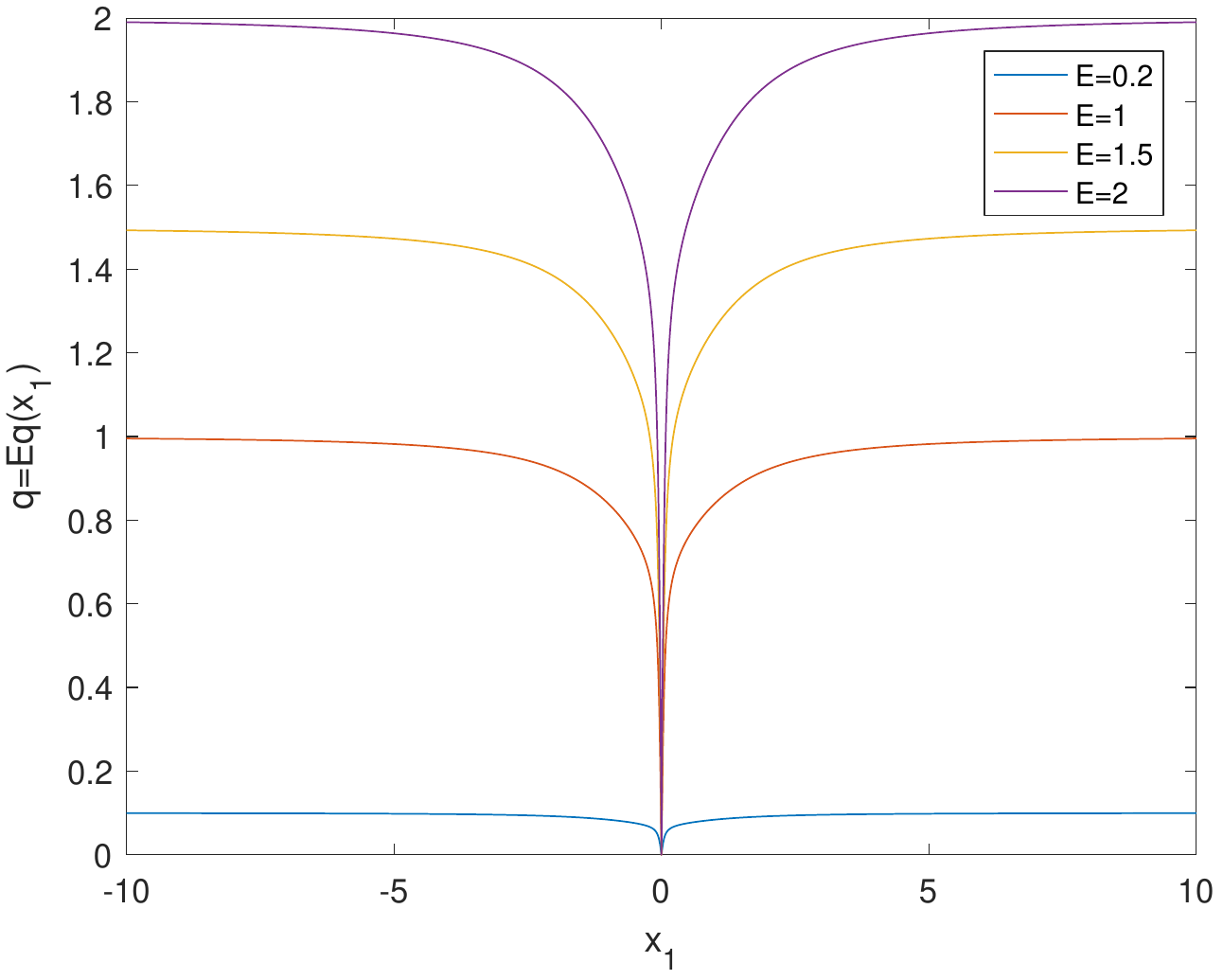}
\subcaption{$x_2=0.95$}
\end{subfigure}
\caption{Plot of the curves of integration for the Bragg transform for varying $E$ and $x_2$. $s_1=0$ is fixed.}
\label{Bgcurve}
\end{figure}

The Bragg transform maps $f$ to its weighted integrals over the set of 1-D curves
\begin{equation}
\label{curve}
\mathcal{Q}=\{(q,x_1)\in\mathbb{R}_+\times\mathbb{R} : q=Eq_1(x_1-s_1), E>0, s_1\in\mathbb{R}\},
\end{equation}
and thus we are modeling the Bragg intensity as the weighted integrals of the differential cross section over the set of curves in the plane $\mathcal{Q}$. See figure \ref{Bgcurve} for a visualization of the curves of integration.

We now have our first main theorem which proves the injectivity of the Bragg transform.
\begin{theorem}
\label{main1}
Let $W:\mathbb{R}_+\times\mathbb{R}\times \mathcal{I}\to\mathbb{R}$ be a separable weighting of the form $W(E,x_1,x_2)=W_1(E)W_2(x_1,x_2)$, where $W_1\in C^1(\mathbb{R}_+)$ and $W_2\in C^1(\mathbb{R}\times \mathcal{I})$ are non-vanishing. Let $W_2$ be symmetric about $x_1=0$ and further let $W_2(\cdot,x_2)\in C^1(\mathbb{R})$ be bounded and have bounded first order derivative on $\left[-w(x_2),w(x_2)\right]$ for any $x_2\in\mathcal{I}$. {Then $\mathfrak{B}$ is injective.}
\begin{proof}
{{We use the shorthand notation $w=w(x_2)$ as specified in section \ref{sec3.1}}}. Substituting $x_1\to x_1+s_1$ in equation \eqref{bg11} yields
\begin{equation}
\label{bg111}
\mathfrak{B}f(E,s_1,\Phi(x_2))=\int_{\mathbb{R}}\chi_{\left[-w(x_2),w(x_2)\right]}(x_1)W(E,x_1,x_2)f(Eq_1(x_1),x_1+s_1,x_2)\mathrm{d}x_1.
\end{equation}
Then, after taking the Fourier transform of both sides in the $s_1$ variable, we have
\begin{equation}
\label{bf1}
\begin{split}
\widehat{\mathfrak{B}f}(E,\eta_1,\Phi(x_2))&=\int_{\mathbb{R}}\chi_{\left[-w(x_2),w(x_2)\right]}(x_1)W(E,x_1,x_2)\hat{f}(Eq_1(x_1),\eta_1,x_2)e^{i\eta_1 x_1}\mathrm{d}x_1\\
&=2W_1(E)\int_{0}^{w(x_2)}W_2(x_1,x_2)\cos(\eta_1 x_1)\hat{f}(Eq_1(x_1),\eta_1,x_2)\mathrm{d}x_1,
\end{split}
\end{equation}
where $\eta_1$ is dual to $s_1$ and the second step of \eqref{bf1} follows due to symmetry of $q_1$ and $W_2$ about $x_1=0$ (see the definition of $q_1$ in \ref{bcur}). 

Let $z=q_1(x_1)$ and let $g(z)=q^{-1}_1(z)=x_1$ be the inverse map of $q_1$. We will now find the closed-form expression for $g$. Refer to figure \ref{figpf}.
\begin{figure}[!h]
\centering
\begin{tikzpicture}[scale=4]
\draw [thin, dashed]  (-1,0.7)node[left] {$x_2$}--(1,0.7);
\draw  (-1,1)node[left] {$x_2=1$}--(1,1);
\draw  [->] (0,0.45)--(0.7,0.45)node[right] {$x_1$};
\draw [->]  (0,0.45)node[below] {$O$}--(0,1.15)node[right] {$x_2$};
\draw (-0.3,0.45) circle [radius=0.625];
\draw (0,-0.1)--(0.27,0.7)node[right] {$\vx$};
\draw (0.27,0.7)--(0.4,1);
\draw (0.27,0.7)--(0,1);
\coordinate (origo1) at (0.27,0.7);
\coordinate (pivot1) at (0,1);
\coordinate (bob1) at (0.4,1);
\draw pic[draw=orange, <->,"$\omega$", angle eccentricity=1.5] {angle = bob1--origo1--pivot1};
\draw (-0.3,0.45)node[left] {$\vc$}--(0,1);
\draw [<->] (-0.3,0.45)--(0,0.45);
\coordinate (origo) at (-0.3,0.45);
\coordinate (pivot) at (0,1);
\coordinate (bob) at (0,0.45);
\draw pic[draw=orange, <->,"$\omega$", angle eccentricity=1.5] {angle = bob--origo--pivot};
\node at (-0.17,0.75) {$r$};
\draw (-1,-0.1) node[left] {$x_2=-1$} --(1,-0.1);
\end{tikzpicture}
\caption{A circle with center $\vc=(-\sqrt{r^2-1},0)$, radius $r$ is pictured. The circle intersects the dashed line at two points, whose $x_1$ coordinates are the solutions to \eqref{quad}.}
\label{figpf}
\end{figure}
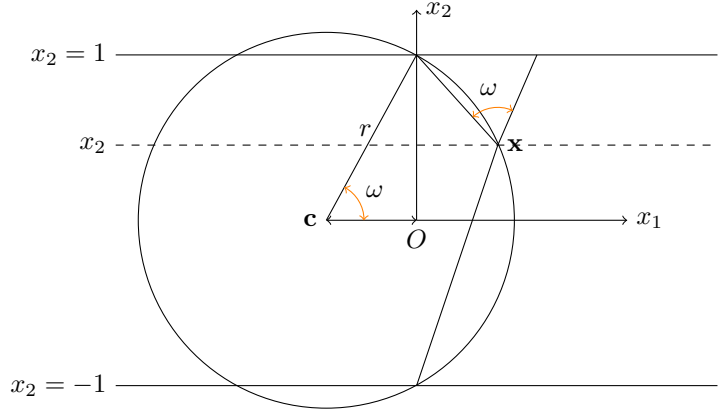
\\
Then $z=\sin\frac{\omega}{2}$ and we have
$$r\sin\omega=1\iff r=\frac{1}{2z\sqrt{1-z^2}}\iff \sqrt{r^2-1}=\frac{1-2z^2}{2z\sqrt{1-z^2}},$$
and the quadratic
\begin{equation}
\label{quad}
x_2^2+\paren{x_1+\sqrt{r^2-1}}^2=r^2\iff x_1^2+\frac{1-2z^2}{z\sqrt{1-z^2}}x_1+(x_2^2-1)=0
\end{equation}
to solve for $x_1$ for each $x_2\in\mathcal{I}$. Solving \eqref{quad} for positive $x_1=g(z)$ yields
\begin{equation}
\label{gBST}
\boxed{g(z)=\frac{\sqrt{1-4x_2^2z^2(1-z^2)}-(1-2z^2)}{2z\sqrt{1-z^2}}.}
\end{equation}
After some calculations we can derive the first order derivative of $g$,
\begin{equation}
\label{g'=gh}
\boxed{g'(z)=\frac{\mathrm{d}}{\mathrm{d}z}g(z)=g(z)h(z),\ \text{where}\ h(z)=\frac{(1-4x_2^2z^2(1-z^2))^{-1/2}}{z(1-z^2)}.}
\end{equation}
See appendix \ref{appA1} for more detail. 

Making the substitution $q=Eq_1(x_1)$ in equation \eqref{bf1} yields
\begin{equation}
\label{bgK}
\begin{split}
\frac{E}{2W_1(E)}\widehat{\mathfrak{B}f}(E,\eta_1,\Phi(x_2))&=\int_{E_m}^{c_2E}{x_1} h\left(\frac{q}{E}\right)W_2\paren{x_1,x_2}\cos\left(\eta_1 {x_1}\right)\\
&\hspace{4cm}\times\hat{f}(q,\eta_1,x_2)\mid_{x_1=g\left(\frac{q}{E}\right)}\mathrm{d}q\\
&=\int_{E_m}^{c_2E}K(q,E)\hat{f}(q,\eta_1,x_2)\mathrm{d}q,
\end{split}
\end{equation}
where $0<c_2=q_1\paren{w(x_2)}<1$. Note that $f=0$ for $q<E_m$, and hence we can set the lower integral limits in \eqref{bgK} to $E_m$. Further note that we can divide by $W_1$ in the first step by assumption that $W_1$ is non-vanishing. 

We define $\widehat{\mathcal{B}f}$ as
\begin{equation}
\begin{split}
\widehat{\mathcal{B}f}(E,\eta_1,\Phi(x_2))&=\frac{\mathrm{d}}{\mathrm{d}E}\left(\frac{E}{2W_1(E)}\widehat{\mathfrak{B}f}(E,\eta_1,\Phi(x_2))\right)\\
&=\int_{E_m}^{c_2E}K_1(q,E)\hat{f}(q,\eta_1)\mathrm{d}q+c(\eta_1)\hat{f}(c_2E,\eta_1,x_2),
\end{split}
\end{equation}
where $K_1(q,E)=\frac{\mathrm{d}}{\mathrm{d}E}K(q,E)$. Here
\begin{equation}
\label{ceta1}
c(\eta_1)=c_2w(x_2)h(c_2)W_2\paren{w(x_2),x_2}\cos\paren{\eta_1 w(x_2)}
\end{equation}
is non-zero for almost all $\eta_1\in\mathbb{R}$ by assumption that $W_2$ is non-vanishing. Let $\eta_1\in\mathbb{R}\backslash\left\{\eta_1=\frac{\pi/2+j\pi}{w(x_2)} : j\in\mathbb{Z}\right\}$ be such that $c(\eta_1)\neq 0$. Then
\begin{equation}
\label{csD}
\begin{split}
\frac{1}{c(\eta_1)}\widehat{\mathcal{B}f}\left(\frac{E}{c_2},\eta_1,\Phi(x_2)\right)&=\frac{1}{c(\eta_1)}\int_{E_m}^{E}K_1\left(q,\frac{E}{c_2}\right)\hat{f}(q,\eta_1,x_2)\mathrm{d}q+\hat{f}(E,\eta_1,x_2)\\
&=\frac{1}{c(\eta_1)}\int_{E_m}^{E}K_2\left(q,E\right)\hat{f}(q,\eta_1,x_2)\mathrm{d}q+\hat{f}(E,\eta_1,x_2).
\end{split}
\end{equation}
We will now show that the kernel $K_2(q,E)$ is bounded on the triangle $T=\{E_m<E<E_M, E_m< q<E\}$, thus proving {the injectivity} of the weighted Bragg transform by Theorem \ref{volt}. We can write
$$K(q,E)=h\paren{\frac{q}{E}}G\paren{g\paren{\frac{q}{E}}},$$
where
$$G(x_1)=x_1 W_2(x_1,x_2)\cos(\eta_1x_1).$$
After {differentiating $K$ with respect to $E$ through application of the chain and product rules, and making the substitution $E\to \frac{E}{c_2}$ we have}
\begin{equation}
K_2(q,E)=-\frac{c^2_2q}{E^2}\left[h'\paren{\frac{c_2q}{E}}G(x_1)+g'\paren{\frac{c_2q}{E}}h\paren{\frac{c_2q}{E}}G'(x_1)\right]_{x_1=g\paren{\frac{c_2q}{E}}},
\end{equation}
where
$$G'(x_1)=W_2(x_1,x_2)\cos(\eta_1x_1)+x_1\paren{W'_2(x_1,x_2)\cos(\eta_1x_1)-W_2(x_1,x_2)\sin(\eta_1x_1)},$$
where $W'_2$ denotes the first order derivative of $W_2$ in $x_1$. By assumption, $W_2(\cdot,x_2)$ and $W'_2(\cdot,x_2)$ are bounded on $\left[-w(x_2),w(x_2)\right]$ for any $x_2\in\mathcal{I}$, and we have $\left|\frac{c^2_2q}{E^2}\right|<\frac{c^2_2E_M}{E_m^2}$ on $T$. Hence it suffices to show that $h'\left(\frac{c_2q}{E}\right)$, $g'\left(\frac{c_2q}{E}\right)$, $h\left(\frac{c_2q}{E}\right)$ and $g\left(\frac{c_2q}{E}\right)$ are bounded on $T$ to show the boundedness of $K_2$ on $T$. 

First we have the derivative of $h$
\begin{equation}
\label{h'}
h'(z)=\frac{20x_2^2z^6-28x_2^2z^4+(8x_2^2+3)z^2-1}{z^2(1-z^2)^2(1-4x_2^2z^2(1-z^2))^{3/2}}.
\end{equation}
Let $c_1=\frac{c_2E_m}{E_M}>0$. Then $c_1<\frac{c_2q}{E}<c_2$ for $(q,E)\in T$ and we have
$$\left|g\left(\frac{c_2q}{E}\right)\right|<g(c_2)=w(x_2),\ \ \ \ \forall(q,E)\in T.$$
This follows since $g$ is monotone increasing (it is the inverse of a monotone increasing function, namely $q_1$). Noting that $\sqrt{1-4x^2_2z^2(1-z^2)}\geq \sqrt{1-x^2_2}$ for $z\in(0,1)$, we can derive upper bounds for $h$,
$$\left|h\left(\frac{c_2q}{E}\right)\right|<\frac{1}{c_1(1-c_2^2)\sqrt{1-x_2^2}},\ \ \ \ \forall(q,E)\in T.$$
It follows from equation \eqref{g'=gh} that
$$\left|g'\left(\frac{c_2q}{E}\right)\right|<\frac{w(x_2)}{c_1(1-c_2^2)\sqrt{1-x_2^2}},\ \ \ \ \forall(q,E)\in T,$$
and from equation \eqref{h'} we have
$$\left|h'\left(\frac{c_2q}{E}\right)\right|<\frac{20x_2^2c_2^6+28x_2^2c_2^4+(8x_2^2+3)c_2^2+1}{c^2_1(1-c_2^2)^2(1-x_2^2)^{3/2}},\ \ \ \ \forall(q,E)\in T.$$
This completes the proof. 
\end{proof}
\end{theorem}

\begin{remark}
\label{rem1}
In Theorem \ref{main1}, $\hat{f}$ is determined almost everywhere, namely for all
$$\eta_1\in\mathbb{R}\backslash\left\{\eta_1=\frac{\pi/2+j\pi}{w(x_2)} : j\in\mathbb{Z}\right\}.$$
While this is sufficient to prove injectivity, the effects on stability are not discussed. We can expect the recovery of $\hat{f}$ to be more problematic with noise for $\eta_1\to \frac{\pi/2+j\pi}{w(x_2)}$, as in this case $\frac{1}{c(\eta_1)}\to\infty$ (see equation \eqref{ceta1}) and any noise in the data would be amplified upon multiplication by $\frac{1}{c(\eta_1)}$ in equation \eqref{csD}. For smaller source widths ($w$), the frequency of the cosine term $\cos(\eta_1 w(x_2))$ decreases, and the problem $\eta_1$ areas are more sparse. This implies a more stable inversion. Conversely for larger $w$, the Bragg inversion is less stable. {Hence we expect a difference in problem stability with varying $x_2$, based on the relationship $w(x_2)=(1+x_2)\tan\frac{\beta}{2}$ (as in equation \eqref{wdef}). As $x_2$ gets closer to the source (i.e. as $x_2\to-1$) $w$ decreases and the problem stability increases, and we expect to see the converse effect for $x_2$ closer to the detector array (i.e. as $x_2\to 1$).}
\end{remark}

\begin{corollary}
\label{cor1}
Neglecting attenuative effects, the physical modelling terms of equation \eqref{equBG} constitute a weighting of type $W$ as in Theorem \ref{main1}. 
\begin{proof}
To calculate expressions in terms of $x_1$, $x_2$ and $E$ for the physical modelling terms we define the coordinates $\vs=(0,-1,0)$, $\vd=(0,1,\Phi(x_2))$, $\vv=(0,1,0)$ and $\vx=(x_1,x_2,0)$, which correspond to the $s_1=d_1=0$ case in the notation of section \ref{sec3}. 

The initial source intensity can be modeled by the inverse square law of the source-scatterer distance \cite[page 50]{model}
$$I_0(E,x_1,x_2)=\frac{I_0(E)}{|\vs-\vx|^2}=\frac{I_0(E)}{x_1^2+(x_2+1)^2},$$
where $I_0$ is the initial source spectrum, which we can assume to be $C^1((0,\infty))$ and non-zero. The assumption of $C^1$ would mean avoiding effects such as Bremsstrahlung peaks in tube spectra, but we do not wish to go into such physical concerns here. 

The solid angle is
$$\mathrm{d}\Omega(x_1,x_2)=D_A\times\frac{\left|(\vd-\vx)\cdot\vv\right|}{|\vd-\vx|^3}=\frac{|1-x_2|}{(x_1^2+(1-x_2)^2+\Phi^2(x_2))^{\frac{3}{2}}},$$
where we have set $D_A=1$ in the last step without loss of generality. The polarization factor is
\begin{equation}
\begin{split}
P(x_1,x_2)&=\frac{1+\cos^2 2\theta(\vd,\vs,\vx)}{2}\\
&=\frac{1}{2}\paren{1+\paren{\frac{(\vx-\vs)\cdot(\vd-\vx)}{|(\vx-\vs)||(\vd-\vx)|}}^2}\\
&=\frac{1}{2}\paren{1+\frac{\paren{x_1^2-(x_2+1)(1-x_2)}^2}{(x_1^2+(x_2+1)^2)(x_1^2+(1-x_2)^2+\Phi^2(x_2))}}.
\end{split}
\end{equation}
See equation \eqref{bangle} for the calculation of the Bragg angle $\theta$. Hence we write $W_1(E)=I_0(E)$ and $W_2(x_1,x_2)=Q(x_1,x_2)P(x_1,x_2)$, where
$$Q(x_1,x_2)=\frac{|1-x_2|}{(x_1^2+(x_2+1)^2)(x_1^2+(1-x_2)^2+\Phi^2(x_2))^{\frac{3}{2}}}.$$
By assumption on $I_0(E)$ the weighting $W_1$ satisfies the conditions of Theorem \ref{main1}. Further it is clear that $W_2$ is symmetric about $x_1=0$ and non-vanishing. We have the first order derivatives of $P$ and $Q$ with respect to $x_1$,
$$P'(x_1,x_2)=\frac{x_1(1-x_2^2+x_1^2)\left[\Phi^2(x_2)(x^2_1+x_2^2+4x_2+3)+4(x^2_1-x^2_2+1)\right]}{(x_1^2+(x_2+1)^2)^2(x_1^2+(1-x_2)^2+\Phi^2(x_2))^2}$$
and
$$Q'(x_1,x_2)=-\frac{x_1|1-x_2|\left[\Phi^2(x_2)+5x^2_1+5x_2^2+2x_2+5\right]}{(x_1^2+(x_2+1)^2)^2(x_1^2+(1-x_2)^2+\Phi^2(x_2))^{\frac{5}{2}}}.$$
{The derivatives $P'$ and $Q'$ above were calculated using Wolfram Mathematica 10 mathematical software.} It follows that $Q(\cdot,x_2)$, $Q'(\cdot,x_2)$, $P(\cdot,x_2)$ and $P'(\cdot,x_2)$ are bounded on $[-w(x_2),w(x_2)]$ for any $x_2\in\mathcal{I}$ and hence $W_2(\cdot,x_2)$ and $W'_2(\cdot,x_2)$ are bounded on $[-w(x_2),w(x_2)]$ for any $x_2\in\mathcal{I}$ by the product rule. This completes the proof.
\end{proof}
\end{corollary}
\begin{corollary}
Under the further assumption that $W_1$ is bounded on $\mathfrak{E}$ and $\Phi'$ is bounded on $\mathcal{I}$, where $\Phi$ is as described in section \ref{sec3.1} and equation \ref{phi}, the operator $\mathfrak{B} : L^2_0(\mathfrak{E}\times\mathbb{R}\times\mathcal{I})\to L^2(\mathfrak{E}\times\mathbb{R}\times\Phi(\mathcal{I}))$ is bounded.
\begin{proof}
From Theorem \ref{main1}, equation \eqref{bgK}, we have
\begin{equation}
\label{newK}
\begin{split}
\widehat{\mathfrak{B}f}(E,\eta_1,\epsilon)&=\frac{2W_1(E)}{E}\int_{E_m}^{E}K(q,E)\chi_{\left[0,w(x_2)\right]}\paren{x_1}\hat{f}(q,\eta_1,\Phi^{-1}(\epsilon))\mid_{x_1=g\paren{\frac{q}{E}}}\mathrm{d}q,
\end{split}
\end{equation}
where $\frac{2W_1(E)}{E}K(q,E)\chi_{\left[0,w(x_2)\right]}\paren{g\paren{\frac{q}{E}}}<M$ is bounded on $T=\{E_m<E<E_M, E_m< q<E\}$ by assumption on $W_1$ and by the proof of boundedness of $K$ given in Theorem \ref{main1}. Note that we have included $\chi_{\left[0,w(x_2)\right]}$ in equation \eqref{bgK} and altered the integral limits accordingly. We do this for convenience so that the integrals are taken over $\mathfrak{E}=[E_m,E_M]$ to align with the domain and range of $\mathfrak{B}$. 

Let $\Omega_1=\mathfrak{E}\times\mathbb{R}\times\Phi(\mathcal{I})$ and let $\Omega_2=\mathfrak{E}\times\mathbb{R}\times \mathcal{I}$. From equation \eqref{newK}, it follows that
\begin{equation}
\label{equb1}
\begin{split}
\|\widehat{\mathfrak{B}f}\|^2_{L^2(\Omega_1)}&\leq M^2\int_{\Phi(\mathcal{I})}\int_{\mathbb{R}}\int_{E_m}^{E_M}\paren{\int_{E_m}^{E}|\hat{f}(q,\eta_1,\Phi^{-1}(\epsilon))|\mathrm{d}q}^2\mathrm{d}E\mathrm{d}\eta_1\mathrm{d}\epsilon\\
&\leq C\int_{\Phi(\mathcal{I})}\int_{\mathbb{R}}\int_{E_m}^{E_M}|\hat{f}(q,\eta_1,\Phi^{-1}(\epsilon))|^2\mathrm{d}q\mathrm{d}\eta_1\mathrm{d}\epsilon\ \ \ \text{(by Cauchy-Schwartz)}\\
&=C\int_{\Phi(\mathcal{I})}\int_{\mathbb{R}}\int_{E_m}^{E_M}|f(q,x_1,\Phi^{-1}(\epsilon))|^2\mathrm{d}q\mathrm{d}x_1\mathrm{d}\epsilon \ \ \ \text{(by Theorem \ref{plan})}\\
&=C\int_{-1}^1|\Phi'(x_2)|\int_{\mathbb{R}}\int_{E_m}^{E_M}|f(q,x_1,x_2)|^2\mathrm{d}q\mathrm{d}x_1\mathrm{d}x_2\\
&\leq M_1C\|f\|^2_{L^2(\Omega_2)},
\end{split}
\end{equation}
where $C=M^2(E_M-E_m)^2$, $|\Phi'|<M_1$ and we note that $f(q, \cdot, x_2)\in L^2(\mathbb{R})$ for almost all $(q,x_2)\in \mathbb{R}_+\times\mathbb{R}$ by Fubini's theorem (and by assumption that $f$ is $L^2$), so we can apply Plancherel's theorem to the partial Fourier transform in the third step of \eqref{equb1}. The above proves that $\widehat{\mathfrak{B}f}\in L^2(\mathfrak{E}\times\mathbb{R}\times\Phi(\mathcal{I}))$ and hence $\widehat{\mathfrak{B}f}(E,\cdot,\epsilon)\in L^2(\mathbb{R})$ for almost all $(E,\epsilon)\in \mathbb{R}_+\times\Phi(\mathcal{I})$ (by Fubini's theorem). It follows that Plancherel's theorem applies to the partial Fourier transform (in the variable $s_1$) and we have
$$\|\mathfrak{B}f\|^2_{L^2(\mathfrak{E}\times\mathbb{R}\times\Phi(\mathcal{I}))}=\|\widehat{\mathfrak{B}f}\|^2_{L^2(\mathfrak{E}\times\mathbb{R}\times\Phi(\mathcal{I}))},$$
which completes the proof.
\end{proof}
\end{corollary}

\section{Extension to $\mathbb{R}^{n+1}$ and to the surfaces of revolution of $C^2$ curves}
\label{ndim}
Here we present a generalization of our results to Radon transforms which describe the integrals over the $n$-dimensional surfaces of revolution of $C^2$ curves embedded in $\mathbb{R}^{n+1}$, for $n\geq 1$. When $n=1$, the integrals are taken over $q_1$ and its reflection about $x_1=0$ (as in figure \ref{Bgcurve}), to be clear on what we mean by a surface of revolution in $\mathbb{R}^2$. 

\subsection{The generalized Bragg transform}
\begin{figure}[!htb]
\centering
\begin{tikzpicture}[scale=3.5]
\begin{scope}[blend group=soft light]
\draw (-1,0,-1) -- (1,0,-1);
\draw (-1,0,-1) -- (-1,0,1);
\draw (-1,0,1) -- (1,0,1);
\draw (1,0,-1) -- (1,0,1);

\draw (-1,1.5,-1) -- (1,1.5,-1);
\draw (-1,1.5,-1) -- (-1,1.5,1);
\draw (-1,1.5,1) -- (1,1.5,1);
\draw (1,1.5,-1) -- (1,1.5,1);

\draw [thick,red] (-1,0.75,-1) -- (1,0.75,-1);
\draw [thick,red] (-1,0.75,-1) -- (-1,0.75,1);
\draw [thick,red] (-1,0.75,1) -- (1,0.75,1);
\draw [thick,red] (1,0.75,-1) -- (1,0.75,1);

\draw (0.45-0.5,1.5,0.45-0.5)--(0.55-0.5,1.5,0.45-0.5);
\draw (0.45-0.5,1.5,0.45-0.5)--(0.45-0.5,1.5,0.55-0.5);
\draw (0.55-0.5,1.5,0.55-0.5)--(0.45-0.5,1.5,0.55-0.5);
\draw (0.55-0.5,1.5,0.55-0.5)--(0.55-0.5,1.5,0.45-0.5);

\coordinate (O) at (0.5,0,0.5);

    \def\rx{0.6}
    \def\ry{0.2}
    \def\z{0.75}

    \path [name path = ellipse]    (0,\z) ellipse ({\rx} and {\ry});
    \path [name path = horizontal] (-\rx,\z-\ry*\ry/\z)
                                -- (\rx,\z-\ry*\ry/\z);
    \path [name intersections = {of = ellipse and horizontal}];

    \draw[fill = gray!75, gray!75] (intersection-1) -- (0,0)
      -- (intersection-2) -- cycle;
    \draw[fill = gray!50, densely dashed] (0,\z) ellipse ({\rx} and {\ry});

\path (0,0,0) coordinate (S);
\path (0,1.5,0) coordinate (D);
\path (0,0.75,0) coordinate (D1);
\path (0.6,0.75,0) coordinate (w);
\path (-0.35,0.75,-0.05) coordinate (w1);

\draw (S)--(D1);
\draw [snake it, ->](S)node[left]{\textcolor{black}{$\hat{\vs}$}}--(w1);
\draw [snake it, ->](w1)--(D);

  \end{scope}

  \draw [<-] (0.3,0.4) -- (0.775,0.2)node[below]{cone beam};

\path (0,0,0) coordinate (S);
\path (0,1.5,0) coordinate (D);
\path (0,0.75,0) coordinate (D1);
\path (0.6,0.75,0) coordinate (w);
\path (-0.35,0.75,-0.05) coordinate (w1);
\path (-0.49,1.05,-0.07) coordinate (w2);

\draw [dashed] (w1)--(w2);
\draw pic[draw=blue,thick, <->,"$\omega$", angle eccentricity=1.7] {angle = D--w1--w2};

\node at (-0.4,0.75,-0.05) {$\hat{\vx}$};
\draw [thick,<->] (D1)node[left]{$\textbf{b}$}--(w);
\node at (0.3,0.8,0) {$w$};
\draw pic[draw=blue,thick, <->,"$\beta$", angle eccentricity=1.5] {angle = w--S--D1};
\draw [thick,->] (S)--(0.3,0,0)node[below]{$x_1$};
\draw [thick,->] (S)--(0,0,0.3)node[left]{$x_2$};
\draw [thick,->] (S)--(0,1.7,0)node[left]{$x_3$};
\draw [->] (1.2,0,1.2)node[right]{$\hat{\vs}$-plane ($\{x_3=-1\}$)}--(1,0,1);
\draw [->] (1.1,0.75,-1.1)node[right]{crystal plane {$f$}}--(1,0.75,-1);
\draw [->] (-1.2,1.5,-1.2)node[left]{$\hat{\vd}$-plane ($\{x_3=1\}$)}--(-1,1.5,-1);
\draw [->] (0.65-0.5,1.5,0.35-0.5)node[right]{$\hat{\vd}$}--(0.55-0.5,1.5,0.45-0.5);

\end{tikzpicture}
\caption{Bragg scanning modality in the $n=2$ case. A square detector $\hat{\vd}=(\vs,1)$ is shown opposite a source $\hat{\vs}=(\vs,-1)$, and collects photons (shown as wavy lines) scattered from points $\hat{\vx}=(\vx,1-|\textbf{b}-\hat{\vs}|)$ on the crystal plane. The crystal sample (the red plane) is placed between, and is parallel to the $\hat{\vs}$-plane and $\hat{\vd}$-plane. The center of the base of the cone is $\textbf{b}$, the source opening angle is $\beta$ (as in figure \ref{fig1}), and the source width is $w=|\textbf{b}-\hat{\vs}|\tan\beta$. The momentum transfer is $q=Eq_1(|\vx-\vs|)=E\sin\frac{\omega}{2}$.}
\label{fig3D}
\end{figure}
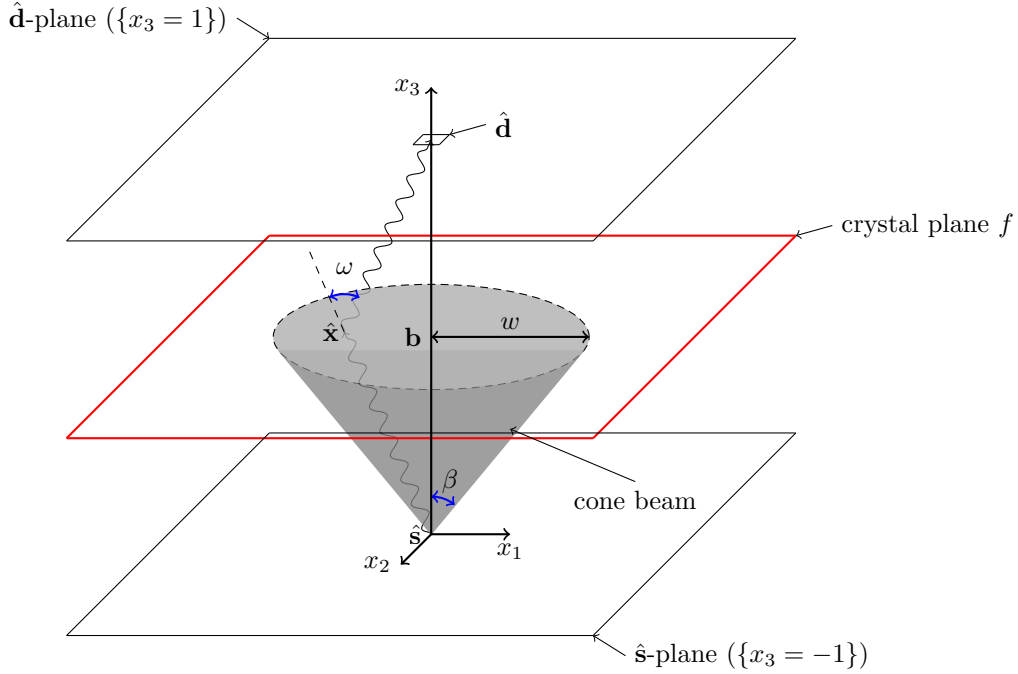

We define the generalized Bragg transform $\mathfrak{B}_n : L^2_0(\mathfrak{E}\times\mathbb{R}^{n})\to L^2(\mathbb{R}_+\times\mathbb{R}^n)$ as
\begin{equation}
\label{bg2}
\begin{split}
\mathfrak{B}_nf(E,\vs)&=\int_{\mathbb{R}^{n}}W(E,|\vx-\vs|)\chi_{B^n_w(\vs)}(\vx)f(Eq_1(|\vx-\vs|),\vx)\mathrm{d}\vx,\\
\end{split}
\end{equation}
where $B^n_w(\vs)=\{\vx\in\mathbb{R}^n : |\vx-\vs|\leq w\}${, and $w$, in this section, remains constant (i.e. $w$ does not depend on $x_2$ in this section).} To clarify, the analysis presented here differs from that of Theorem \ref{main1} in the sense that we consider general $q_1\in C^2([0,\infty))$. Theorem \ref{main1} covers the specific case in BST when $q_1$ is described by equation \eqref{bcur} and $n=1$. The $n=1,2$ cases are those to which we can ascribe physical meaning at this stage. When $n=2$, the scattering is restricted to planes of crystallites in $\mathbb{R}^3$. This is analogous to the $n=1$ case in figure \ref{figmain}, where the scatter is constrained to lines in $\mathbb{R}^2$. The $n>2$ case is of interest theoretically and adds to the works of \cite{gen1,gen2,gen3,gen4}.  For an illustration of the BST scanning geometry when $n=2$, see figure \ref{fig3D}. In this case $\hat{\vs}$ and $\hat{\vd}$ are translated opposite one another on parallel planes. The source is cone-beam and images a plane sample of crystals (at $\{x_3=1-|\textbf{b}-\hat{\vs}|\}$). The source cone and crystal plane intersection is $\{(\vx,1-|\textbf{b}-\hat{\vs}|) : \vx\in B^n_w(\vs)$\} (a disc in $\mathbb{R}^3$). After setting $q_1$ as in equation \eqref{bcur} (replacing $x_2$ with $x_2\to1-|\textbf{b}-\hat{\vs}|$), and $n=2$, $\mathfrak{B}_nf$ models the Bragg signal in the geometry of figure \ref{fig3D}. {In the $n=1$ case (as depicted in figure \ref{figmain}) the scatter is restricted to planes (or lines in 2-D) using the Venetian blind collimation technology, which is possible since the detectors can be offset in 3-D (i.e. in the $x_3$ direction). In the $n=2$ case of figure \ref{fig3D} we cannot offset the detectors in the same way (i.e. in the $x_4$ direction), and thus we instead restrict the scanning target ($f$) itself to a plane (e.g. the red plane of figure \ref{fig3D}) in order to restrict the scatter to planes and for the model \ref{bg2} to apply. In practice $f$ would represent a thin film of crystalline material placed parallel to the source and detector planes, as in figure \ref{fig3D}. Thus} we anticipate that the scanning modality of figure \ref{fig3D} will have practical application in materials characterization of thin films. {In \cite[page 150]{ballantine1996acoustic} the authors discuss some of the current methods in materials characterization of thin films, none of which fall into the framework considered in this paper.}

Refer back to the discussions in the introduction in paragraph 5. We now show that the inversion of $\mathfrak{B}_n$ is not covered by the theorems of Palamodov \cite{gen1}. The generating function $\phi$ (using the notation of \cite{gen1}) for the surfaces of integration considered here is
\begin{equation}
\label{genf}
\phi(E,\vs,q,\vx)=q-Eq_1(|\vx-\vs|).
\end{equation}
Then the surfaces of revolution are the $n$-dimensional submanifolds of $\mathcal{Z}=\{(E,\vs,q,\vx) : \phi(E,\vs,q,\vx)=0\}$,
for fixed $(E,\vs)\in\mathbb{R}_+\times \mathbb{R}^n$. For the theory of \cite{gen1} to apply we require $\phi$ to be a regular generating function. See \cite[page 2]{gen1} for the definition of a regular generating function. The $\phi$ of equation \eqref{genf} does not satisfy regularity since there are conjugate points. For example, take $E=1$, $\vs=\textbf{0}_n$ and $\vx_1,\vx_2\in S^{n-1}$ with $\vx_1\neq\vx_2$, and let $q=q_1(|\vx_1|)=q_1(|\vx_2|)$, where $\textbf{0}_n$ is the zero vector length $n$. Then $\phi(E,\vs,q,\vx_1)=\phi(E,\vs,q,\vx_2)=0$ which violates the conjugacy conditions of \cite{gen1}. Hence the reconstruction formulae of \cite{gen1} do not apply to $\mathfrak{B}_n$.

Now we have our second main theorem which proves the injectivity of $\mathfrak{B}_n$ for $n\geq 2$. {Note that in this section $w$ does not depend on $x_2$ and remains constant in the proofs presented.}
\begin{theorem}
\label{main2}
Let $n\geq 2$, $w>0$ and let $q_1\in C^2([0,\infty))$ be strictly monotone increasing with $g=q^{-1}_1$. Let $q_1$ satisfy $q_1(0)=0$ and $q'_1>0$ on $[g(c_1),w]$, where $c_1=\frac{c_2E_m}{E_M}$ and $c_2=q_1(w)$. Let $W:\mathbb{R}_+\times\mathbb{R}\to\mathbb{R}$ be a separable weighting of the form $W(E,t)=W_1(E)W_2(t)$, where $W_1\in C^1(\mathbb{R}_+)$ and $W_2\in C^1(\mathbb{R})$ are non-vanishing. Further let $W_2$ be bounded and have bounded first order derivative on $[-w,w]$. {Then $\mathfrak{B}_n$ is injective.}
\begin{proof}
Let $\vx=t\Theta$, where $t=|\vx|\in\mathbb{R}_+$ and $\Theta\in S^{n-1}$ is the direction of $\vx$. Let $\mathrm{d}\Omega_{n-1}$ denote the surface measure on $S^{n-1}$. Then equation \eqref{bg2} becomes
\begin{equation}
\label{bg2.1}
\begin{split}
\mathfrak{B}_nf(E,\vs)&=\int_{S^{n-1}}\int_0^{\infty}t^{n-1}W(E,t)\chi_{[0,w]}(t)f(Eq_1(t),t\Theta+\vs)\mathrm{d}t\mathrm{d}\Omega_{n-1}.
\end{split}
\end{equation}
Taking the Fourier transform of both sides of \eqref{bg2.1} in $\vs$ yields
\begin{equation}
\begin{split}
\widehat{\mathfrak{B}_nf}(E,\boldsymbol{\eta})&=W_1(E)\int_0^{\infty}\left[\int_{S^{n-1}}\mathrm{d}\Omega_{n-1}e^{it\Theta\cdot\boldsymbol{\eta}}\right]t^{n-1}W_2(t)\chi_{[0,w]}(t)\hat{f}(Eq_1(t),\boldsymbol{\eta})\mathrm{d}t\\
&=\int_0^{w}L(\boldsymbol{\eta},t)t^{n-1}W_2(t)\hat{f}(Eq_1(t),\boldsymbol{\eta})\mathrm{d}t,
\end{split}
\end{equation}
where $\boldsymbol{\eta}$ is dual to $\vs$ and
\begin{equation}
\label{K1}
\boxed{\begin{split}
L(\boldsymbol{\eta},t)&=\int_{S^{n-1}}\cos(t\Theta\cdot\boldsymbol{\eta})\mathrm{d}\Omega_{n-1}\\
&=V_{n-2}\int_{0}^{\pi}\sin^{n-2}(\varphi)\cos(t|\boldsymbol{\eta}|\cos\varphi)\mathrm{d}\varphi,\ \ \ \text{(by spherical symmetry)},
\end{split}}
\end{equation}
where $V_{n-2}=\frac{2\pi^{\frac{n-1}{2}}}{\Gamma(\frac{n-1}{2})}$
is the volume of $S^{n-2}$ for $n\geq 2$. The second step of \eqref{K1} follows after setting $\boldsymbol{\eta}=(\textbf{0}_{n-1},|\boldsymbol{\eta}|)$ in step one, where \textbf{0}$_{n-1}$ is the zero vector length $n-1$. Note that we can do this without loss of generality as we are taking the integral over the whole sphere in the first step of \eqref{K1}, and hence $L(\cdot,t)$ is radial for any $t\in\mathbb{R}_+$.

Making the substitution $q=Eq_1(t)$ yields
\begin{equation}
\label{VoltRn}
\begin{split}
\frac{E}{W_1(E)}\widehat{\mathfrak{B}_nf}(E,\boldsymbol{\eta})&=\int_{E_m}^{c_2E}\frac{L(\boldsymbol{\eta},t)}{q_1'({t})}{t}^{n-1}W_2\paren{{t}}\hat{f}(q,\boldsymbol{\eta})\mid_{t=g\left(\frac{q}{E}\right)}\mathrm{d}q\\
&=\int_{E_m}^{c_2E}K(q,E)\hat{f}(q,\boldsymbol{\eta})\mathrm{d}q,
\end{split}
\end{equation}
noting that we can divide by $W_1$ in the first step by assumption that $W_1$ is non-vanishing. We define $\widehat{\mathcal{B}_nf}$ as
\begin{equation}
\label{eqce}
\begin{split}
\widehat{\mathcal{B}_nf}(E,\boldsymbol{\eta})&=\frac{\mathrm{d}}{\mathrm{d}E}\paren{\frac{E}{W_1(E)}\widehat{\mathfrak{B}_nf}(E,\boldsymbol{\eta})}=\int_{E_m}^{c_2E}K_1(q,E)\hat{f}(q,\boldsymbol{\eta})\mathrm{d}q+c(\boldsymbol{\eta})\hat{f}(c_2E,\boldsymbol{\eta}),
\end{split}
\end{equation}
where $K_1(q,E)=\frac{\mathrm{d}}{\mathrm{d}E}K(q,E)$ and
$$c(\boldsymbol{\eta})=c_2\frac{L(\boldsymbol{\eta},w)}{q_1'(w)}w^{n-1}W_2\paren{w}.$$
Let $u=t\abs{\boldsymbol{\eta}}$, let $J(u)=V_{n-2}\int_{0}^{\pi}\sin^{n-2}(\varphi)\cos(u\cos\varphi)\mathrm{d}\varphi$ and let $\mathcal{E}_n=\{\boldsymbol{\eta}\in\mathbb{R}^n : J(w|\boldsymbol{\eta}|)=0\}$. We will now show that $J$ is non-zero almost everywhere, thus allowing us to divide through by $c(\boldsymbol{\eta})$ in equation \eqref{eqce} for $\boldsymbol{\eta}\in\mathbb{R}^n\backslash \mathcal{E}_n$. Using the Maclaurin series of $\cos u$, we have
\begin{equation}
\label{J(u)}
\boxed{J(u) =  \sum_{m=0}^\infty \frac{(-1)^m u^{2m}}{2m!}\paren{\int_0^\pi
\sin^{n-2}(\vp)\cos^{2m}(\vp)\,d\vp}.}
\end{equation}
Switching the sum and integral is justified because the integrals in
parenthesis are uniformly bounded by $\pi$, say, so for any $M>0$, the
sum converges uniformly absolutely for $\abs{u}\leq M$.  This means
that $J$ is an entire real analytic function and hence, by analytic continuation, $J=0$ on at most a set of
measure zero. In fact in the $n=2$ case $J$ reduces to a Bessel function $J_0$ of the first kind and $J(u)=2\pi J_0(u)$.

Now, by assumption that $W_2$ is non-vanishing, $c(\boldsymbol{\eta})\neq 0$ for all $\boldsymbol{\eta}\in\mathbb{R}^n\backslash \mathcal{E}_n$, and hence
\begin{equation}
\label{VoltR}
\begin{split}
\frac{1}{c(\boldsymbol{\eta})}\widehat{\mathcal{B}_nf}\paren{\frac{E}{c_2},\boldsymbol{\eta}}&=\frac{1}{c(\boldsymbol{\eta})}\int_{E_m}^{E}K_1\paren{q,\frac{E}{c_2}}\hat{f}(q,\boldsymbol{\eta})\mathrm{d}q+\hat{f}(E,\boldsymbol{\eta})\\
&=\frac{1}{c(\boldsymbol{\eta})}\int_{E_m}^{E}K_2\paren{q,E}\hat{f}(q,\boldsymbol{\eta})\mathrm{d}q+\hat{f}(E,\boldsymbol{\eta})
\end{split}
\end{equation}
for $\boldsymbol{\eta}\in\mathbb{R}^n\backslash \mathcal{E}_n$. Let $f_1(z)=L\paren{\boldsymbol{\eta},\left(g\left(z\right)\right)}$, $f_2(z)=\paren{q_1'(g\left(z\right))}^{-1}=g'(z)$ and $f_3(z)=g(z)^{n-1}$. Then
$$K(q,E)=G\paren{\frac{q}{E}}W_2\paren{g\left(\frac{q}{E}\right)},$$
where $G(z)=f_1(z)f_2(z)f_3(z)$. It follows that
\begin{equation}
K_2(q,E)=-\frac{c^2_2q}{E^2}\Big[g'(z)W'_2(g(z))G(z)+W_2(g(z))G'(z)\Big]_{z=\frac{c_2q}{E}},
\end{equation}
where
$$G'(z)=f'_1(z)f_2(z)f_3(z)+f_1(z)\paren{f'_2(z)f_3(z)+f'_3(z)f_2(z)}.$$
By assumption, $W_2$ and $W'_2$ are bounded on $[-w,w]$, and we have $\left|\frac{c^2_2q}{E^2}\right|<\frac{c^2_2E_M}{E_m^2}$ on $T=\{E_m<E<E_M, E_m< q<E\}$. Hence to show the boundedness of $K_2\paren{q,E}$ on $T$, it is sufficient to show that the $f_j(t)$ and their first order derivatives are bounded on $[c_1,c_2]$, where $c_1=\frac{c_2E_m}{E_M}$. The boundedness of $f_1$ and $f_3$ on $[c_1,c_2]$ is trivial. By assumption we know that $q'_1>0$ and continuous on $[g(c_1),w]$. Hence $q'_1$ is bounded away from zero on $[g(c_1),w]$, and $|f_2|<M$ is bounded on $[c_1,c_2]$.

The derivatives of the $f_j$ are
\begin{equation}
\begin{split}
f'_1(z)&=-f_2(z)\int_{S^{n-1}}(\Theta\cdot\boldsymbol{\eta})\sin(g(z)\Theta\cdot\boldsymbol{\eta})\mathrm{d}\Omega_{n-1},\\
\end{split}
\end{equation}
$f'_2(z)=-f_2(z)^3q_1''(g\left(z\right))$ and $f'_3(z)=f_2(z)g(z)^{n-2}$. By assumption $q''_1$ is continuous on $[g(c_1),w]$ and hence $q''_1\circ g$ is bounded on $[c_1,c_2]$. Further $|f'_1|<M\frac{\pi}{2w}V_{n-1}$, where $V_{n-1}$ is the surface area of $S^{n-1}$ and it clear that the $f'_j$ are bounded on $[c_1,c_2]$. By Theorem \ref{volt} and the convergence of the Neumann series of \eqref{VoltR}, it follows that we have a unique recovery of $\hat{f}(q,|\boldsymbol{\eta}|\xi)$ for all $q\in \mathfrak{E}$ and $\xi\in S^{n-1}$, and for almost all $|\boldsymbol{\eta}|\in\mathbb{R}_+$. Hence we can {reconstruct $f$} by inverse Fourier transform.
\end{proof}
\end{theorem}

\subsection{The special case when $n=1$}
\begin{figure}[!h]
\begin{subfigure}{0.32\textwidth}
\includegraphics[width=0.9\linewidth, height=4cm, keepaspectratio]{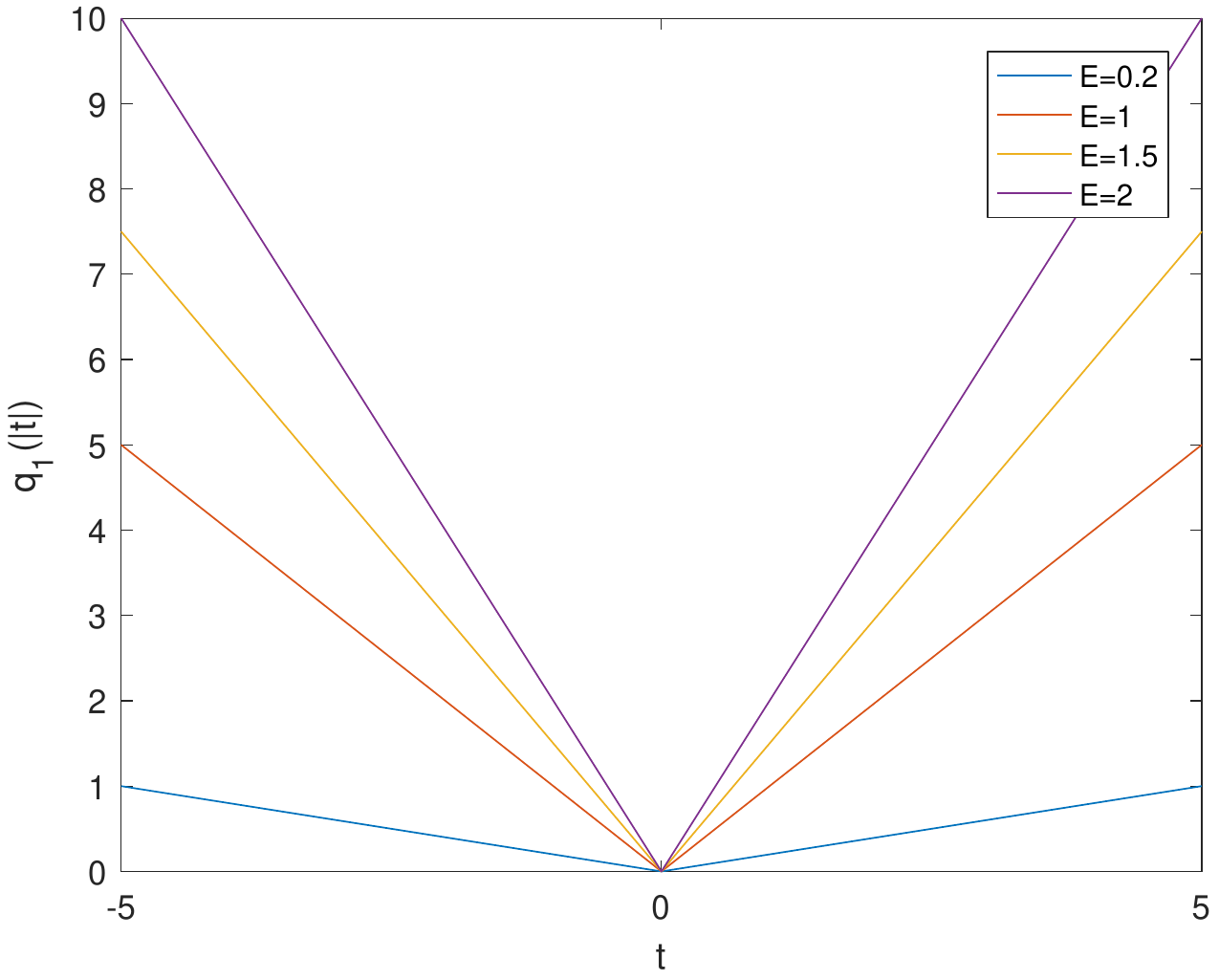}
\subcaption{$q_1(t)=t$}
\end{subfigure}
\begin{subfigure}{0.32\textwidth}
\includegraphics[width=0.9\linewidth, height=4cm, keepaspectratio]{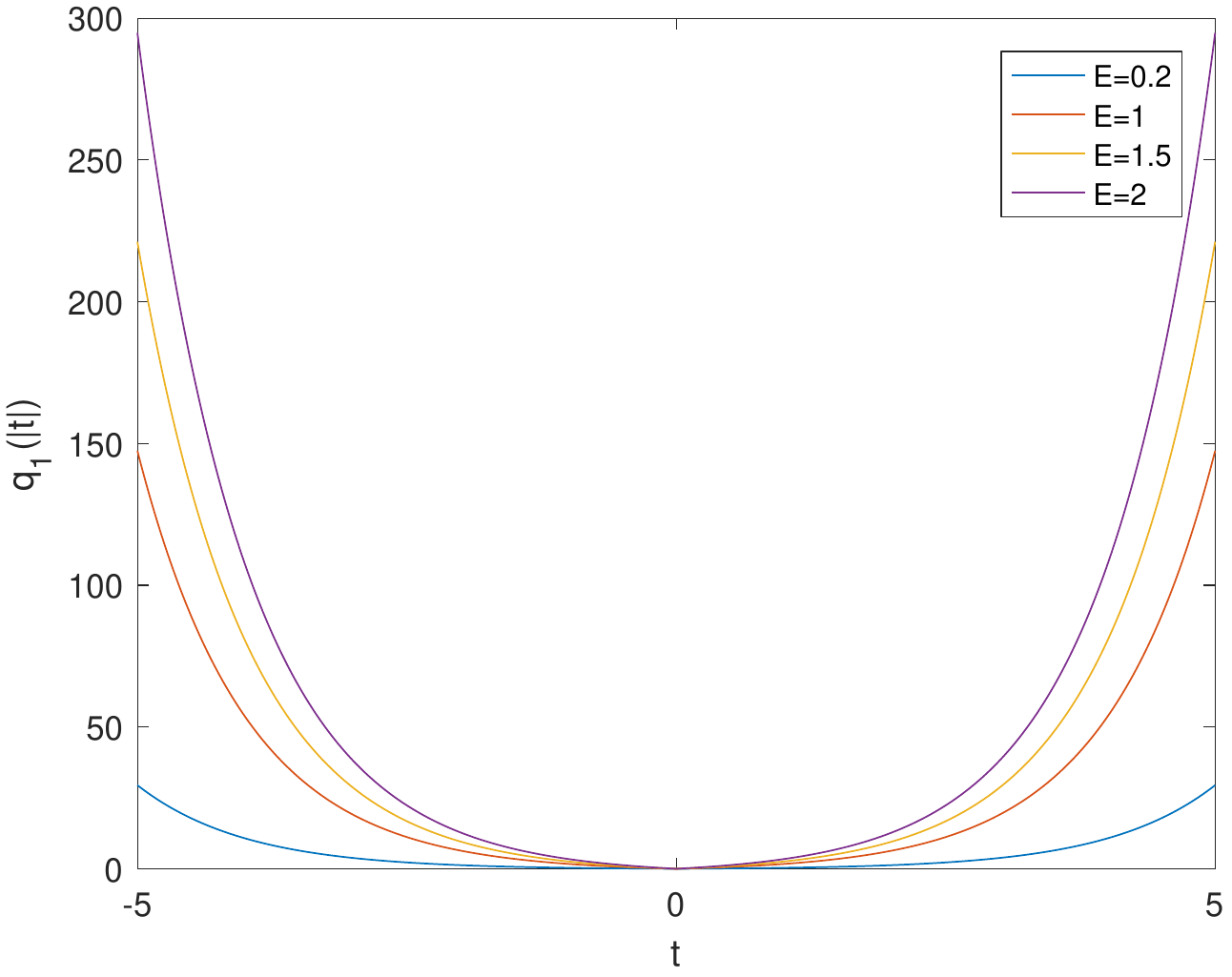} 
\subcaption{$q_1(t)=e^t-1$}
\end{subfigure}
\begin{subfigure}{0.32\textwidth}
\includegraphics[width=0.9\linewidth, height=4cm, keepaspectratio]{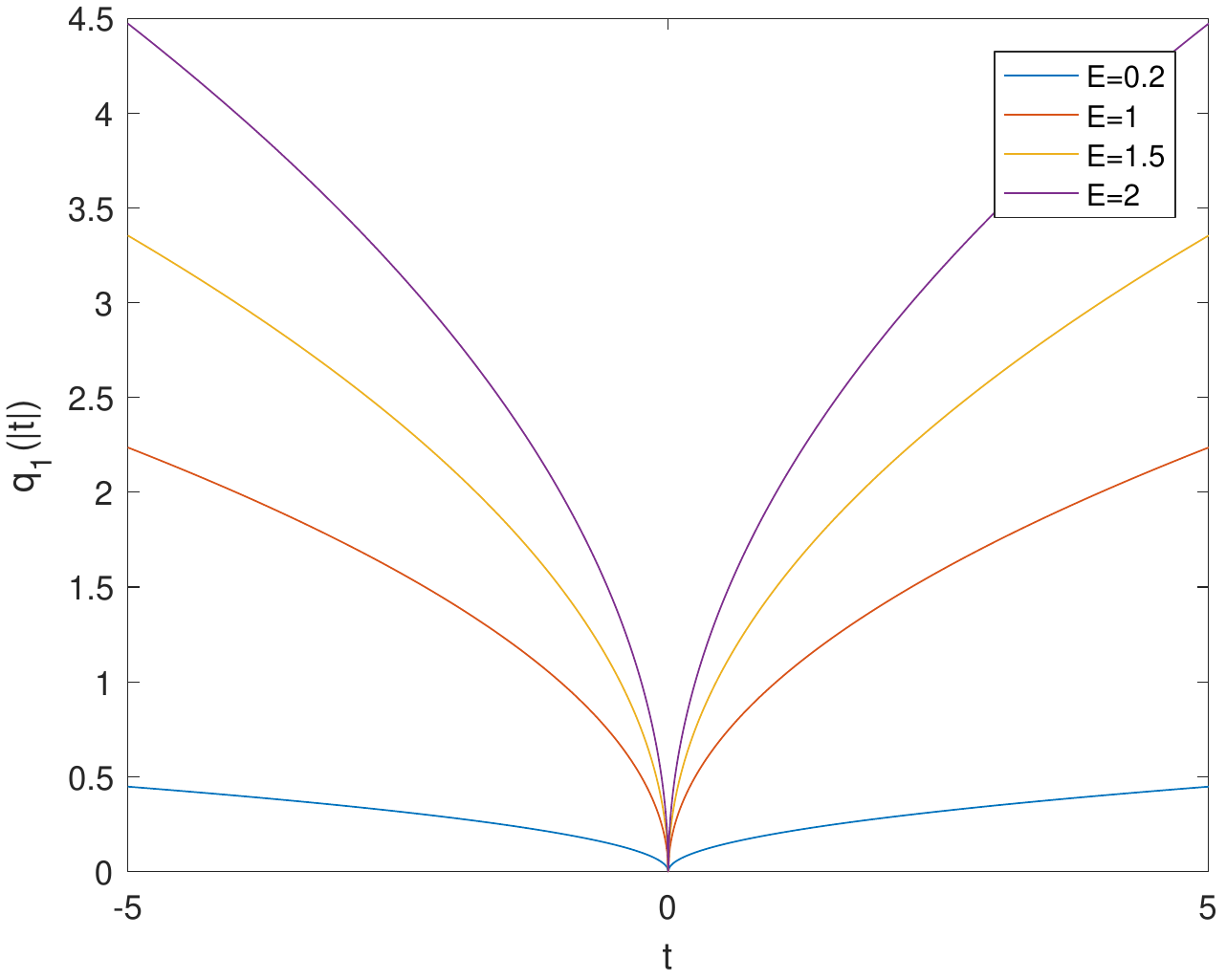}
\subcaption{$q_1(t)=\sqrt{t}$}\label{gc3}
\end{subfigure}
\begin{subfigure}{0.32\textwidth}
\includegraphics[width=0.9\linewidth, height=4cm, keepaspectratio]{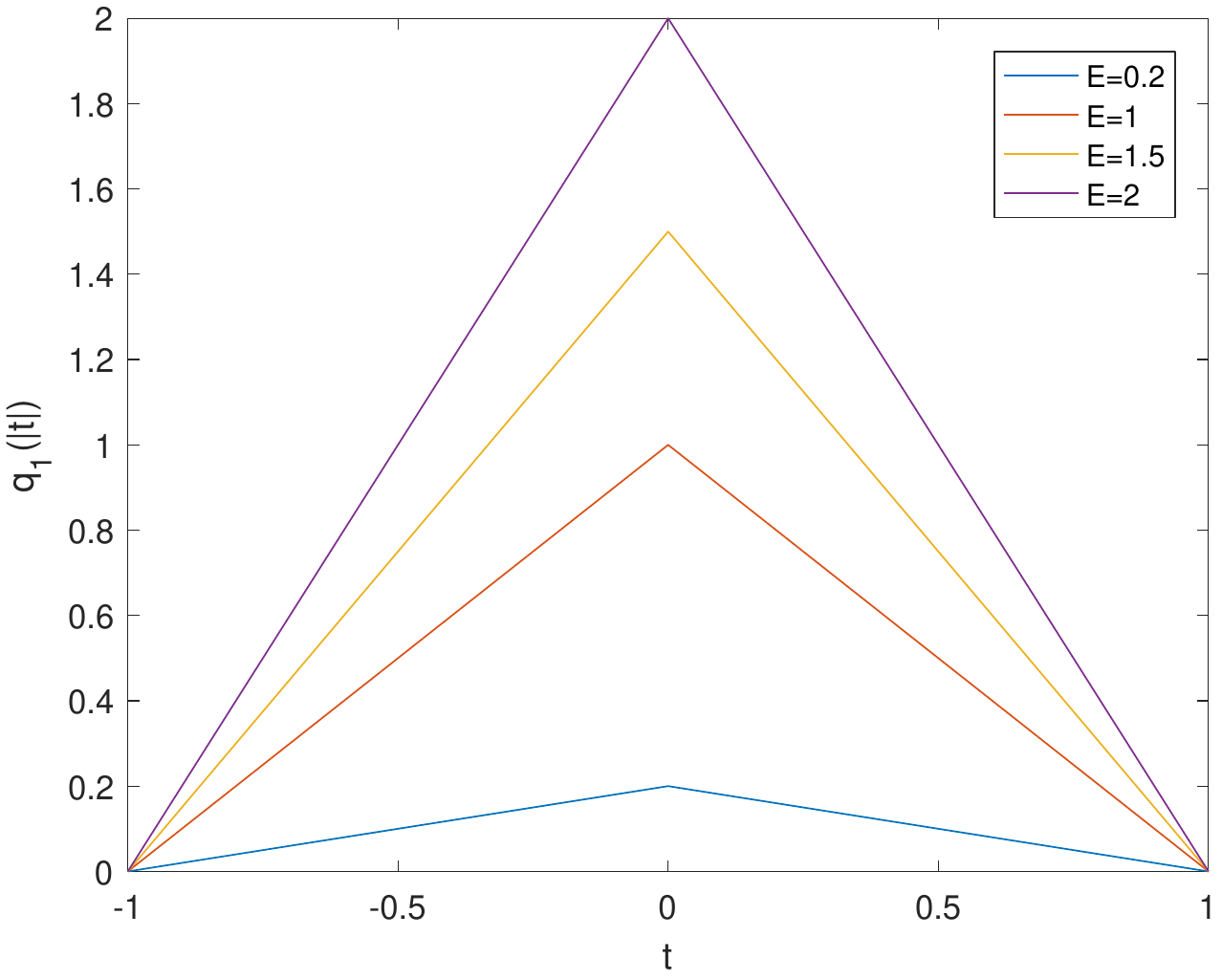}
\subcaption{$q_1(t)=1-t$}
\end{subfigure}
\begin{subfigure}{0.32\textwidth}
\includegraphics[width=0.9\linewidth, height=4cm, keepaspectratio]{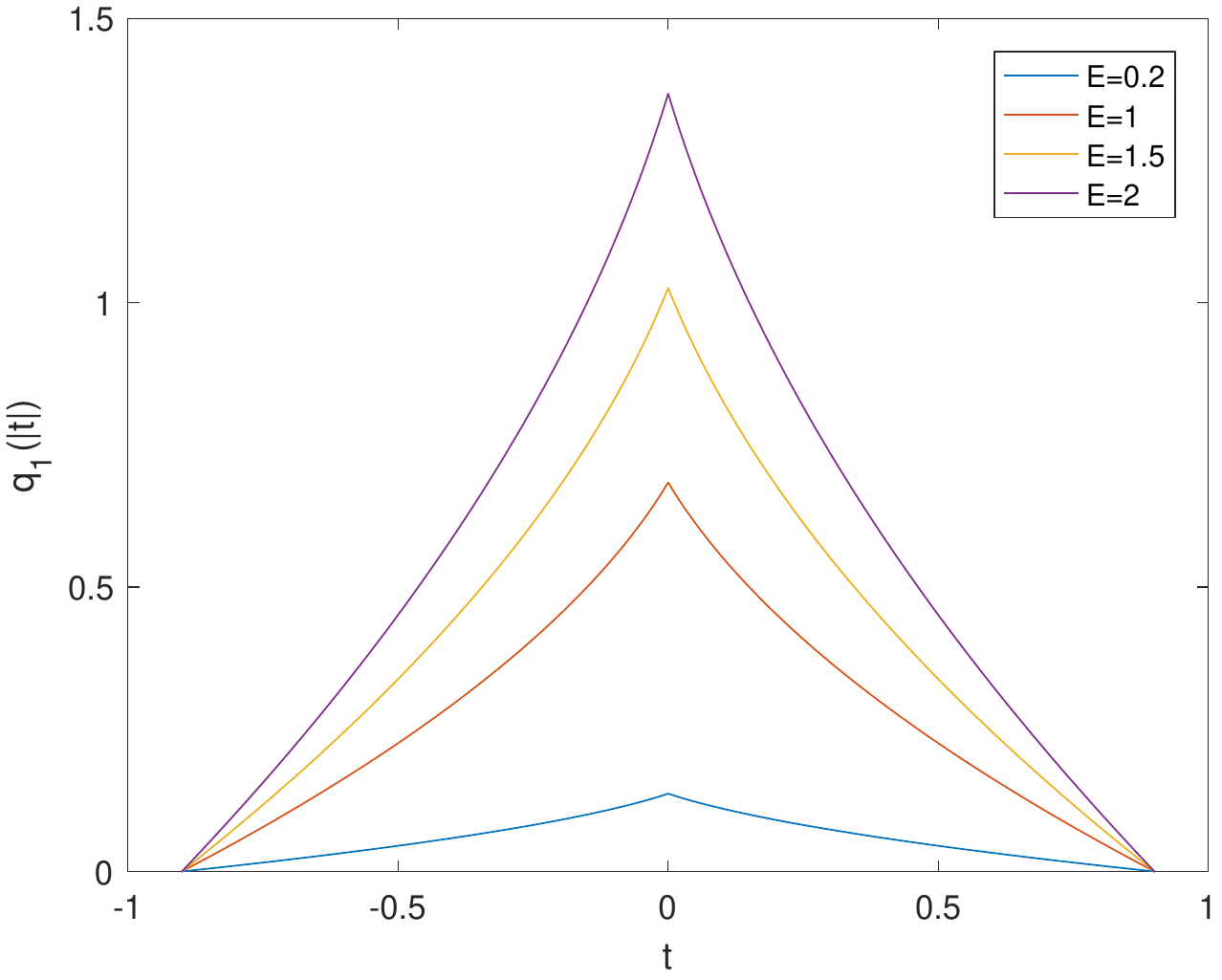} 
\subcaption{$q_1(t)=1-\sqrt{t+\frac{1}{10}}$}
\end{subfigure}
\begin{subfigure}{0.32\textwidth}
\includegraphics[width=0.9\linewidth, height=4cm, keepaspectratio]{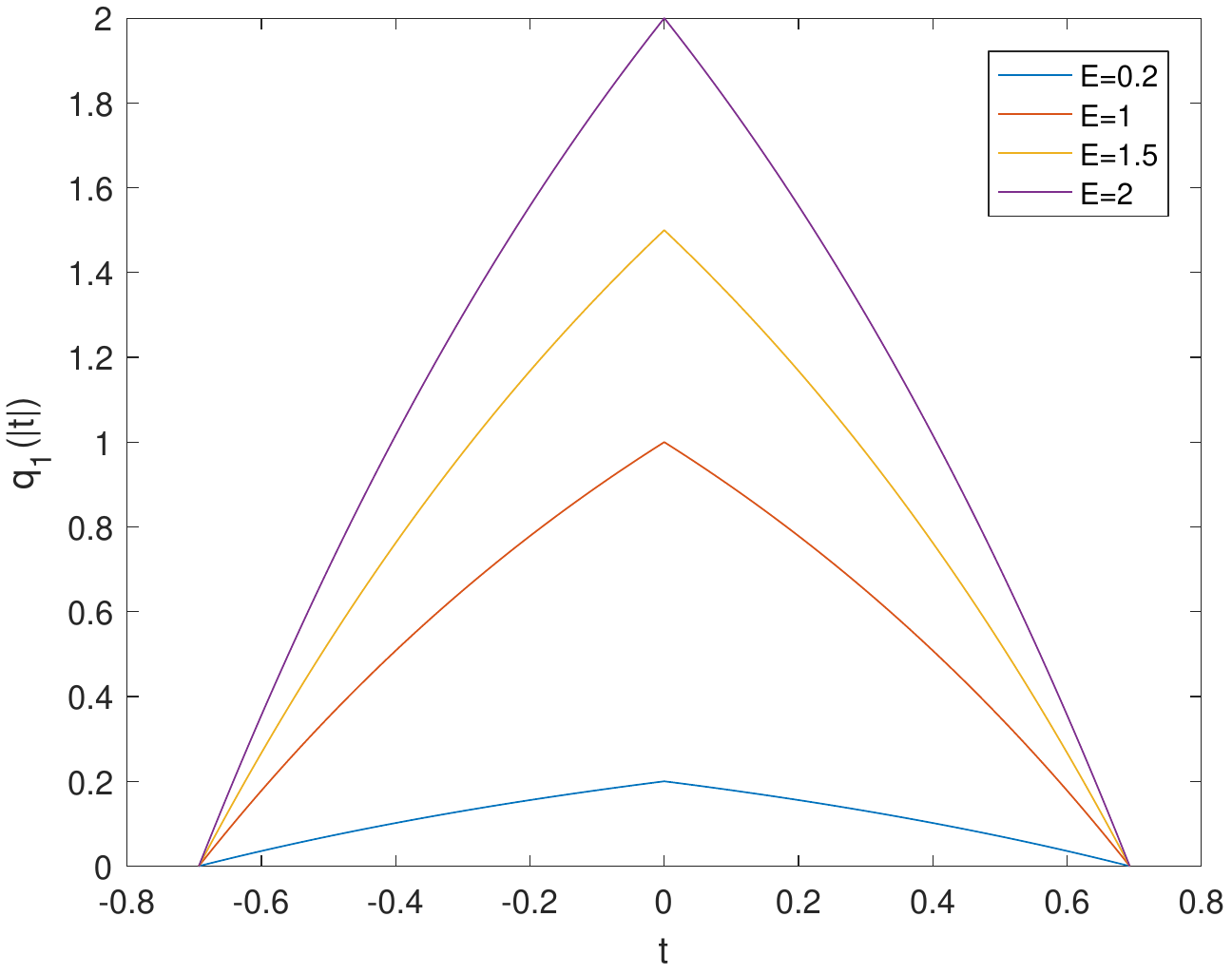}
\subcaption{$q_1(t)=2-e^t$}
\end{subfigure}
\caption{$q_1$ curve examples. For the decreasing curves displayed we would have to choose $w>1$ for the injectivity of $\mathfrak{B}_1$ to hold.}
\label{GC2}
\end{figure}

In this case we allow the curves of integration $q_1$
to be increasing or decreasing and assume they satisfy the following conditions:
\begin{enumerate}[label=(\roman*)]
\item \label{cond1} increasing $q_1$ -- $q_1(0)=0$ and $q'_1>0$ on $[g(c_1),w]$, where $c_1=\frac{c_2E_m}{E_M}$ and $c_2=q_1(w)$. 
\item \label{cond2} decreasing $q_1$ -- $q_1(0)>0$, $\lim_{t\to\infty}q_1(t)<0$, $w>g(0)$ and $q'_1<0$ on $[0,g(c_1)]$, where $c_1=\frac{c_2E_m}{E_M}$ and $c_2=q_1(0)$.
\end{enumerate}
\begin{theorem}
\label{main3}
Let $q_1\in C^2([0,\infty))$ be strictly monotone and let $g=q^{-1}_1$. Let $W:\mathbb{R}_+\times\mathbb{R}\to\mathbb{R}$ be a separable weighting with the same properties as in Theorem \ref{main2}. Then under the conditions \ref{cond1} or \ref{cond2} above for increasing or decreasing $q_1$ respectively, {$\mathfrak{B}_1$ is injective.}
\begin{proof}
Let $q_1$ be decreasing and let condition \ref{cond2} be satisfied. Then setting $n=1$ and taking the Fourier transform of \eqref{bg2} in $s_1$ gives
\begin{equation}
\label{VoltRn1}
\begin{split}
\frac{1}{2W_1(E)}\widehat{\mathfrak{B}_1f}(E,\eta_1)&=\int_0^{w}W_2(t)\cos(t\eta_1)\hat{f}(Eq_1(t),\eta_1)\mathrm{d}t\\
&=-\int_{E_m}^{c_2E}\frac{\cos\left(\eta_1 {t}\right)}{q_1'({t})}W_2\paren{{t}}\hat{f}(q,\eta_1)\mid_{t=g\left(\frac{q}{E}\right)}\mathrm{d}q\\
&=-\int_{E_m}^{c_2E}K(q,E)\hat{f}(q,\eta_1)\mathrm{d}q,
\end{split}
\end{equation}
noting that we can replace the lower limit $Eq_1(w)<0$ (strictly negative since $w>g(0)$) by $E_m$ for all $E>\frac{E_m}{c_2}$, since $f=0$ for $q<E_m$. We define $\widehat{\mathcal{B}_1f}$ as
\begin{equation}
\begin{split}
-\widehat{\mathcal{B}_1f}(E,\eta_1)&=\frac{\mathrm{d}}{\mathrm{d}E}\paren{\frac{E}{2W_1(E)}\widehat{\mathfrak{B}_1f}(E,\eta_1)}\\
&=\int_{E_m}^{c_2E}\frac{\mathrm{d}}{\mathrm{d}E}K(q,E)\hat{f}(q,\eta_1)\mathrm{d}q+c_2\frac{\cos\left(\eta_1 w\right)}{q_1'(w)}W_2(w)\hat{f}(c_2E,\eta_1)\\
&=\int_{E_m}^{c_2E}K_1(q,E)\hat{f}(q,\eta_1)\mathrm{d}q+c(\eta_1)\hat{f}(c_2E,\eta_1),
\end{split}
\end{equation}
where $c(\eta_1)=c_2\frac{1}{q_1'(0)}W_2(w)\neq 0$ since $q'(0)<0$ and $W_2$ are non-vanishing. It follows that
\begin{equation}
\begin{split}
\frac{1}{c(\eta_1)}\widehat{\mathcal{B}_1f}\paren{\frac{E}{c_2},\eta_1}&=\frac{1}{c(\eta_1)}\int_{E_m}^{E}K_1\paren{q,\frac{E}{c_2}}\hat{f}(q,\eta_1)\mathrm{d}q+\hat{f}(E,\eta_1)\\
&=\frac{1}{c(\eta_1)}\int_{E_m}^{E}K_2\paren{q,E}\hat{f}(q,\eta_1)\mathrm{d}q+\hat{f}(E,\eta_1)
\end{split}
\end{equation}
for $\eta_1\in\mathbb{R}\backslash\left\{\eta_1=\frac{\pi/2+j\pi}{w} : j\in\mathbb{Z}\right\}$. Let $f_4(t)=\cos(\eta_1 g(t))$. Then
$$K(q,E)=G\paren{\frac{q}{E}}W_2\paren{g\left(\frac{q}{E}\right)},$$
where $G(z)=f_4(z)f_2(z)f_3(z)$, and $f_2$ and $f_3$ are as in the $n\geq 2$ case. We have $f'_4(z)=-\eta_1f_2(z)\sin(\eta_1 g(z))$, and hence $f_4$ and its first order derivative are bounded on $[c_1,c_2]$, where $c_1=\frac{c_2E_m}{E_M}$ as before. By the same arguments as in the $n\geq 2$ case, it follows that $K_2\paren{q,E}$ is bounded on $T=\{E_m<E<E_M, E_m< q<E\}$. Hence we can recover the Fourier transform $\hat{f}(q,\eta_1)$ for all $q\in[0,E_M]$ and almost all $\eta_1\in\mathbb{R}$ by Theorem \ref{volt}. The injectivity of the Bragg transform follows. 

For $q_1$ increasing under condition \ref{cond1} the argument is the same as in Theorem \ref{main2} except in the calculations, $L$ of equation \eqref{K1} is replaced by $L(\eta_1,t)=\cos(\eta_1 t)$.
\end{proof}
\end{theorem}
\begin{remark}
We note that in the $n\geq 2$ case for decreasing curves we cannot prove injectivity in the same way as Theorem \ref{main3} since the kernel
$$K(c_2E,E)=\frac{L(\boldsymbol{\eta},g(c_2))}{q_1'(g(c_2))}g^{n-1}(c_2)W_2\paren{g(c_2)}\equiv 0\ \ \  (\text{noting that $g(c_2)=0$}) $$
of equation \eqref{VoltRn} is zero on $\{q=c_2E\}$, and hence violates the conditions of \cite[page 15]{Tric}. Thus at this stage we prove the injectivity of $\mathfrak{B}_n$ for decreasing curves when $n=1$. In the case of decreasing $q_1$, we see from the proof of Theorem \ref{main3} that only certain source widths (i.e. $w>g(0)$) are sufficient for a unique solution. Without such conditions on $w$, the lower limit of the integral in the third step of \eqref{VoltRn1} would vary with $E$, hence leading to an integral equation that is not of classical Volterra type. While this does not imply non-uniqueness, it does prevent us however from using the well established theory on linear Volterra equations \cite{Tric}. Hence at this stage we require the $w>g(0)$ condition to prove injectivity.  

For some example generalized Bragg curves, see figure \ref{GC2}. We note that $q_1(t)=\sqrt{t}$ (figure \ref{gc3}) fails to be $C^2$ on $[0,\infty)$ (at the origin). However since the support of $f$ is assumed to be bounded away from the origin, we can consider such curves for inversion.
\end{remark}
\begin{corollary}
Under the additional assumption that $W_1$ is bounded on $\mathfrak{E}$, the generalized Bragg transform $\mathfrak{B}_n : L^2_0(\mathfrak{E}\times\mathbb{R}^{n})\to L^2(\mathfrak{E}\times\mathbb{R}^n)$ is a bounded operator for $n\geq1$.
\begin{proof}
This follows similar ideas to the proof of Corollary \ref{cor1} and is a consequence of the boundedness of the kernel $K$ on $T$ as in equation \eqref{VoltRn1} ($n=1$ case) and equation \eqref{VoltRn} ($n\geq 2$ case). Note also that we have restricted the codomain of $\mathfrak{B}_{n}$ to $\mathfrak{E}$ to prove boundedness as in Corollary \ref{cor1}.
\end{proof}
\end{corollary}
\section{Dealing with the out-of-plane detectors}
\label{sdoff}
In section \ref{bgsec1} we provided {injectivity results} for $\mathfrak{B}$ under the $\epsilon=0$ assumption. That is, under the assumption of a negligible source-detector offset. Here we cover the $\epsilon>0$ case and modify our inversion results to address the out-of-plane detectors. 

\subsection{The offset Bragg transform}
\label{6.2}
Let $f\in L^2_0(\mathfrak{E}\times \mathbb{R}\times\mathcal{I})$ and let $\Phi$ be as described at the start of section \ref{sec3.1} (see equation \ref{phi}). Let
\begin{equation}
\label{q13}
q_1(x_1)=\frac{1}{\sqrt{2}}\sqrt{1+\frac{x_1^2-(1-x^2_2)}{\sqrt{x_1^2+(x_2+1)^2}\sqrt{x_1^2+(1-x_2)^2+\epsilon^2}}},
\end{equation}
where in this case the source-detector offset $\epsilon=\Phi(x_2)$ is included in the calculation. Then the offset Bragg transform $\mathfrak{B}_{\epsilon} : L^2_0(\mathfrak{E}\times\mathbb{R}\times \mathcal{I})\to L^2(\mathbb{R}_+\times\mathbb{R}\times\Phi(\mathcal{I}))$ is defined
\begin{equation}
\begin{split}
\label{bg4}
\mathfrak{B}_{\epsilon}f(E,s_1,\epsilon)=\int_{\mathbb{R}}\chi_{\left[-w(x_2),w(x_2)\right]}(x_1-s_1)&W\paren{E,x_1-s_1,\Phi^{-1}(\epsilon)}\\
&\times f(Eq_1(x_1-s_1),x_1,\Phi^{-1}(\epsilon))\mathrm{d}x_1.
\end{split}
\end{equation}
We now have our final main theorem which proves {the injectivity} of the offset Bragg transform.
\begin{theorem}
\label{main4}
Let $W:\mathbb{R}_+\times\mathbb{R}\times \mathcal{I}\to\mathbb{R}$ be a separable weighting of the form $W(E,x_1,x_2)=W_1(E)W_2(x_1,x_2)$, where $W_1,W_2$ have the same properties as in Theorem \ref{main1}. Further let  $c_{\epsilon}E_M<c_2E_m$ be satisfied for all pairs $(x_2,\epsilon)=(x_2,\Phi(x_2))$ for $x_2\in\mathcal{I}$, where $c_2=q_1\paren{w(x_2)}$ and $c_{\epsilon}=q_1(0)$. {Then $\mathfrak{B}_{\epsilon}$ is injective.}
\begin{proof}
Letting $g(z)=q^{-1}_1(z)$ and following the same inversion process as in Theorem \ref{main1} gives
\begin{equation}
\begin{split}
\frac{E}{2W_1(E)}\widehat{\mathfrak{B}_{\epsilon}f}(E,\eta_1,\Phi(x_2))=\int_{c_{\epsilon}E}^E\frac{\cos\left(\eta_1 {x_1}\right)}{q'_1\left({x_1}\right)}W_2\paren{x_1,x_2}&\chi_{\left[0,w(x_2)\right]}\left({x_1}\right)\\
&\times\hat{f}(q,\eta_1,x_2)\mid_{x_1=g\left(\frac{q}{E}\right)}\mathrm{d}q,
\end{split}
\end{equation}
where $q'_1$ denotes the first derivative of $q_1$. Note that in the expression for $q_1$, $\epsilon$ is now substituted for $\Phi(x_2)$. By assumption that $c_{\epsilon}<\frac{c_2E_m}{E_M}$ we have
\begin{equation}
\label{eqb}
\begin{split}
\frac{E}{2c_2W_1\paren{\frac{E}{c_2}}}\widehat{\mathfrak{B}_{\epsilon}f}\paren{\frac{E}{c_2},\eta_1,\Phi(x_2)}&=\int_{\frac{c_{\epsilon}E}{c_2}}^{\frac{E}{c_2}}\frac{\cos\left(\eta_1 {x_1}\right)}{q'_1\left({x_1}\right)}W_2\paren{x_1,x_2}\chi_{\left[0,w(x_2)\right]}\left({x_1}\right)\\
&\hspace{3.5cm}\times\hat{f}(q,\eta_1,x_2)\mid_{x_1=g\left(\frac{c_2q}{E}\right)}\mathrm{d}q\\
&=\int_{E_m}^{E}\frac{\cos\left(\eta_1 {x_1}\right)}{q'_1\left({x_1}\right)}W_2\paren{x_1,x_2}\\
&\hspace{3.5cm}\times\hat{f}(q,\eta_1,x_2)\mid_{x_1=g\left(\frac{c_2q}{E}\right)}\mathrm{d}q,
\end{split}
\end{equation}
for all $E\in [0,E_M]$ and $x_2\in \mathcal{I}$. The lower limit in the second step of \eqref{eqb} reduces to ${E_m}$ since $\frac{c_{\epsilon}E}{c_2}<\frac{c_{\epsilon}E_M}{c_2}<E_m$ for $E\in \mathfrak{E}$ and $\hat{f}=0$ for $q<E_m$ by assumption. 

After some calculations we have the first order derivative of $q_1$,
\begin{equation}
\label{box1}
\boxed{q'_1(x_1)=\frac{h(x_1)}{q_1(x_1)},\ \ \ h(x_1)=\frac{P_1(x_1)}{4h^3_1(x_1)},}
\end{equation}
where
\begin{equation}
\label{hBST1}
\boxed{h_1(x_1)=\sqrt{(x_1^2+(1+x_2)^2)(x_1^2+(1-x_2)^2+\Phi^2(x_2))},}
\end{equation}
and
\begin{equation}
\label{PBST1}
\boxed{P_1(x_1)=4x_1(1-x_2^2+x_1^2)+\Phi^2(x_2)x_1((x_2+1)(x_2+3)+x_1^2).}
\end{equation}
See appendix \ref{appA2} for more detail. Hence
\begin{equation}
\label{Volte}
\begin{split}
\frac{E}{2c_2W_1\paren{\frac{E}{c_2}}}\widehat{\mathfrak{B}_{\epsilon}f}\paren{\frac{E}{c_2},\eta_1,\Phi(x_2)}&=c_2\int_{E_m}^{E}\frac{q\cos\left(\eta_1 {x_1}\right)}{Eh\left({x_1}\right)}W_2\paren{x_1,x_2}\\
&\hspace{3.5cm}\times\hat{f}(q,\eta_1,x_2)\mid_{x_1=g\left(\frac{c_2q}{E}\right)}\mathrm{d}q\\
&=\int_{E_m}^{E}K(q,E)\hat{f}(q,\eta_1,x_2)\mathrm{d}q,
\end{split}
\end{equation}
noting that in the first step we have used
\begin{equation}
\label{box2}
\boxed{g'\paren{z}=\frac{1}{q'_1\paren{g\paren{z}}}=\frac{z}{h\paren{g\paren{z}}},}
\end{equation}
substituting $z=\frac{c_2q}{E}$. Equation \eqref{box2} follows directly from equation \eqref{box1}. We define $\widehat{\mathcal{B}_{\epsilon}f}$ as
\begin{equation}
\begin{split}
\widehat{\mathcal{B}_{\epsilon}f}(E,\eta_1,\Phi(x_2))&=\frac{\mathrm{d}}{\mathrm{d}E}\paren{\frac{E}{2c_2W_1\paren{\frac{E}{c_2}}}\widehat{\mathfrak{B}_{\epsilon}f}\paren{\frac{E}{c_2},\eta_1,\Phi(x_2)}}\\
&=\int_{E_m}^{E}K_1(q,E)\hat{f}(q,\eta_1,x_2)\mathrm{d}q+c(\eta_1)E\hat{f}(E,\eta_1,x_2),
\end{split}
\end{equation}
where $K_1(q,E)=\frac{\mathrm{d}}{\mathrm{d}E}K(q,E)$ and
{$$c(\eta_1)=c_2\frac{\cos\left(\eta_1 w(x_2)\right)}{h\left(w(x_2)\right)}W_2\paren{w(x_2),x_2}\neq 0$$}
away from $\left\{\eta_1=\frac{\pi/2+j\pi}{w(x_2)} : j\in\mathbb{Z}\right\}$. It follows that
\begin{equation}
\begin{split}
\frac{1}{c(\eta_1)E}\widehat{\mathcal{B}_{\epsilon}f}\paren{E,\eta_1,\Phi(x_2)}&=\frac{1}{c(\eta_1)E}\int_{E_m}^{E}K_1\paren{q,E}\hat{f}(q,\eta_1,x_2)\mathrm{d}q+\hat{f}(E,\eta_1,x_2).
\end{split}
\end{equation}
We have
$$K(q,E)=G\paren{\frac{c_2q}{E}}W_2\paren{g\left(\frac{c_2q}{E}\right),x_2},\ \text{where}\ G\paren{z}=\frac{z\cos\left(\eta_1 {g\left(z\right)}\right)}{h\left({g\left(z\right)}\right)},$$
and hence
$$K_1(q,E)=-\frac{c_2q}{E^2}\left[\paren{\frac{c_2q}{E}}^2\frac{W'_2\paren{x_1,x_2}\cos(\eta_1 x_1)}{h^2(x_1)}+G'\paren{\frac{c_2q}{E}}W_2\paren{x_1,x_2}\right]_{x_1=g\left(\frac{c_2q}{E}\right)},$$
where $W'_2$ denotes the first order derivative of $W_2$ in the $x_1$ variable. We will now show that $K_1$ is bounded on $T$, thus proving the invertiblity of the Bragg transform by Theorem \ref{volt}. Using the identity $g'(z)=\frac{z}{h(g(z))}$ (equation \eqref{box2}), and after direct application of the chain and product rules to $G$ we have
\begin{equation}
\begin{split}
G'\paren{z}&=\frac{\cos\paren{\eta_1 g(z)}}{h\paren{g(z)}}-z\paren{\frac{\eta_1g'(z)\sin\paren{\eta_1g(z)}}{h\paren{g(z)}}+\frac{g'(z)h'\paren{g(z)}\cos\paren{\eta_1 g(z)}}{h^2\paren{g(z)}}}\\
&=\frac{h^2\paren{g(z)}\cos\paren{\eta_1 g(z)}}{h^3\paren{g(z)}}-\frac{\eta_1z^2\sin\paren{\eta_1g(z)}}{h^2\paren{g(z)}}-\frac{z^2h'\paren{g(z)}\cos\paren{\eta_1 g(z)}}{h^3\paren{g(z)}}\\
&=\left[\frac{\paren{h(x_1)^2-z^2h'\paren{{x_1}}}\cos\paren{\eta_1 {x_1}}}{h^3\paren{{x_1}}}-\frac{\eta_1z^2\sin\paren{\eta_1 {x_1}}}{h^2\paren{{x_1}}}\right]_{x_1=g(z)}.
\end{split}
\end{equation}
By assumption, $W_2(\cdot,x_2)$ and $W_2'(\cdot,x_2)$ are bounded on $\left[-w(x_2),w(x_2)\right]$ for any $x_2\in \mathcal{I}$, and we have $\left|\frac{c_2q}{E}\right|<\frac{c_2E_M}{E_m}$ on $T$. Hence it suffices to show that $h'(g\paren{\frac{c_2q}{E}})$ is bounded above on $T$ and further that $h(g\paren{\frac{c_2q}{E}})$ is bounded above and away from zero on $T$ to prove the boundedness of $K_1$ on $T$.

Let $c_1=\frac{c_2E_m}{E_M}$ as before. We have $0<h_1(0)<|h_1(x_1)|<h_1\paren{w(x_2)}$ for $x_1\in \left[-w(x_2),w(x_2)\right]$. The polynomial $P_1(x_1)$ is strictly monotone increasing since
$$P'_1(x_1)=12x^2_1+4(1-x^2_2)+\Phi^2(x_2)\paren{3x^2_1+(x_2+1)(x_2+3)}>0$$
for all $x_1\in\mathbb{R}$ and $x_2\in\mathcal{I}$, and hence $0<P_1(g(c_1))<|P_1(x_1)|<P_1\paren{w(x_2)}$ for $x_1\in\left[g(c_1),w(x_2)\right]$. Now since $c_1<\frac{c_2q}{E}<c_2$ for $(q,E)\in T$, it follows that
\begin{equation}
\boxed{0<\frac{P_1(g(c_1))}{4h^3_1\paren{w(x_2)}}<\left|h\paren{g\paren{\frac{c_2q}{E}}}\right|<\frac{P_1\paren{w(x_2)}}{4h^3_1(0)}.}
\end{equation}
for $(q,E)\in T$. Let $P_2(x_1)=3x_1(\Phi(x_2)^2+2(1+x_2^2+x_1^2))$. Then the first derivative of $h$ is
\begin{equation}
\begin{split}
\boxed{h'(x_1)=\frac{P'_1(x_1)}{4h^{3}_1(x_1)}-\frac{P_1(x_1)P_2(x_1)}{4h^{5}_1(x_1)},}
\end{split}
\end{equation}
which is bounded on $\left[g(c_1),w(x_2)\right]$. Hence $h'(g\paren{\frac{c_2q}{E}})$ is bounded on $T$ and thus we can recover $\hat{f}(q,\eta_1,x_2)$ uniquely for all $q\in[0,E_M]$, $\eta_1\in\mathbb{R}\backslash\left\{\eta_1=\frac{\pi/2+j\pi}{w(x_2)} : j\in\mathbb{Z}\right\}$ and $x_2\in\mathcal{I}$, by Theorem \ref{volt}. As the Fourier transform is known almost everywhere, we can recover $f(q,\vx)$ for all $q\in\mathfrak{E}$ and $\vx\in\mathbb{R}\times\mathcal{I}$.
\end{proof}
\end{theorem}
\begin{remark}
\label{rem2}
Theorem \ref{main4} presents a generalization of the results of Theorem \ref{main1}, in the sense that we can set $\Phi=0$ in the calculation of $q_1$ in Theorem \ref{main4} to achieve the inversion results of Theorem \ref{main1}. If the $q_1(0)<\frac{q_1(w)E_m}{E_M}$ condition of Theorem \ref{main4} is not met, then the problem reduces (in equation \ref{eqb}) to the inversion of Volterra type equations, where both the upper and lower integration limits vary as functions of $E$. As such operators have little coverage in the literature, we require the $q_1(0)<\frac{q_1(w)E_m}{E_M}$ condition at this stage to use the theory of classical Volterra operators of the second kind \cite{Tric}.

To give a geometric explanation of the $q_1(0)<\frac{q_1(w)E_m}{E_M}$ condition we consider the curves of integration $\{(q,x_1)\in\mathbb{R}^2 : q=Eq_1(x_1)\}$ for varying $E\in [E_m,\frac{E_M}{q_1(w)}]$, where $x_2=0$, $E_m=0.15$, $w=1$ and $\epsilon=0.01$ is set to some small (relative to the tunnel length) offset. To clarify $q_1$ is defined as in equation \eqref{q13}. See figure \ref{Bgcurveoff}.
\begin{figure}[!h]
\centering
\begin{subfigure}{0.49\textwidth}
\includegraphics[width=0.9\linewidth, height=5.2cm, keepaspectratio]{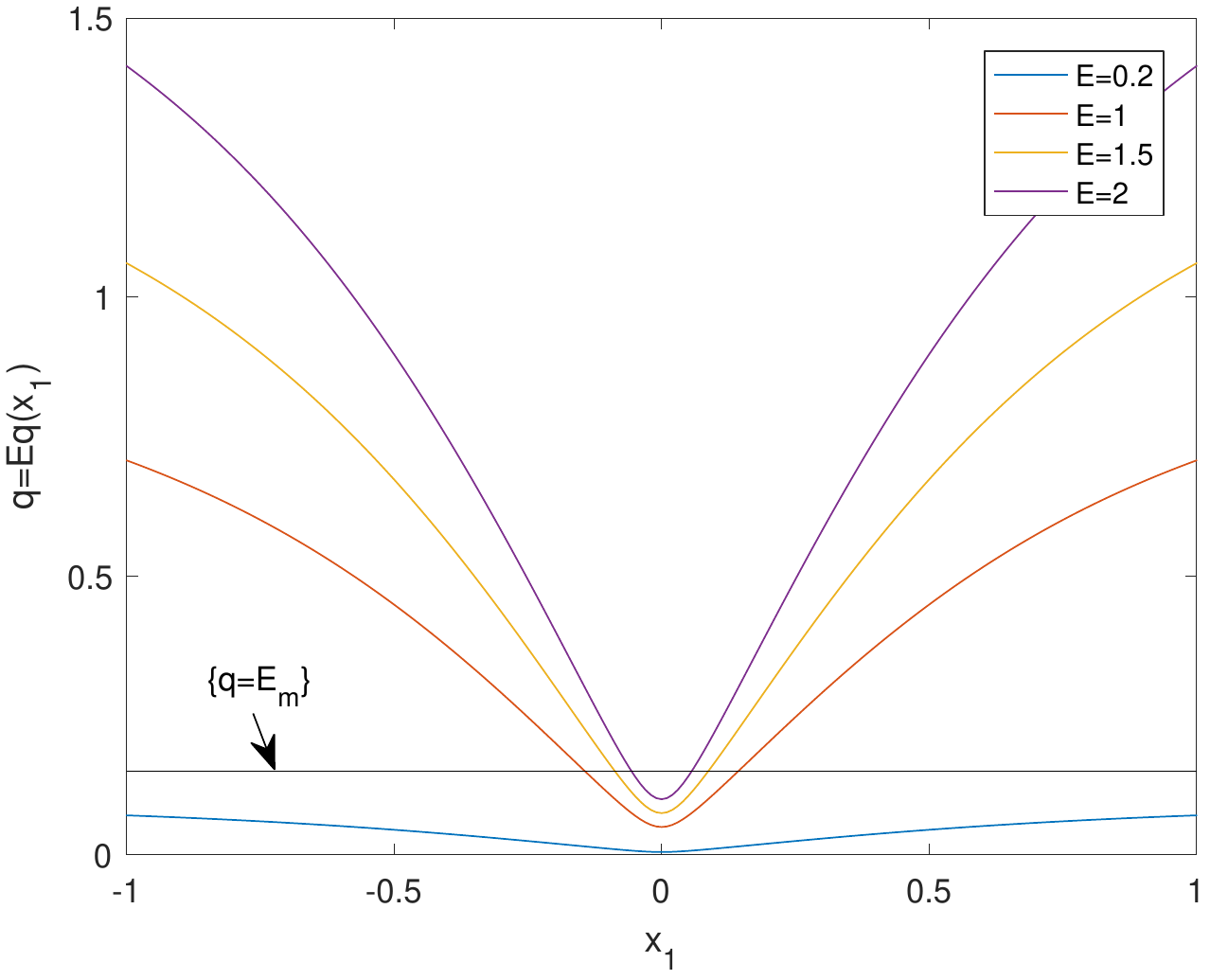}
\subcaption{$\frac{E_M}{q_1(w)}=2$}
\end{subfigure}
\begin{subfigure}{0.49\textwidth}
\includegraphics[width=0.9\linewidth, height=5.2cm, keepaspectratio]{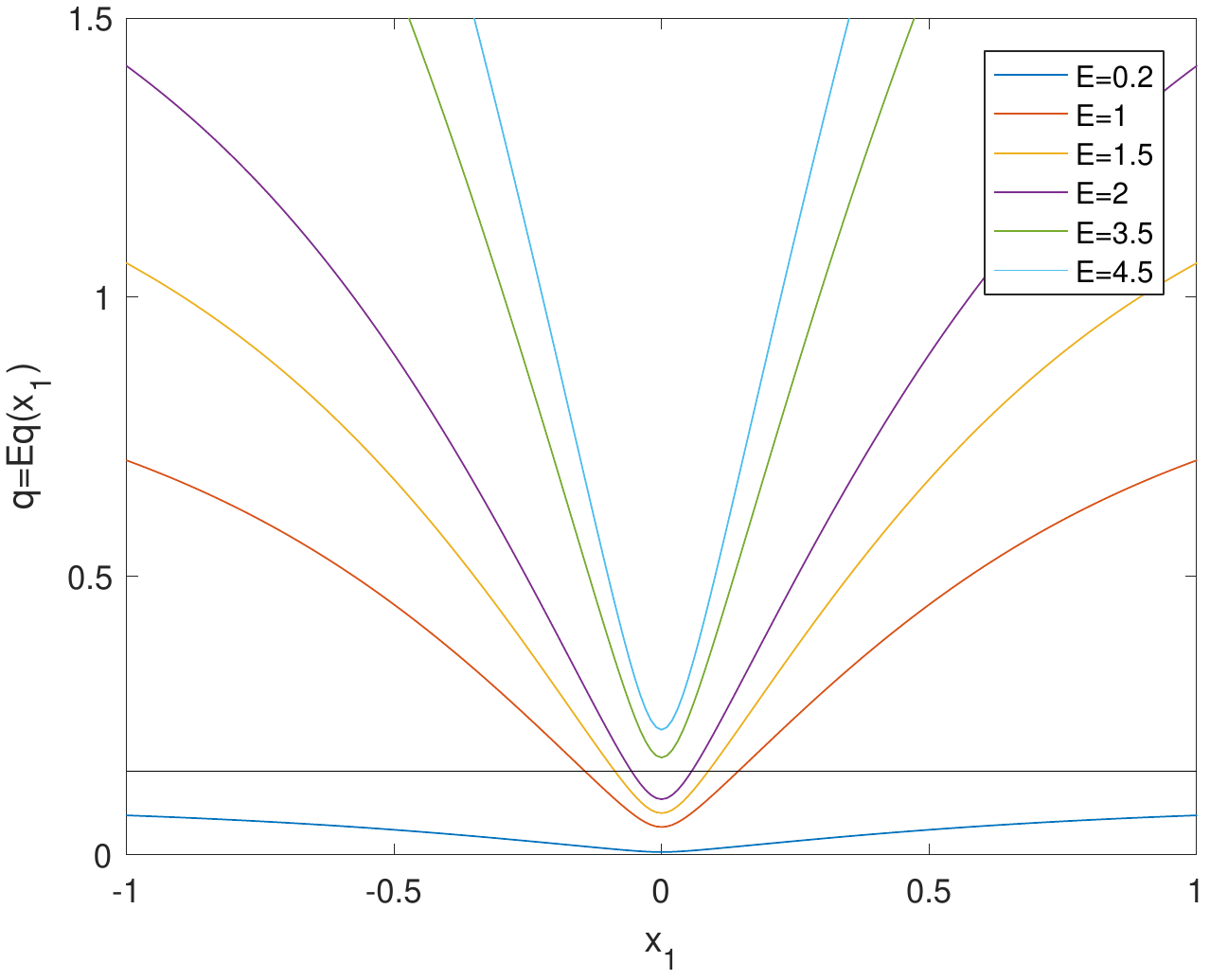}
\subcaption{$\frac{E_M}{q_1(w)}=4.5$}
\end{subfigure}
\caption{Plot of the curves of integration for the offset Bragg transform for varying $E$ and $\frac{E_M}{q_1(w)}$. The line $\{q=E_m\}$ is displayed in black.}
\label{Bgcurveoff}
\end{figure}
In this case the integration curves do not intersect the origin (in contrast to section \ref{bgsec1} and figure \ref{Bgcurve}) and the $q$ coordinate of the curve saddle point $q=Eq_1(0)$ increases monotonically with $E$. The $q_1(0)<\frac{q_1(w)E_m}{E_M}$ condition is equivalent to the statement ``all integration curve saddle points lie below the line $\{q=E_m\}$". That is $q_M=\frac{E_M}{q_1(w)}q_1(0)=\frac{q_1(0)E_M}{q_1(w)}<E_m$. In figure \ref{Bgcurveoff} we give examples of two integration curve sets, one which satisfies $q_M<E_m$ (left-hand figure) and one where $q_M>E_m$ (right-hand figure).
\end{remark}
\begin{corollary}
Under the further assumption that $W_1$ is bounded on $\mathfrak{E}$ and $\Phi'$ is bounded on $\mathcal{I}$, the offset Bragg transform $\mathfrak{B}_{\epsilon} : L^2_0(\mathfrak{E}\times\mathbb{R}\times \mathcal{I})\to L^2(\mathfrak{E}\times\mathbb{R}\times\Phi(\mathcal{I}))$ is a bounded operator.
\begin{proof}
This follows similar ideas to the proof of Corollary \ref{cor1} and is a consequence of the boundedness of the kernel $K$ in equation \eqref{Volte}. Note also that we have restricted the codomain of $\mathfrak{B}_{\epsilon}$ to $\mathfrak{E}$ to prove boundedness as in Corollary \ref{cor1}.
\end{proof}
\end{corollary}

\subsection{Assisting the machine design}
\label{MD}
In Theorem \ref{main4} we require that
\begin{equation}
\label{md}
\boxed{q_1(0)<\frac{q_1\paren{w(x_2)}E_m}{E_M}}
\end{equation}
for all pairs $(x_2,\Phi(x_2))$ for $x_2\in\mathcal{I}$, to acheive sufficient conditions for a unique solution. 


We can use the inequality \eqref{md} to assist in the machine design of the portal scanner, in the sense that we can determine a $x_2\to\epsilon$ map $\Phi$ such that \eqref{md} is satisfied. For example, let us consider energies up to $E_M=1\AA^{-1}$, and let the maximum lattice spacing  be $a_0=10\AA$. Then $E_m=\frac{1}{2a_0}=1/20\AA^{-1}$. 
Let 
\begin{equation}
S=\left\{(x_2,\epsilon) \in \mathcal{I}\times\mathbb{R}_+: q_1(0)<\frac{q_1\paren{w(x_2)}E_m}{E_M}\right\}
\end{equation}
denote the region of $(x_2,\epsilon)$ space for which \eqref{md} is satisfied; that is, the region in design space for which the injectivity condition is satisfied.
\begin{figure}[!h]
\begin{subfigure}{0.32\textwidth}
\includegraphics[width=0.9\linewidth, height=4cm, keepaspectratio]{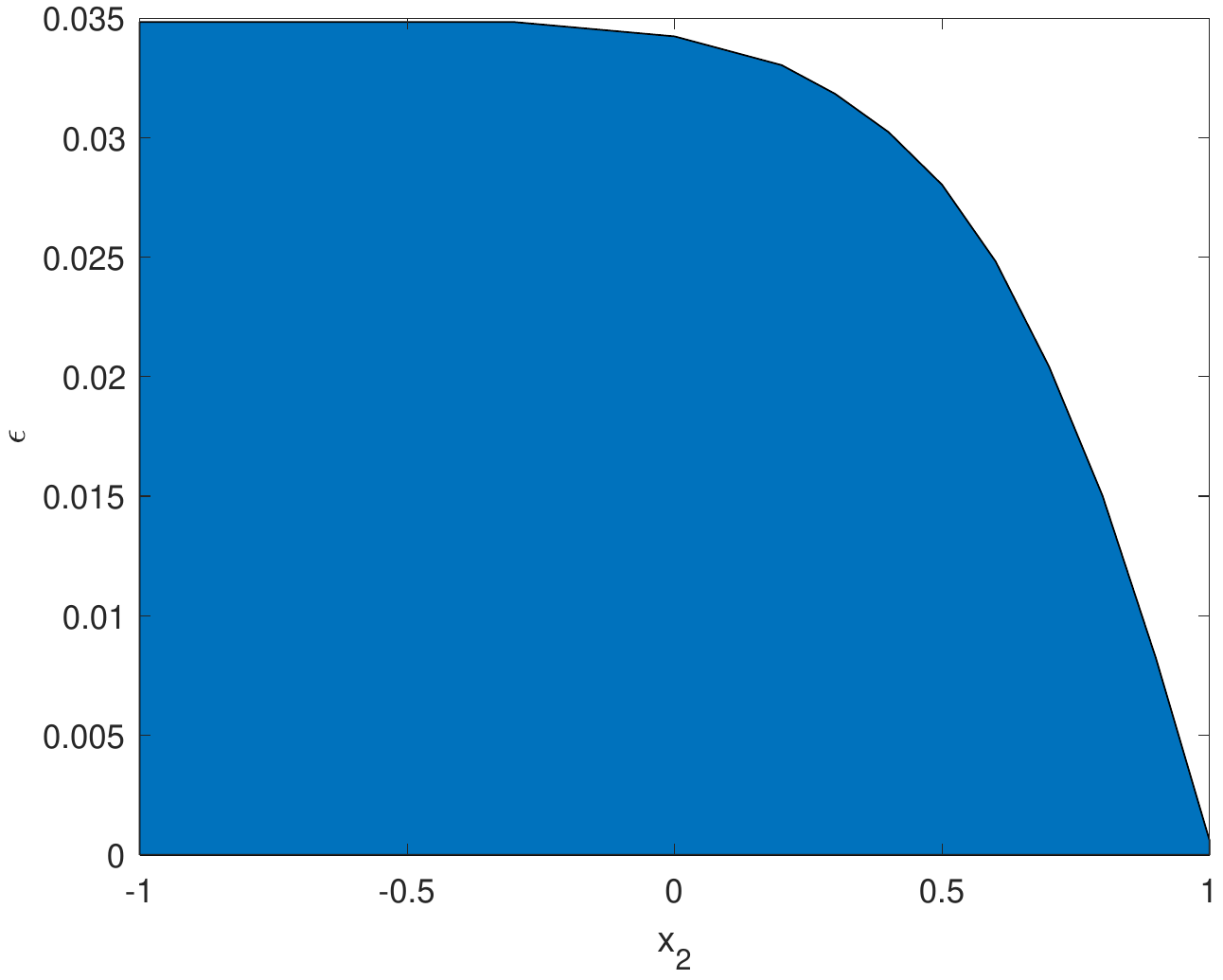}
\end{subfigure}
\begin{subfigure}{0.32\textwidth}
\includegraphics[width=0.9\linewidth, height=4cm, keepaspectratio]{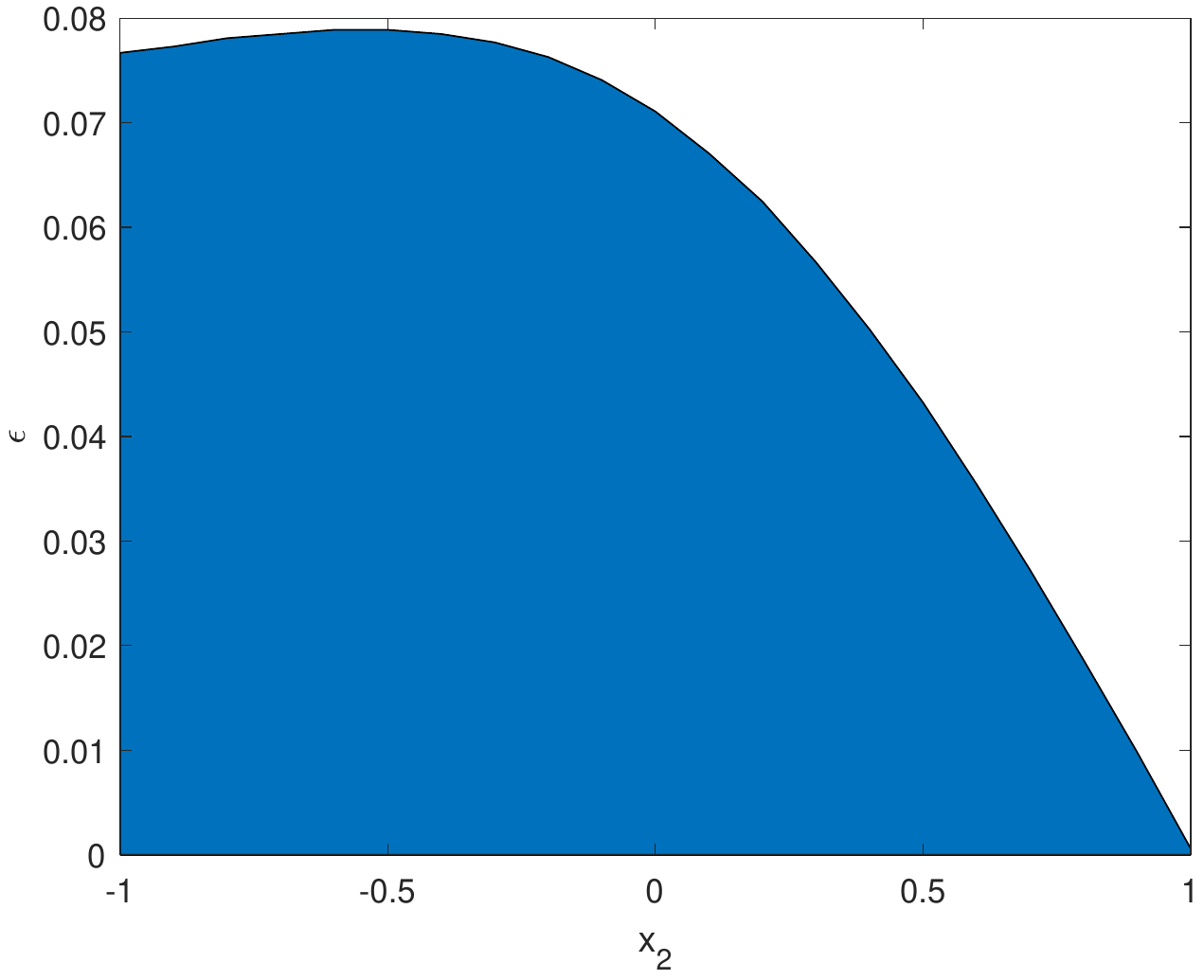} 
\end{subfigure}
\begin{subfigure}{0.32\textwidth}
\includegraphics[width=0.9\linewidth, height=4cm, keepaspectratio]{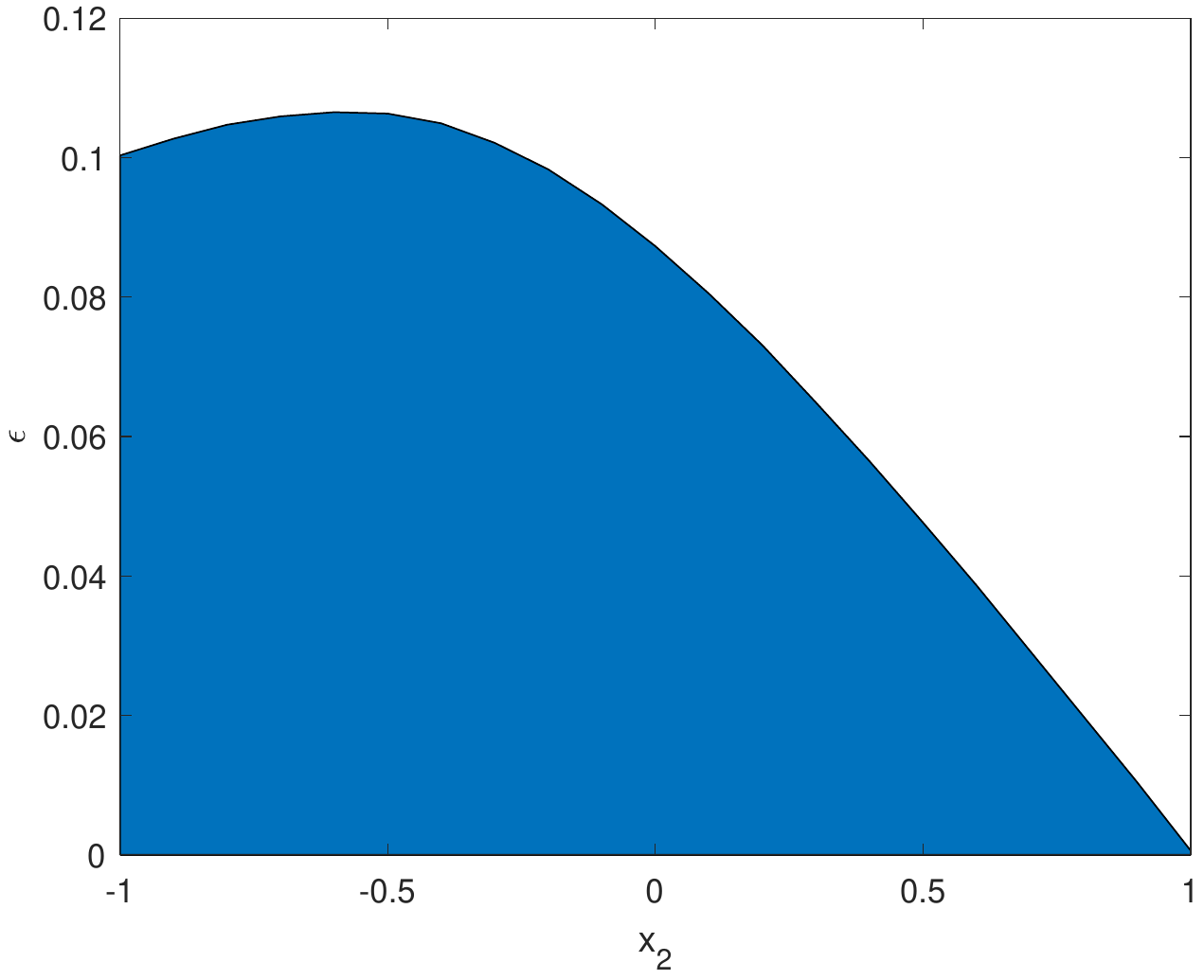}
\end{subfigure}
\begin{subfigure}{0.32\textwidth}
\includegraphics[width=0.9\linewidth, height=4cm, keepaspectratio]{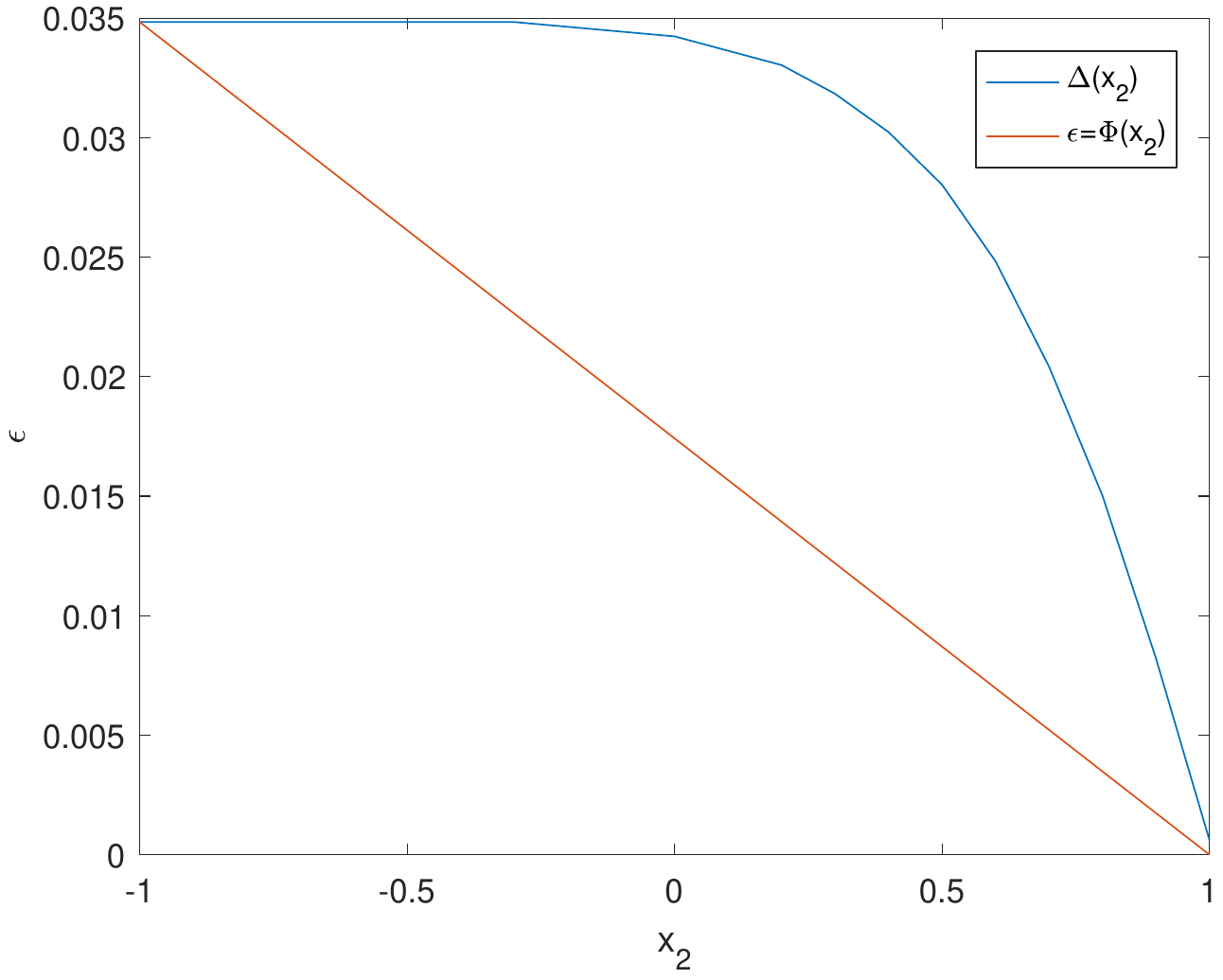}
\caption{$\beta=40^{\circ}$}
\end{subfigure}
\begin{subfigure}{0.32\textwidth}
\includegraphics[width=0.9\linewidth, height=4cm, keepaspectratio]{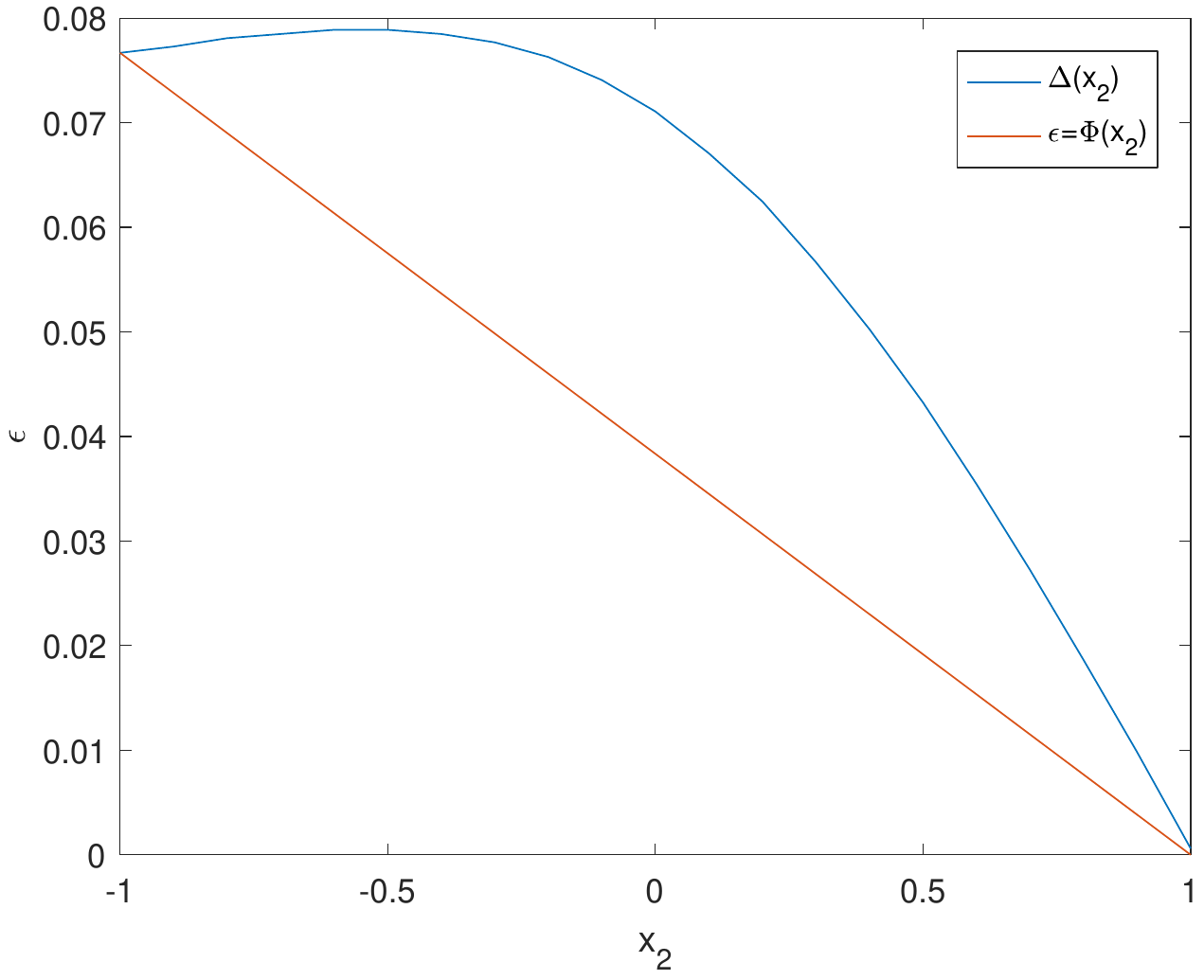} 
\caption{$\beta=90^{\circ}$}
\end{subfigure}
\begin{subfigure}{0.32\textwidth}
\includegraphics[width=0.9\linewidth, height=4cm, keepaspectratio]{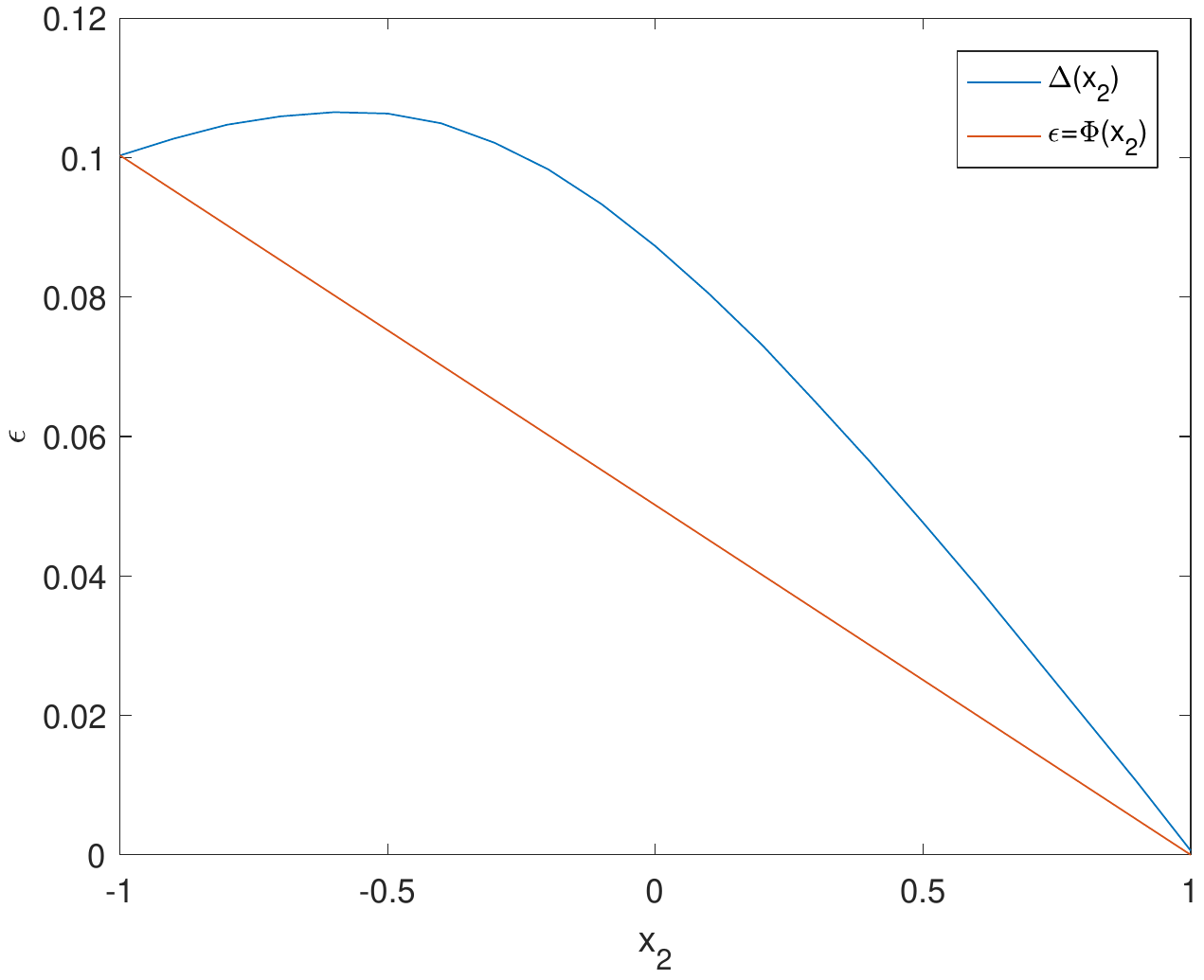}
\caption{$\beta=120^{\circ}$}
\end{subfigure}
\caption{Invertible design regions $S$ (the blue regions on the top row) and possible linear $\Phi$ (the red lines on the bottom row) for varying $\beta$. Note that the $\epsilon$ scales of the figures are different for each $\beta$.}
\label{S2}
\end{figure}
Then to be sure the Bragg transform is invertible we choose $\Phi$ such that $\Phi(\mathcal{I})\subset S$. See figure \ref{S2} where we have shown invertible design regions $S$ for varying source opening angles $\beta$. 

\begin{figure}[!h]
\begin{subfigure}{0.32\textwidth}
\includegraphics[width=0.9\linewidth, height=4cm, keepaspectratio]{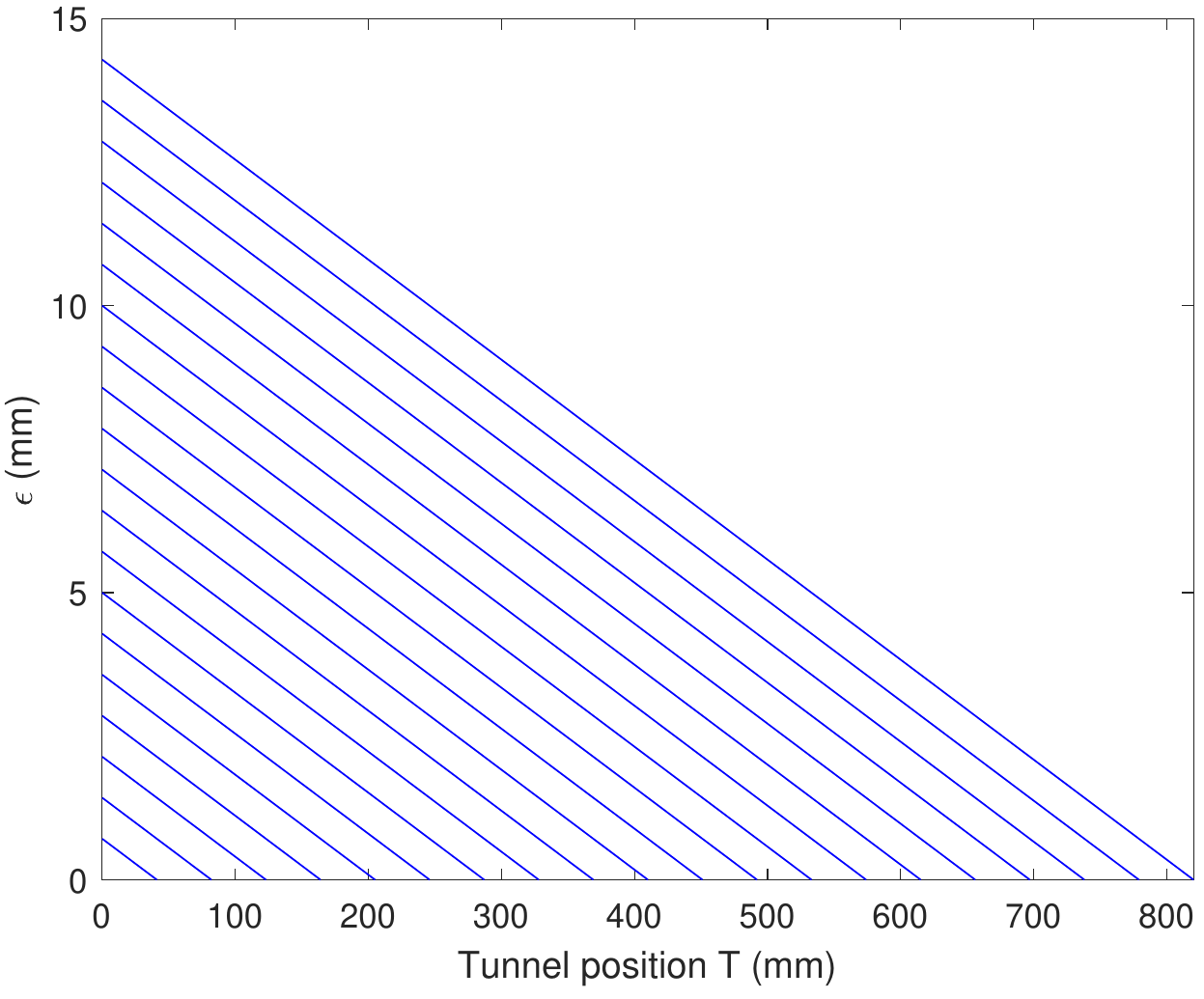}
\caption{$\beta=40^{\circ}$}
\end{subfigure}
\begin{subfigure}{0.32\textwidth}
\includegraphics[width=0.9\linewidth, height=4cm, keepaspectratio]{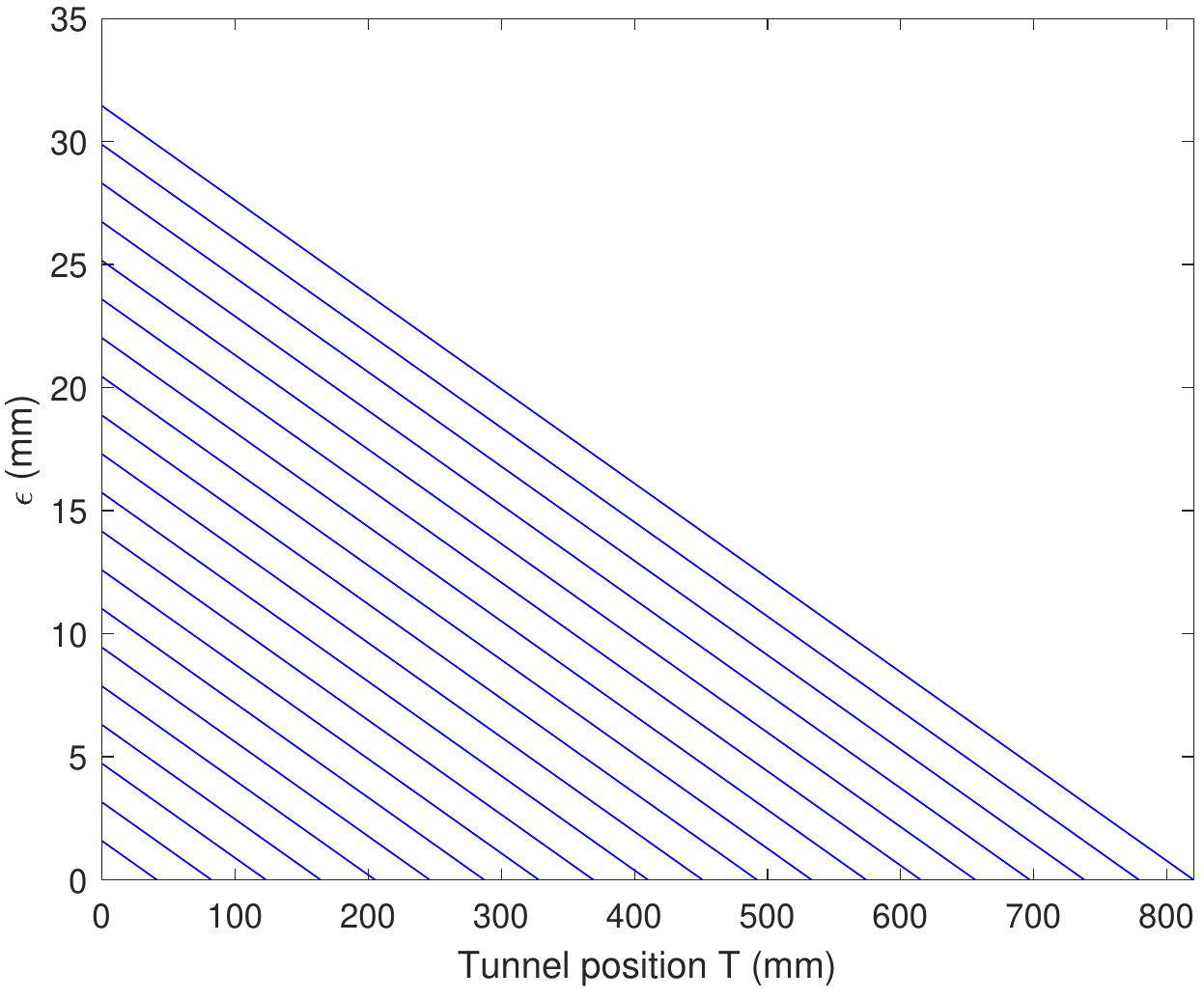} 
\caption{$\beta=90^{\circ}$}
\end{subfigure}
\begin{subfigure}{0.32\textwidth}
\includegraphics[width=0.9\linewidth, height=4cm, keepaspectratio]{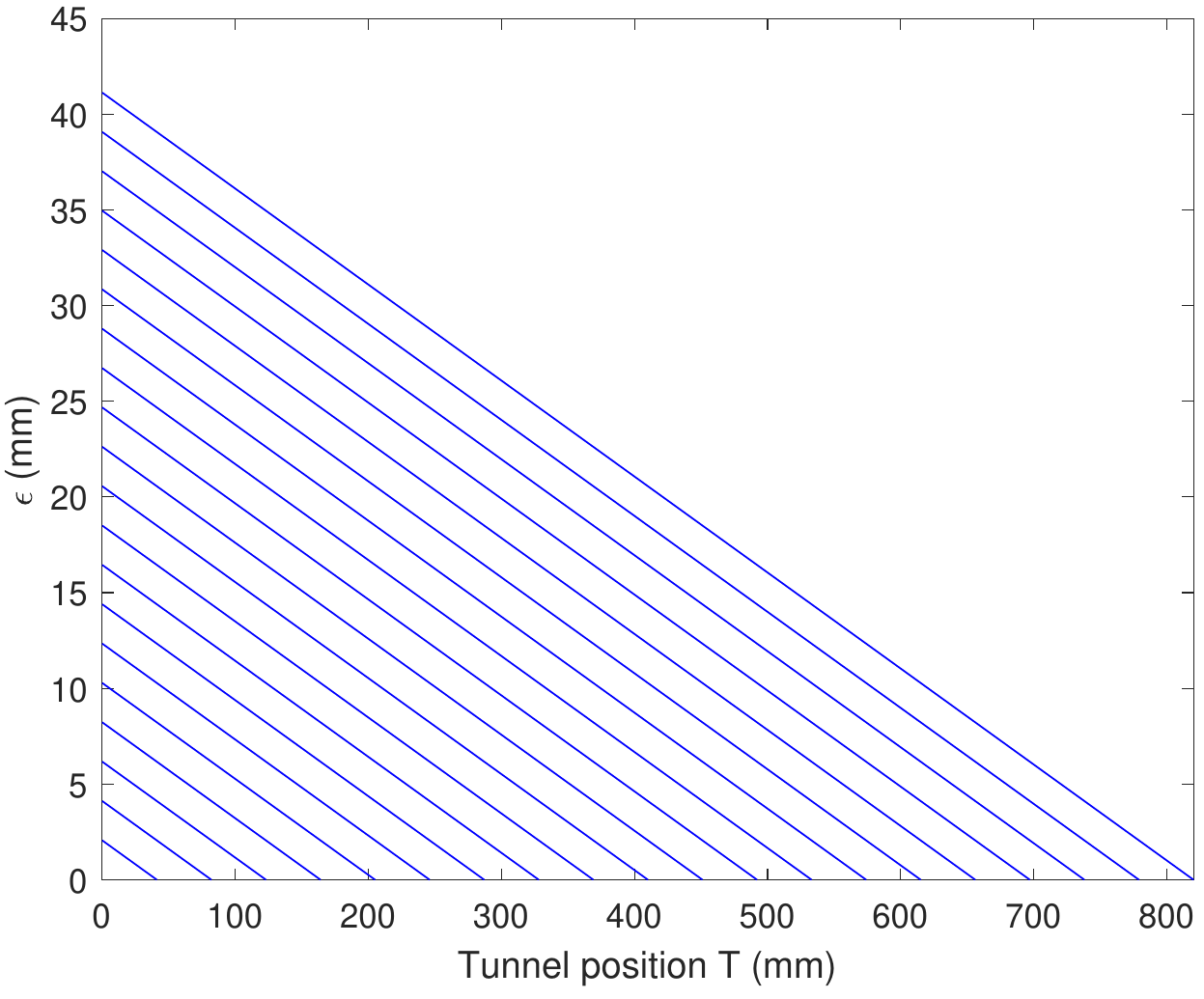}
\caption{$\beta=120^{\circ}$}
\end{subfigure}
\caption{Venetian blind detector configurations for varying source opening angles $\beta$. We show 21 detector arrays at $\epsilon=\Phi(x_2)$ for $x_2\in\left\{-1+\frac{j-1}{10} : 1\leq j\leq 21\right\}$, where for each $\beta$, $\Phi$ is the corresponding straight line relationship of figure \ref{S2}. The blue lines represent the collimation planes which intersect the tunnel at position $T=420(-x_2+1)$, where $x_2=\Phi^{-1}(\epsilon)$.}
\label{S4}
\end{figure}

We now give examples of linear $\Phi$ which satisfy $\Phi(\mathcal{I})\subset S$. Let $\Delta : \mathcal{I}\to\mathbb{R}_+$ be defined as $\Delta(\nu)=\max\paren{S\cap\{x_2=\nu\}}$. So $\Delta$ outputs the maximum value of $\epsilon\in S$ corresponding to a given $x_2=\nu$. Then setting $\Phi$ to be the straight line through $(-1,\Delta(-1))$ and $(1,\Delta(1))$ is sufficient to satisfy \eqref{md}. See figure \ref{S2}. Using the linear $\Phi$ maps we can simulate Venetian blind detector array configurations. To do this we convert to more practically useful units. We let the tunnel position be $T=420(-x_2+1)$mm (previously $x_2$), so the scanner length is 820mm. These dimensions are chosen based on prototype design specifics for the portal scanner of figure \ref{figmain}. The corresponding $\epsilon$ values are scaled in the same fashion. See figure \ref{S4}.
We see that as $\beta$ increases, the range of $\epsilon$ increases and we are allowed more freedom in the scanner design, while ensuring that the inversion results of Theorem \ref{main4} hold. For example when $\beta=40^{\circ}$ the maximum detector offset is $\epsilon=14$mm, which is small relative to the scanning tunnel length (820mm). In this case we may have problems in the construction. Conversely when $\beta=120^{\circ}$ the maximum $\epsilon=41$mm and we have more space to offset the detectors. The drawback of using a larger $\beta$ is that the Bragg inversion is less stable, as discussed in Remark \ref{rem1}. We note that the analysis presented here only covers the case when $E_m=1/20\AA^{-1}$ and $E_M=1\AA^{-1}$. We may decide upon different $E_m$ and $E_M$ in future work (depending on the application) which will alter the invertible design region.

{
\begin{discussion}
\label{dis1}
The injectivity results presented in this paper are proven by direct inversion of the Bragg operators, whereby we show that the solution for the Fourier components of $f$ can be written explicitly as a Neumann series. The proofs presented thus lay out novel inversion methods for BST. 

The injectivity results are important since they remove the concern for image artifacts due to null space in the reconstruction, and, with sufficient data and signal-to-noise ratio, they provide mathematical guarantee that the solution is a sensible representation of the ground truth.  Upon further analysis of the inversion steps we were able to gain some insight into the problem stability in Remark \ref{rem1}, and we have determined sufficient measurements for the portal scanner (e.g. construction of $\Phi$) such that $\mathfrak{B}_{\epsilon}$ is injective in section \ref{MD}. This in turn provides injectivity guarantees for $\mathfrak{B}_a$ (as in \eqref{equBG}) with $A_1$ and $A_2$ set to $A_1=A_2=1$ in the kernel (i.e. with attenuation neglected from the modeling). We would not suggest the inversion methods detailed in the proofs for practical implementation however, for two main reasons:
\begin{enumerate}
\item The inversion methods presented in Theorems \ref{main1}, \ref{main2}, \ref{main3} and \ref{main4}  require that data be discarded. That is, the Bragg scatter data $\mathfrak{B}_af$ \eqref{equBG} measured by the portal scanner is four-dimensional, whereas the inversion methods of the main theorems use a three-dimensional subset of the $\mathfrak{B}_af$ data (setting $s_1=d_1$), namely $\mathfrak{B}_n f$ (or $\mathfrak{B}_{\epsilon} f$), to prove injectivity. At this stage we have not found a way to implement our inversion methods using Volterra equations and Neumann series, which incorporates the full data $\mathfrak{B}_af$.
\item The recovery of the Fourier transform for $|\boldsymbol{\eta}|$ near the roots of $J(w|\boldsymbol{\eta}|)$ (with $J$ defined as in \eqref{J(u)}). This would likely require significantly limited sampling in the Fourier domain (particularly for larger $w$) as we cannot recover the Fourier transform near the roots of $J(w|\boldsymbol{\eta}|)$ due to division by zero (or values close to the machine precision).
\end{enumerate}
The main purpose  of the theorems presented in this paper is to prove injectivity of the Bragg transforms (i.e. to show the solution is unique) and to provide insight into the problem stability and the machine dimensions required for unique inversion. For a more rigorous analysis of the stability of the Bragg operators from a microlocal perspective, see \cite{webber2020microlocal}.

For practical reconstruction, we would recommend to use discrete methods. This way we can use the complete four-dimensional data, and factor in the a-priori information regarding the noise model (e.g. Poisson noise for photon arrivals with low counts) and target function (e.g. non-negativity) to better regularize the solution. In the next section we detail an algebraic reconstruction algorithm which implements a Poisson log-likelihood objective, TV regularization and non-negativity constraints, using $\mathfrak{B}_af$ as input.
\end{discussion}
}

{

\section{Image reconstructions}
\label{results}
Here we present image reconstructions from the Bragg transform data $\mathfrak{B}_{a}f$.  We choose to show reconstructions from $\mathfrak{B}_{a}f$ as this is the most practically relevant example considered in this paper (i.e. when compared to reconstruction from $\mathfrak{B}_{n}f$ or $\mathfrak{B}_{\epsilon} f$ data). As discussed in the introduction we operate under the assumption to neglect the attenuative effects, and set $A_1=A_2=1$ in \eqref{equBG}. First we establish the imaging parameters and reconstruction method, before giving our results.

\subsection{The imaging parameters}
\label{IP}
In this section we use the scaled measurement system introduced in section \ref{MD}, motivated by current considerations in system design for airport security screening applications, where $x_2\in(0,820)$mm. The transformation mapping $\mathcal{I}\to (0,820)$ is given by $x_2\to 420(-x_2+1)$ as in the figure \ref{S4} caption. Additionally we set the portal scanner width in the $x_1$ direction as $600$mm with $x_1\in[-300,300]$mm. We set an upper $q$ limit of $q_{\text{max}}=1\AA^{-1}$ since $F(q,Z)\approx 0$ ($F$ is as defined as in equation \eqref{Bgm1})  for $q>1$, for $Z<20$, which are typically the atomic numbers of interest in threat detection and security screening applications (e.g. H-C-N-O compounds). Therefore the domain of the reconstruction target $f(q,x_1,x_2)$, using the scaled machine dimensions, is $[0,1]\AA^{-1}\times[-300,300]\text{mm}\times(0,820)\text{mm}$. We consider three values of $x_2$ (scanning profiles) here, namely $x_2\in\{205,410,615\}$mm, and we present separate 2-D image reconstructions of $f(\cdot,\cdot,205)$, $f(\cdot,\cdot,410)$, and $f(\cdot,\cdot,615)$. 

We make use of the full three-dimensional Bragg data here, for each line profile $x_2$ considered, and we employ a range of detector and source positions ($d_1$ and $s_1$ respectively), and photon energies $E$. For this study we use 31 source positions $s_1\in\{-300+20(j-1) : 1\leq j\leq 31\}$mm spaced at 20mm intervals, with source opening angle of $\beta=120^{\circ}$ (see figure \ref{figmain}). We choose $\beta= 120^{\circ}$ as the largest fan width considered in section \ref{MD}, as this allows for the most freedom in machine design, and is thus the most practical. Furthermore, the analysis of Remark \ref{rem1} indicates that the Bragg inversion is less stable for larger $\beta$, and hence we are considering the worst case scenario in terms of stability. 

The sources are modeled as polychromatic (e.g. an X-ray tube), and we consider 29 spectrum energies $E\in\{j : 1\leq j\leq 29\}$keV, with energy bin width 1keV. We use 600 detector positions $d_1\in\{-300+(j-1) : 1\leq j\leq 600\}$mm spaced at 1mm intervals. The $\epsilon$ coordinate of the detectors varies with $x_2$, and is determined by 
\begin{equation}
\epsilon=\Phi(x_2)=\frac{75}{820}x_2.
\end{equation}
The detectors are assumed to be energy-resolved with energy resolution 1keV, so we are able to distinguish between all 29 energies in our range $E\in[1,29]$keV. The $\Phi$ chosen here is based on preliminary measurements for the portal scanner.

\subsection{The materials considered}
\begin{figure}[!h]
\centering
\begin{subfigure}{0.4\textwidth}
\includegraphics[ width=1\linewidth, height=1\linewidth, keepaspectratio]{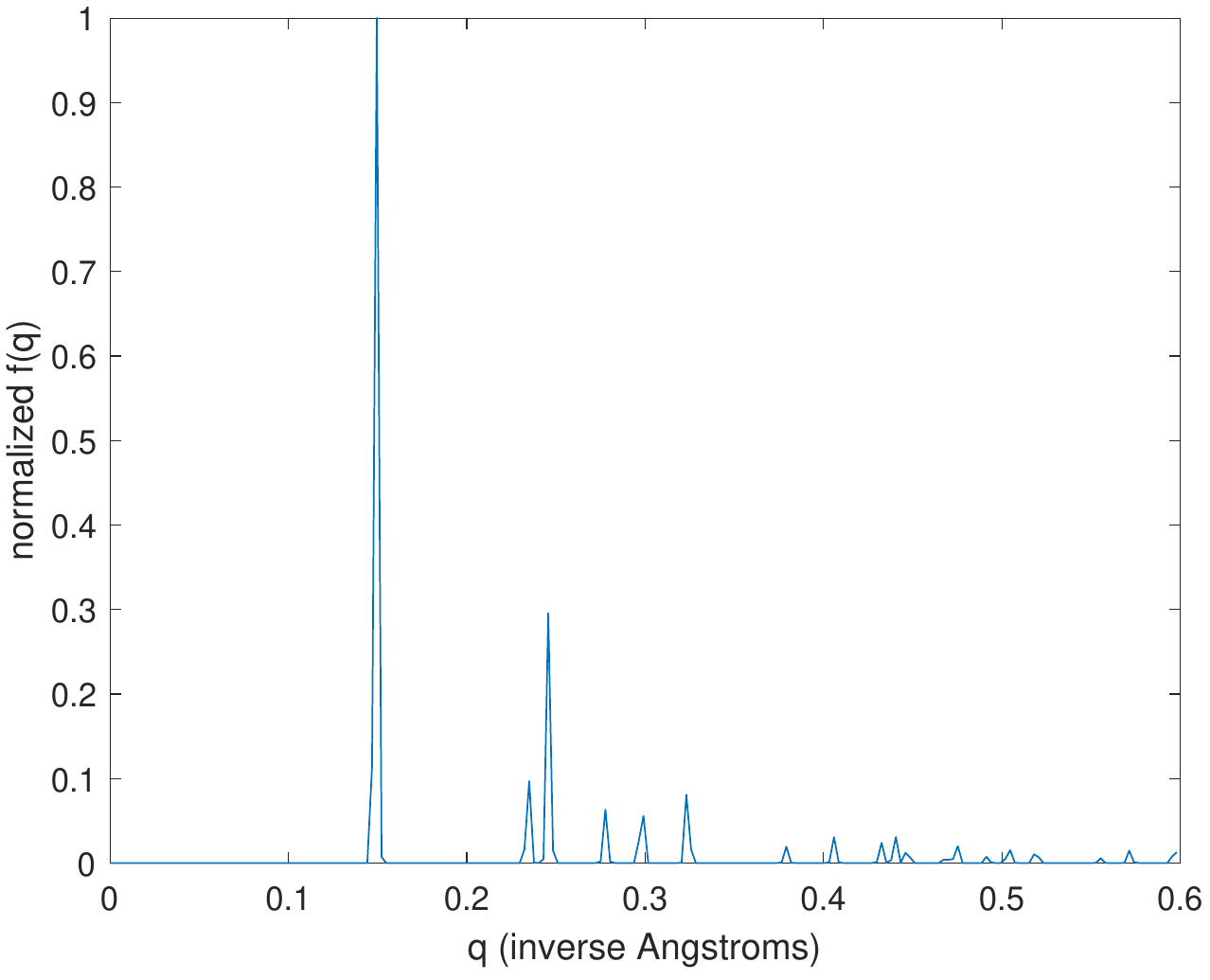}
\subcaption{C-graphite ($Z=6$)}\label{fqa}
\end{subfigure}
\begin{subfigure}{0.4\textwidth}
\includegraphics[ width=1\linewidth, height=1\linewidth, keepaspectratio]{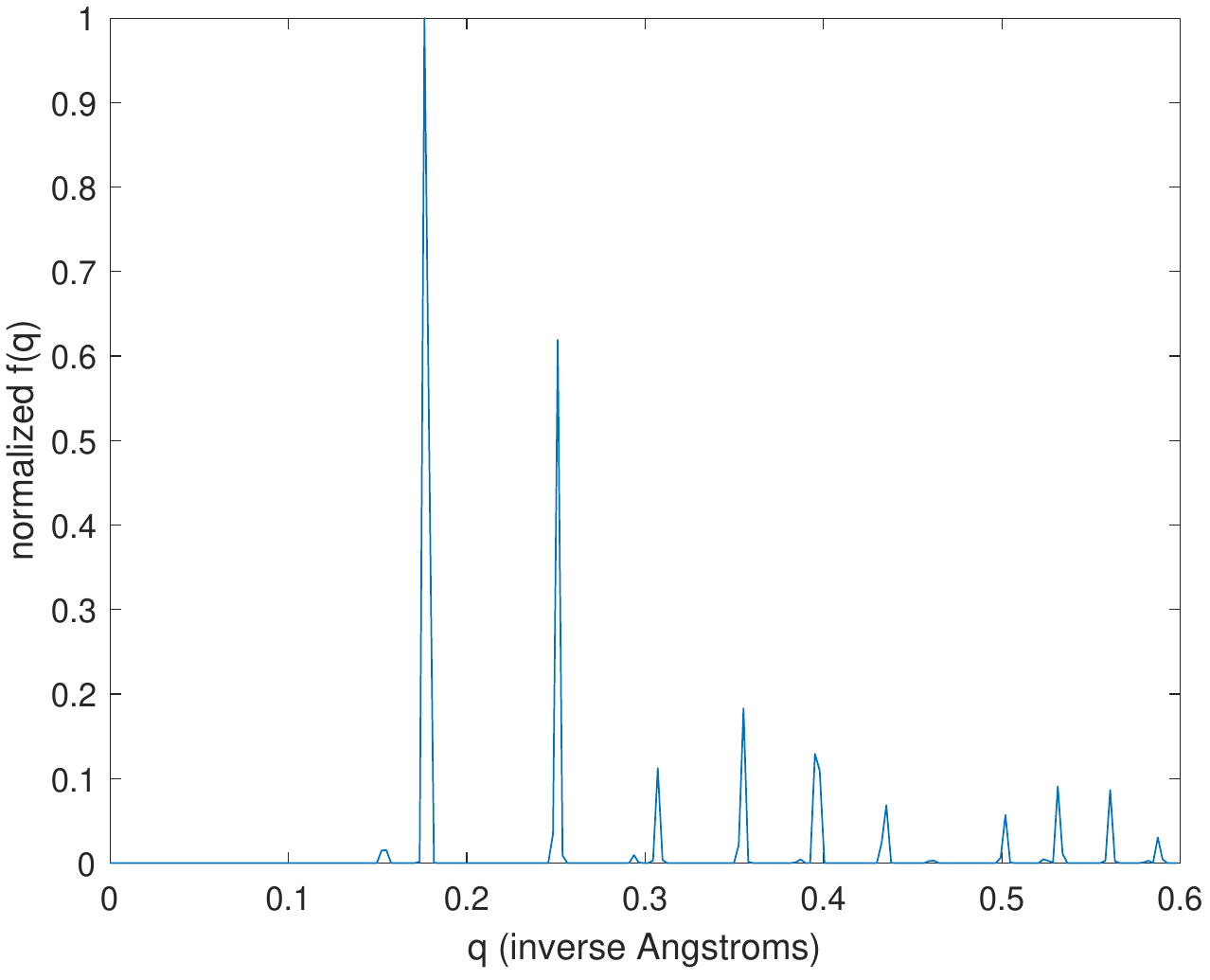} 
\subcaption{NaCl ($Z\approx 15$)}\label{fqb}
\end{subfigure}
\begin{subfigure}{0.4\textwidth}
\includegraphics[ width=1\linewidth, height=1\linewidth, keepaspectratio]{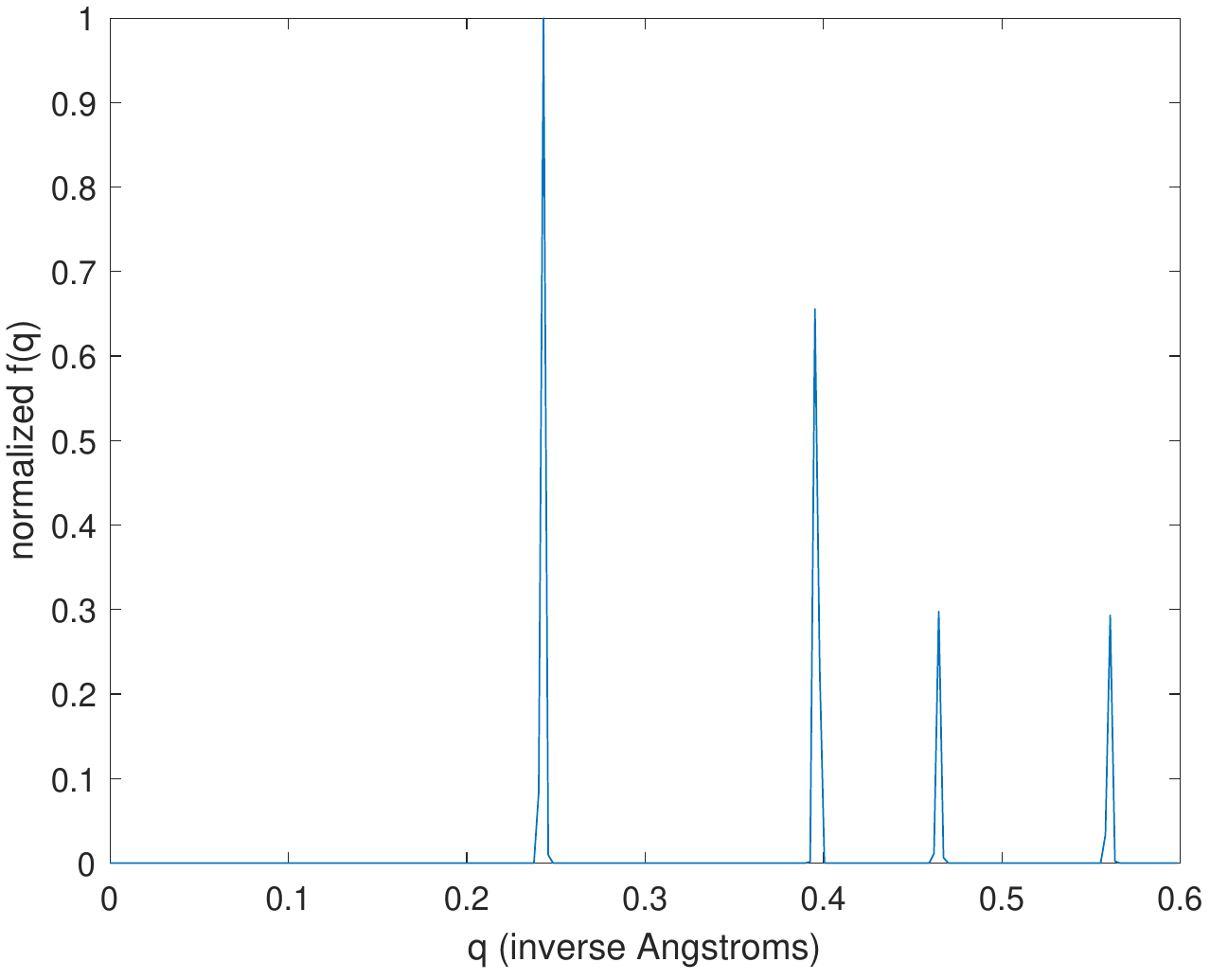}
\subcaption{C-diamond ($Z=6$)}\label{fqc}
\end{subfigure}
\begin{subfigure}{0.4\textwidth}
\includegraphics[ width=1\linewidth, height=1\linewidth, keepaspectratio]{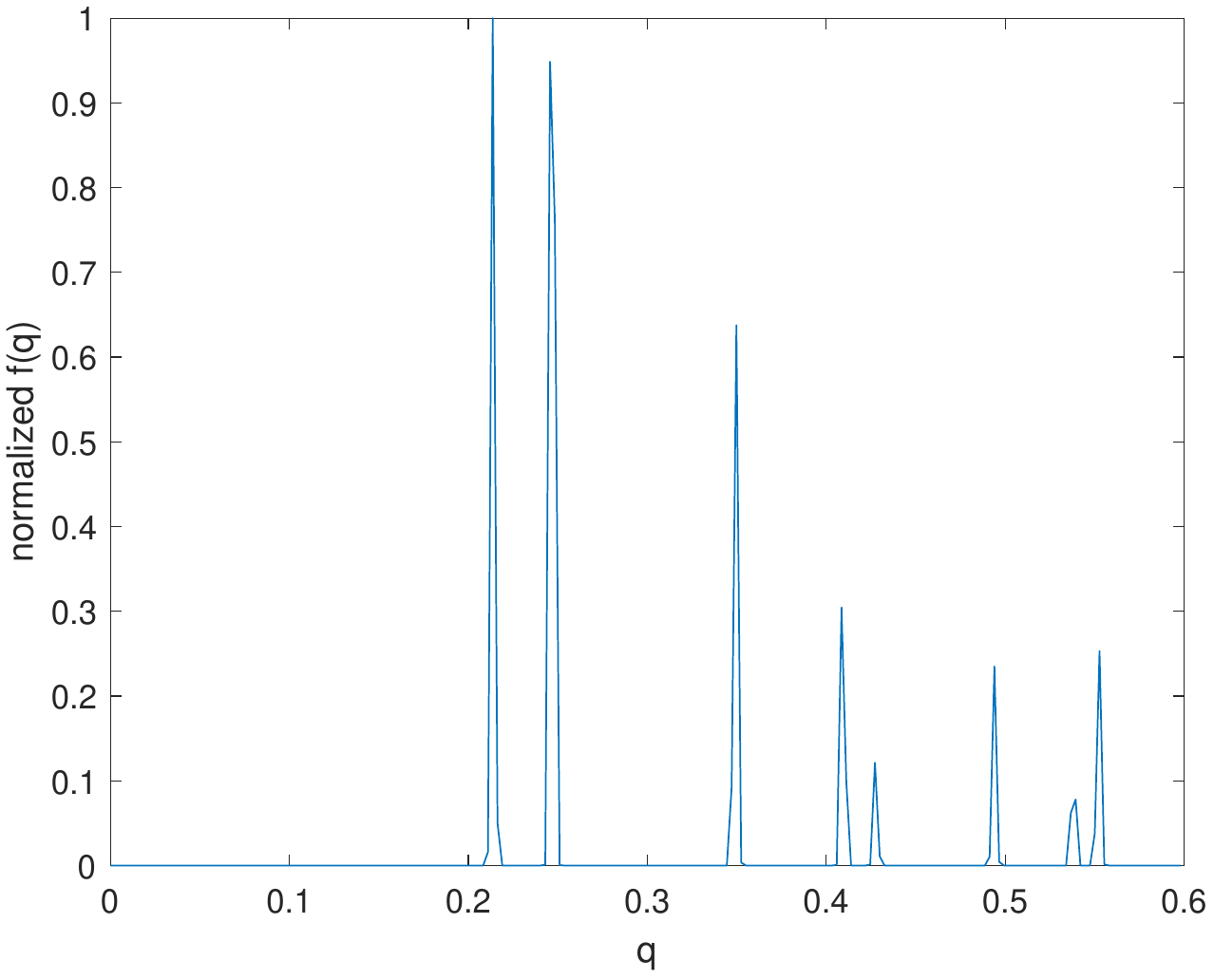}
\subcaption{Al ($Z=13$)}\label{fqd}
\end{subfigure}
\caption{$F(\cdot,Z)$ plots for varying $Z$. The plots are normalized in $L^{\infty}$ (by max value).}
\label{Fq1}
\end{figure}
We consider the $F(q,Z)$ curves for four crystalline powders here, namely Carbon with a graphite structure (denoted C-graphite), NaCl (salt), Carbon with a diamond structure (denoted C-diamond), and Aluminum. Here $F(q,Z)$ is defined as in equation \eqref{Bgm1}. For a given $Z$ let $Q=\text{supp}(F(\cdot,Z))\cap \{q<1\}=\{q_1,\ldots,q_{n_q}\}$ be set of $q$ values for which $F(\cdot,Z)$ is non-zero in the range $q\in[0,1]$, with $|Q|=n_q$. Then, in the simulations conducted here, we model $F$ as the Gaussian mixture
\begin{equation}
\label{BgG}
F(q,Z)\propto\sum_{j=1}^{n_q}F(q_j,Z)e^{-\frac{(q-q_j)^2}{\sigma^2}},
\end{equation}
where $\sigma^2=10^{-6}$ is chosen to be small relative to $q_{\text{max}}$ so that the Gaussians of \eqref{BgG} accurately represent the delta functions of \eqref{Bgm1}. See figure \ref{Fq1} for plots of the $F$ curves for C-graphite, NaCl, C-diamond and Al.

\subsection{The imaging phantoms}
We consider two imaging phantoms for reconstruction. See figure \ref{figPhan}. The ``2-spherical" phantom (on the top row of figure \ref{figPhan}) is comprised of an NaCl (salt) and C-diamond sphere with centers at $x_1=-100$mm and $x_1=100$mm respectively, both with radius $r=15$mm. The $x_2$ coordinate of the spheres is the same as the $x_2$ coordinate of the scanning profile. 
\begin{figure}[!h]
\centering
\begin{subfigure}{0.4\textwidth}
\includegraphics[ width=1\linewidth, height=1\linewidth, keepaspectratio]{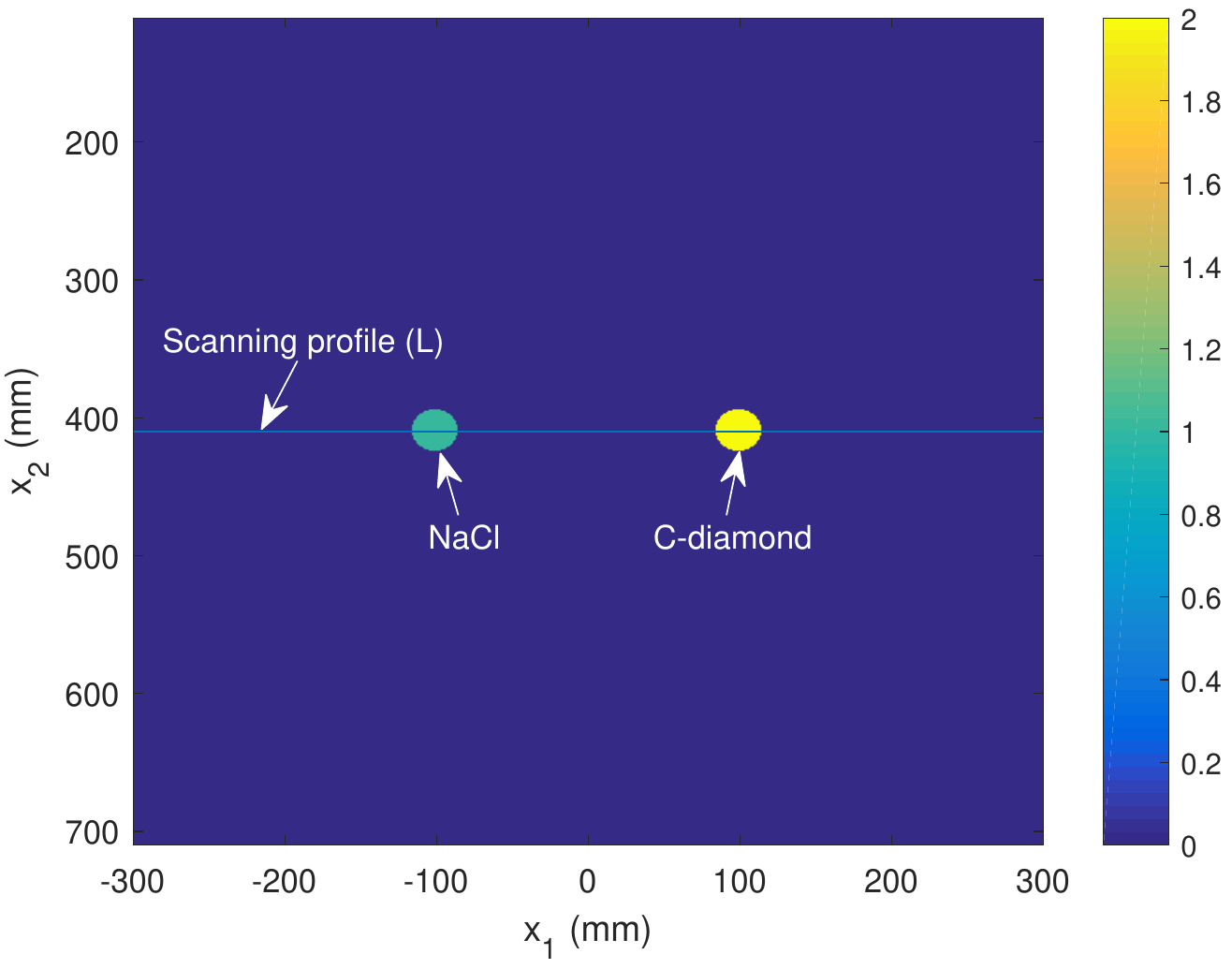} 
\subcaption{$(x_1,x_2)$ space representation}\label{2sphX}
\end{subfigure}
\begin{subfigure}{0.4\textwidth}
\includegraphics[ width=1\linewidth, height=1\linewidth, keepaspectratio]{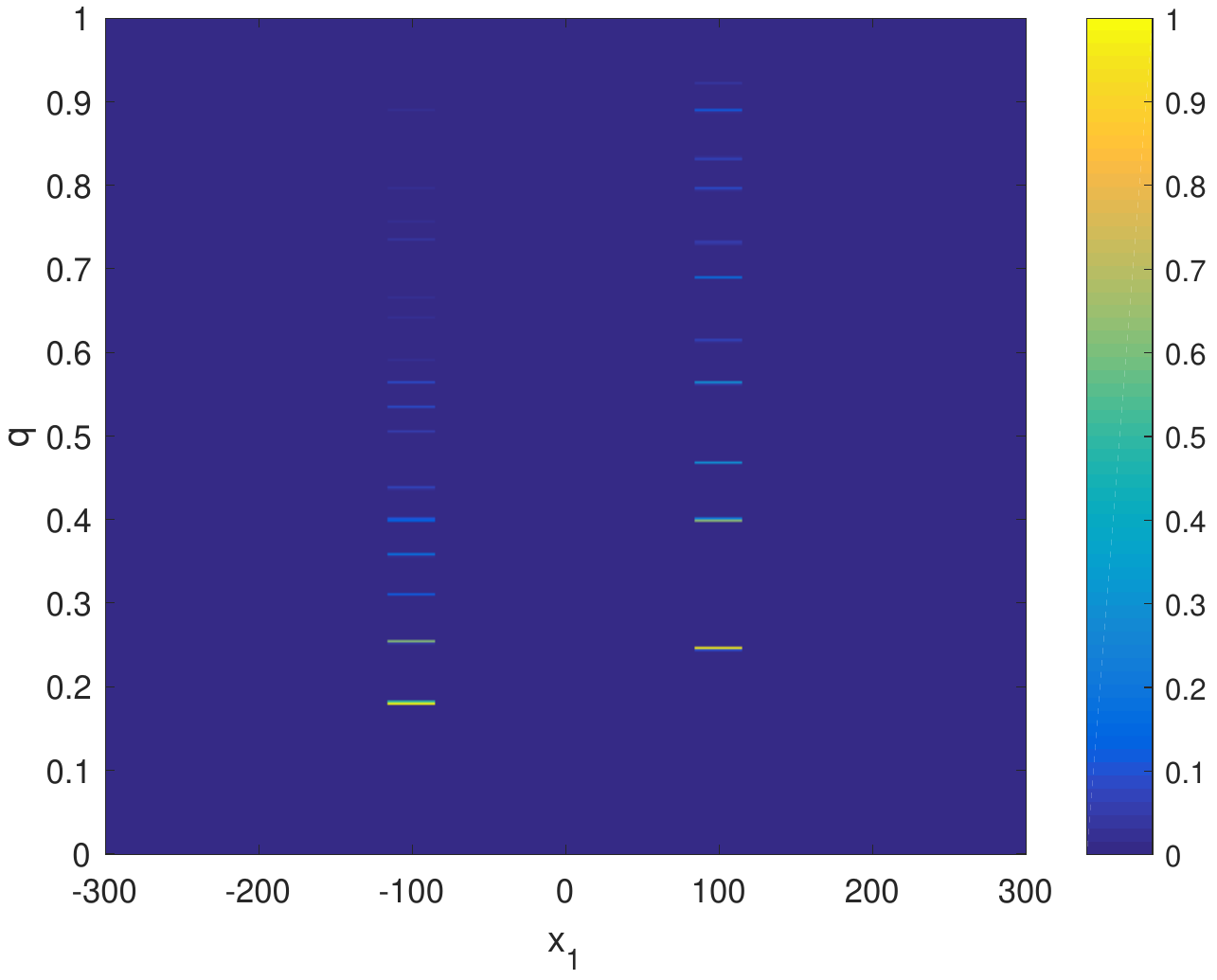} 
\subcaption{$f(\cdot,\cdot,410)$}\label{2sph}
\end{subfigure}
\begin{subfigure}{0.4\textwidth}
\includegraphics[ width=1\linewidth, height=1\linewidth, keepaspectratio]{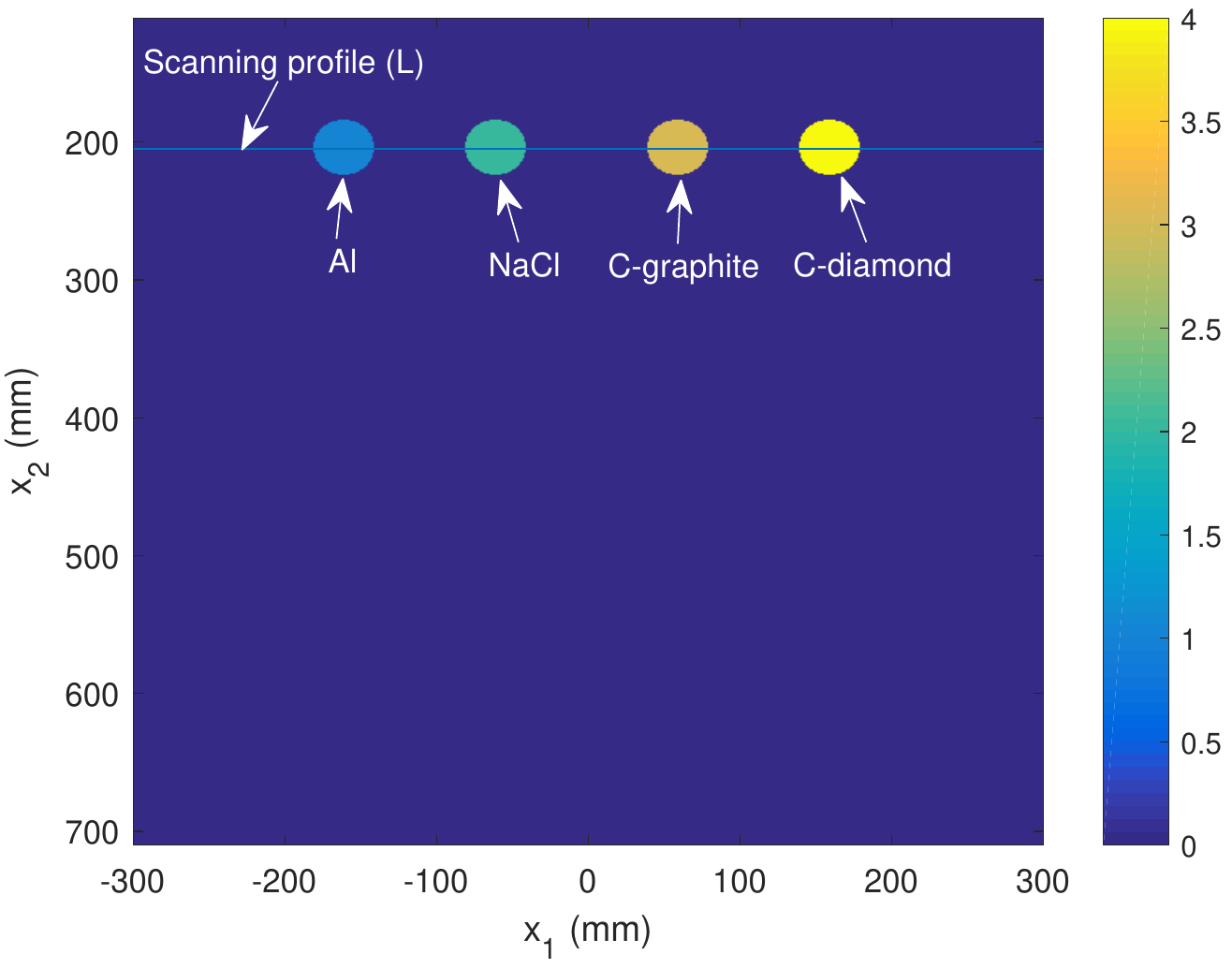} 
\subcaption{$(x_1,x_2)$ space representation}\label{4sphX}
\end{subfigure}
\begin{subfigure}{0.4\textwidth}
\includegraphics[ width=1\linewidth, height=1\linewidth, keepaspectratio]{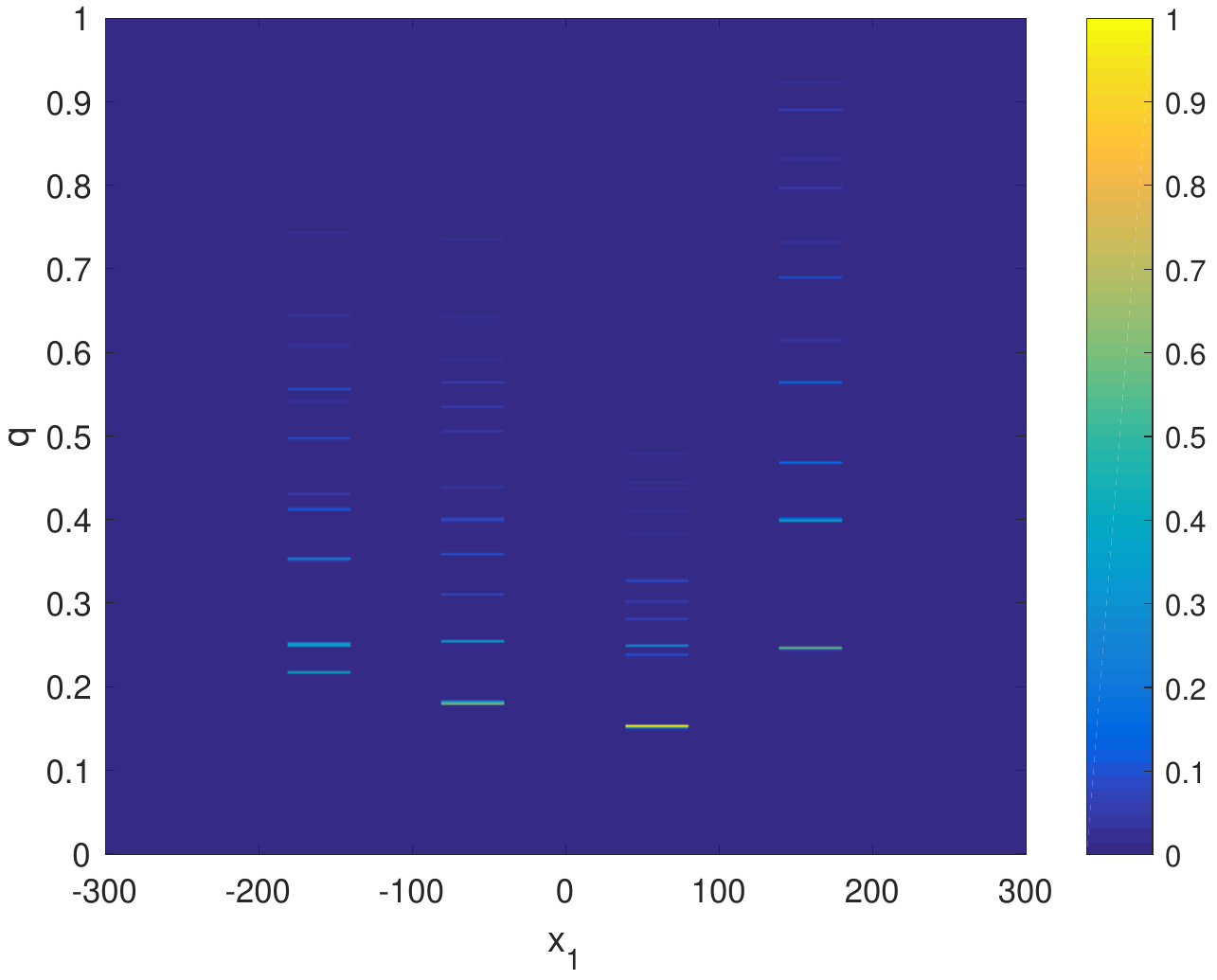} 
\subcaption{$f(\cdot,\cdot,205)$}\label{4sph}
\end{subfigure}
\caption{Top row -- 2-spherical phantom. Bottom row -- 4-spherical phantom. The image values in the left column are included for visualization (e.g. to distinguish between different materials) and have no physical meaning.}
\label{figPhan}
\end{figure}
Figure \ref{2sphX} shows the 2-spherical phantom in $(x_1,x_2)$ space, and the locations and sizes of the spheres. The sphere center $x_2$ coordinate and scanning profile is set to $x_2=0$ as an example. The corresponding $(q,x_1)$ space representation ($f(\cdot,\cdot,410)$) is shown in figure \ref{2sph}. The vertical line profiles of figure \ref{2sph} at $x_1=-100$mm and $x_1=100$mm correspond to the Bragg spectra of figures \ref{fqb} and \ref{fqc}, respectively.

The ``4-spherical" phantom (on the bottom row of figure \ref{figPhan}) is comprised of an Al, NaCl (salt), C-graphite and C-diamond sphere with centers at $x_1=-160$mm, $x_1=-60$mm, $x_1=60$mm and $x_1=160$mm, respectively, all with radius $r=20$mm. In this case, figure \ref{4sphX} shows the locations of the spheres in $(x_1,x_2)$ space, with the sphere center $x_2$ coordinate and scanning profile set to $x_2=205$mm. The corresponding $(q,x_1)$ space representation ($f(\cdot,\cdot,205)$) is shown in figure \ref{4sph}. To clarify, we consider three scanning profiles $x_2\in\{205,410,615\}$mm in total for reconstruction, and translate the sphere center $x_2$ coordinates to the scanning profile $x_2$ in each case. 

The 2-spherical phantom is included as a simple reconstruction target, with smaller and more spaced out objects, and the 4-spherical phantom is a more complex reconstruction target with larger objects in close proximity. Therefore we expect to see better results on the 2-spherical phantom and vice-versa. 

For the phantoms considered $E_m=0.15\AA^{-1}$, $E_M=1\AA^{-1}$, and we can check that the $\Phi$, $\beta$, $E_m$, and $E_M$ used here is sufficient to satisfy condition \ref{md}. Thus the recovery of $f$ from $\mathfrak{B}_{\epsilon}f$ is unique for all reconstructions presented in this section. 

\subsection{Quantitative analysis}
As a quantitative measure of the accuracy of our results we use the $F_1$ score on the image gradient, as is done
in \cite{webber2020joint,andrade2019shearlets} (see also the DICE metric \cite{taha2015metrics}). The gradient (or edge) $F_1$ score is a measure of how well we have recovered the gradient image $\nabla f$ (i.e. the edge map of $f$) and the locations of the Bragg peaks (i.e. the spikes in the Bragg spectra of figure \ref{Fq1}) in the reconstruction. We focus on the locations of the Bragg peaks as they are the most important for materials characterization \cite{DS}. 
The $F_1$ score values are in the range $[0,1]$. An $F_1$ score closer to one indicates a more accurate gradient reconstruction and an $F_1$ score closer to zero indicates a less accurate reconstruction.

\subsection{Data generation}
Let $A\in \mathbb{R}^{p\times (nm)}$ denote the discretized Bragg operator $\mathfrak{B}_{\epsilon}$, let $\vy\in\mathbb{R}^{nm}$ denote the discretized $f(\cdot,\cdot,x_2)$, for some fixed $x_2\in\{205,410,615\}$, and let $\vb\in\mathbb{R}^p$ denote the Bragg integral data $\mathfrak{B}_{\epsilon}f$. Note that the operator $A$ varies with $x_2$, so we are considering three different $A$ operators in total, which correspond to $x_2\in\{205,410,615\}$.  Here $p$ is the number of data points, $n$ is the number of $q$ samples, and $m$ is the number of $x_1$ samples. So $\vy$ is a vectorized $n\times m$ image which represents $f(\cdot,\cdot,x_2)$. Then the noisy Bragg data $\vb_{\textbf{e}}$ is distributed according to the Poisson model
\begin{equation}
\vb_{\textbf{e}}\sim \text{Poisson}\paren{c_{\text{avg}}p\times\frac{A\vy}{\sum_{j=1}^p\paren{A\vy}_j}}.
\end{equation}
That is, the exact data $\vb=A\vy$ is scaled so that the mean entry is $\frac{\sum_{j=1}^pb_j}{p}=c_{\text{avg}}$, where it is then used as mean to a multivariate Poisson distribution.
\begin{figure}
\setlength{\tabcolsep}{5pt}
\begin{tabular}{ c|cc }  
  &  $c_{\text{avg}}=10$ & $c_{\text{avg}}=1$ \\ \hline \\[-0.4cm] 
\rotatebox{90}{\hspace{2.5cm} GT}   & 
   \includegraphics[ width=0.4\linewidth, height=0.4\linewidth, keepaspectratio]{2sph} &
   \includegraphics[ width=0.4\linewidth, height=0.4\linewidth, keepaspectratio]{2sph} \\
 \rotatebox{90}{\hspace{1.75cm} $x_2=0$mm}  & 
   \includegraphics[ width=0.4\linewidth, height=0.4\linewidth, keepaspectratio]{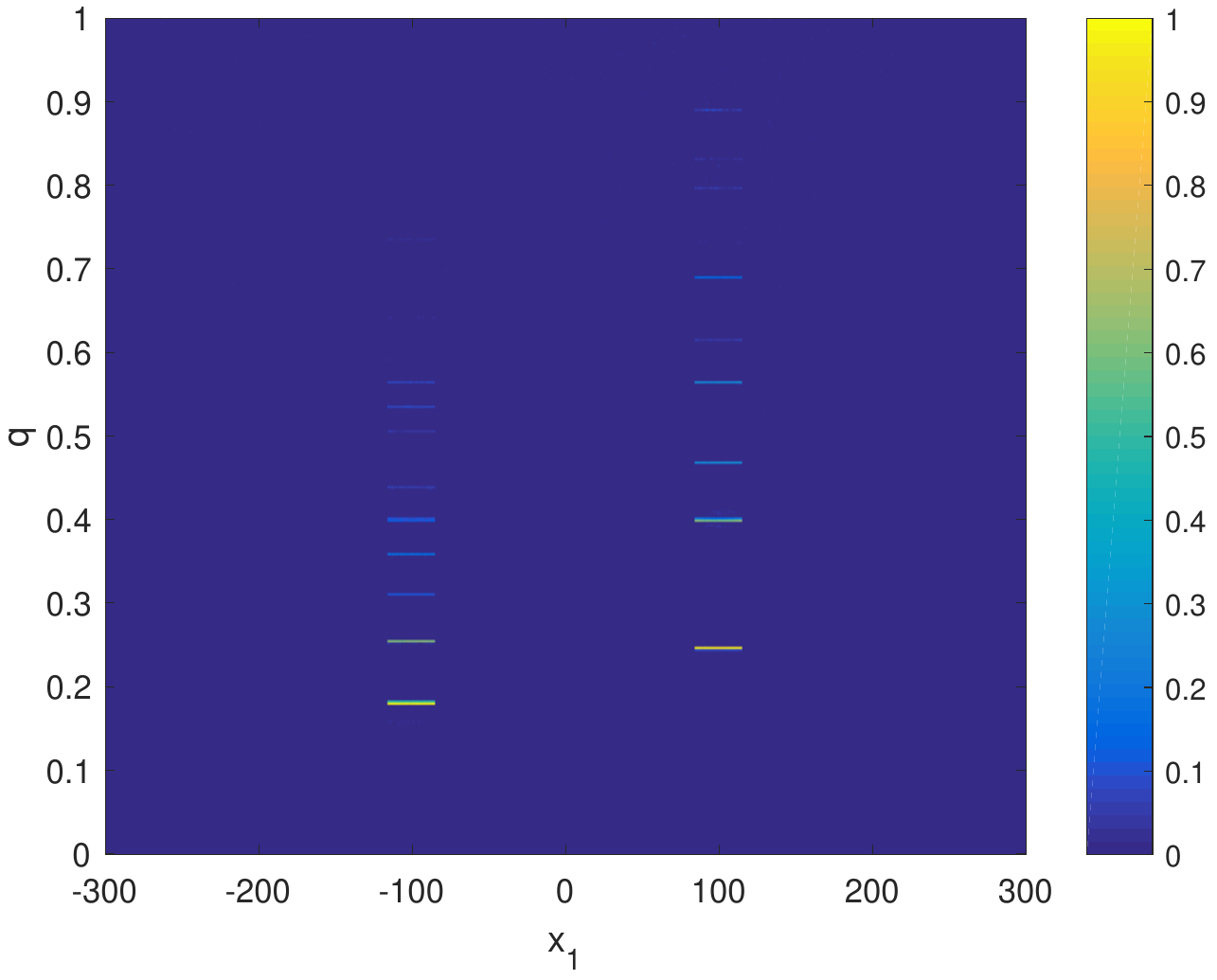} &
   \includegraphics[ width=0.4\linewidth, height=0.4\linewidth, keepaspectratio]{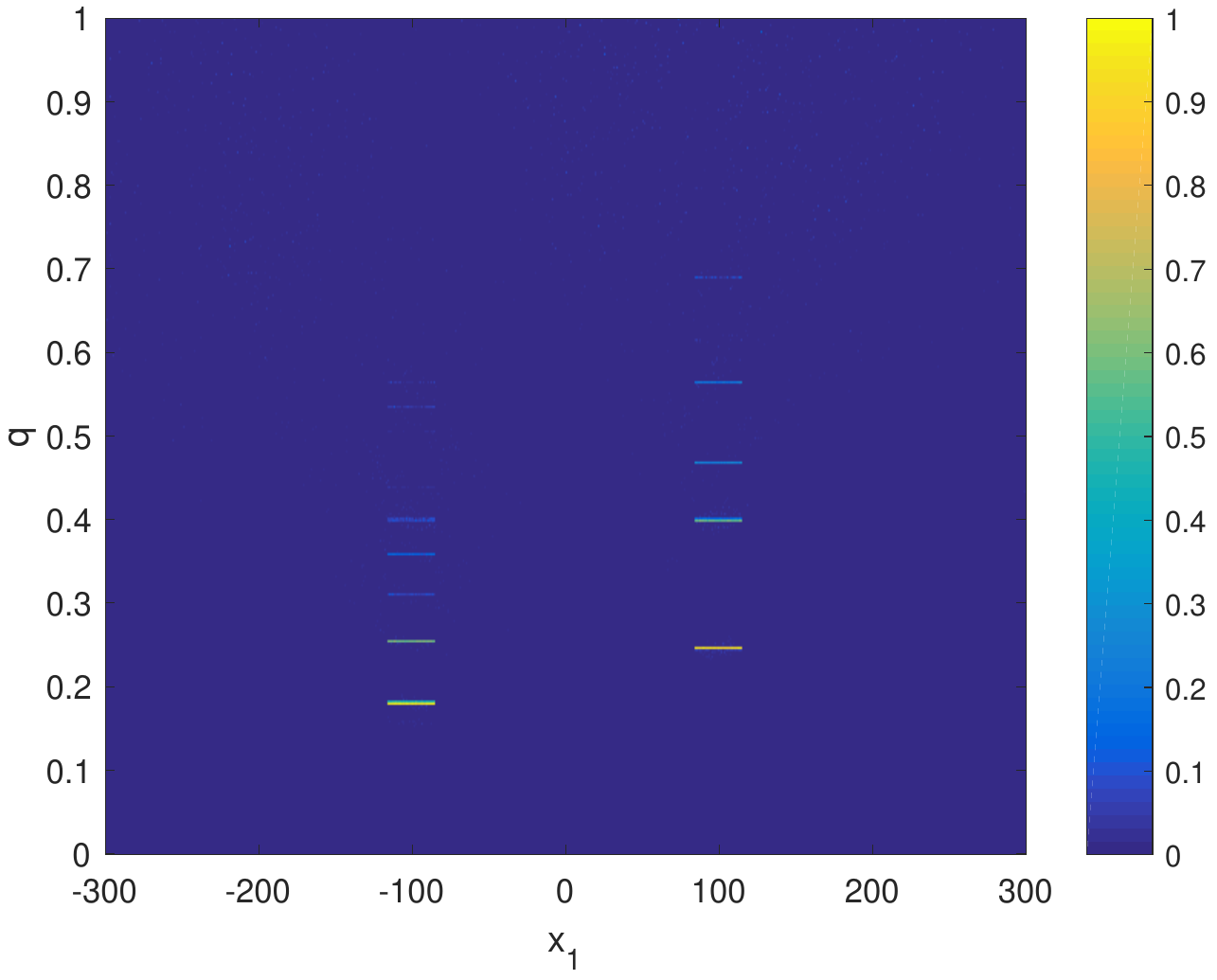} \\
 \rotatebox{90}{\hspace{1.5cm} $x_2=615$mm} &
  \includegraphics[ width=0.4\linewidth, height=0.4\linewidth, keepaspectratio]{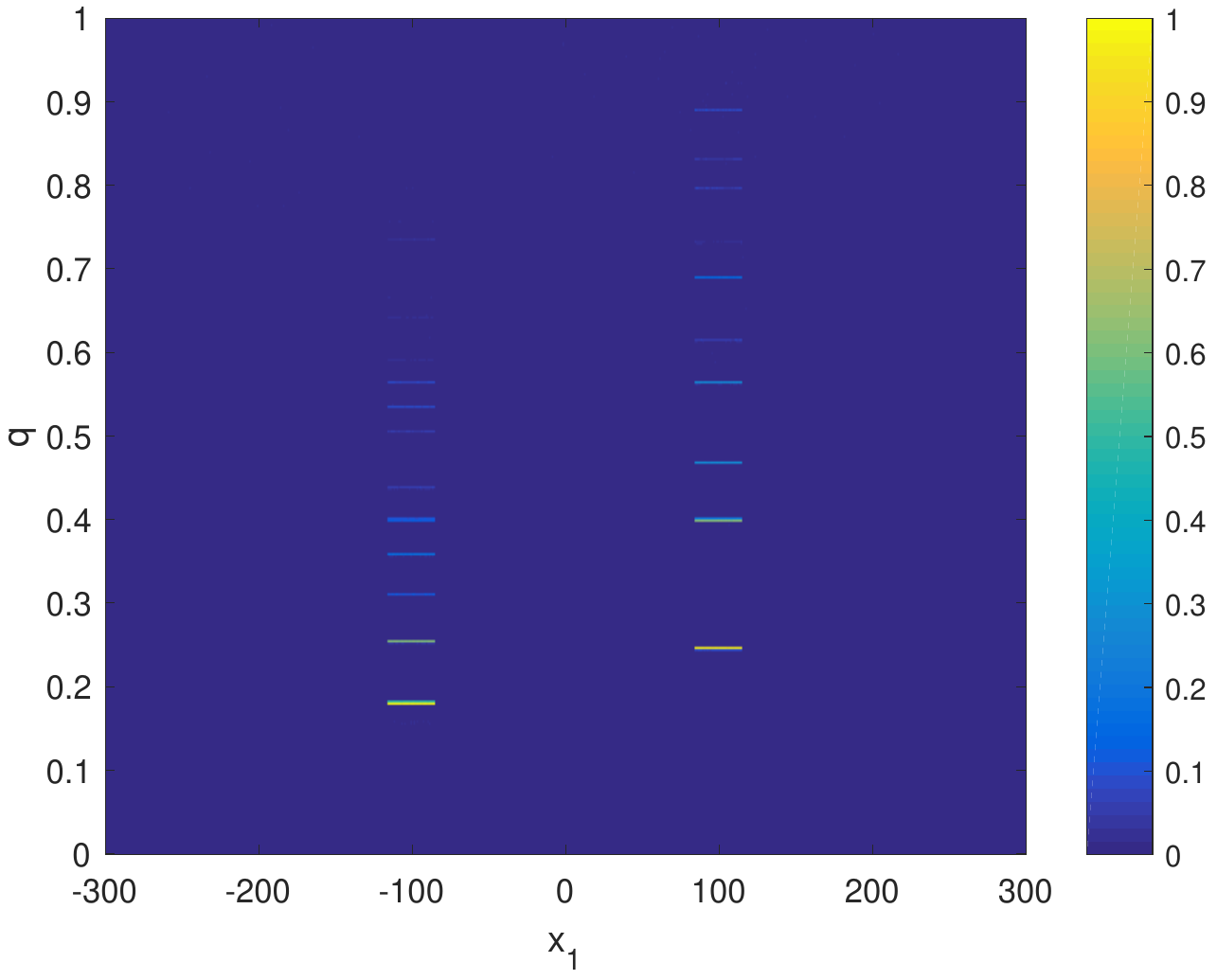} &
   \includegraphics[ width=0.4\linewidth, height=0.4\linewidth, keepaspectratio]{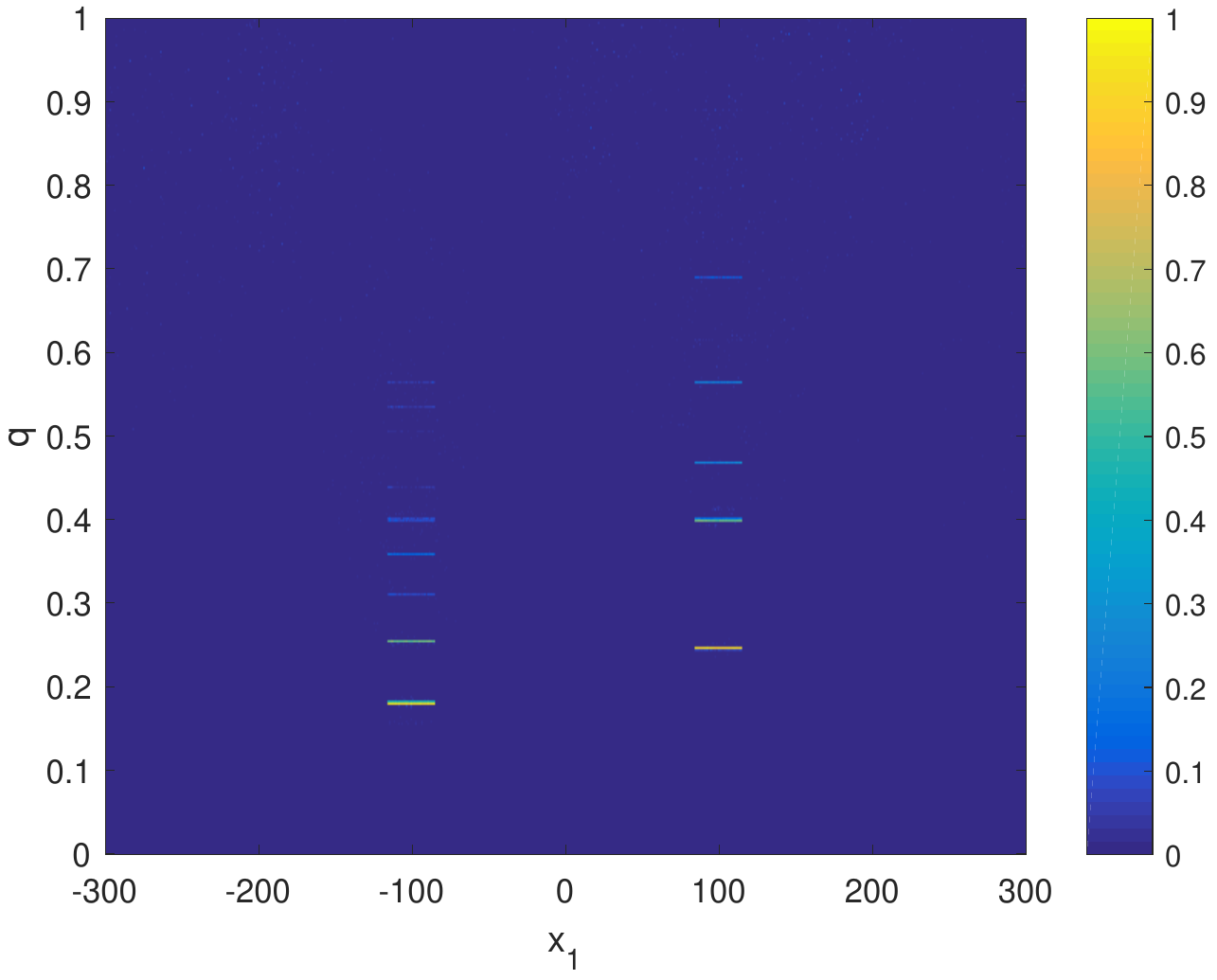} \\
 \rotatebox{90}{\hspace{1.5cm} $x_2=205$mm} &
   \includegraphics[ width=0.4\linewidth, height=0.4\linewidth, keepaspectratio]{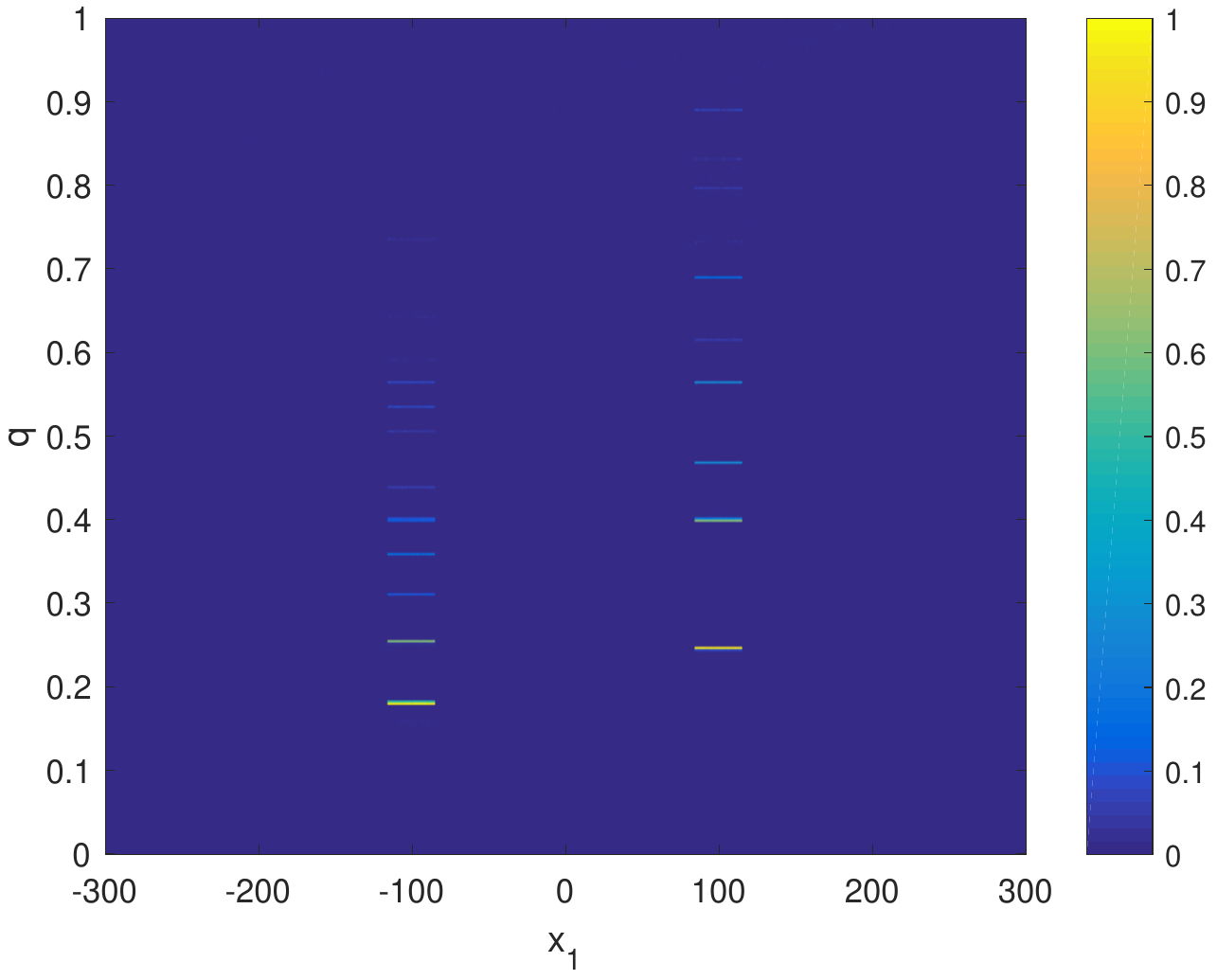} &
   \includegraphics[ width=0.4\linewidth, height=0.4\linewidth, keepaspectratio]{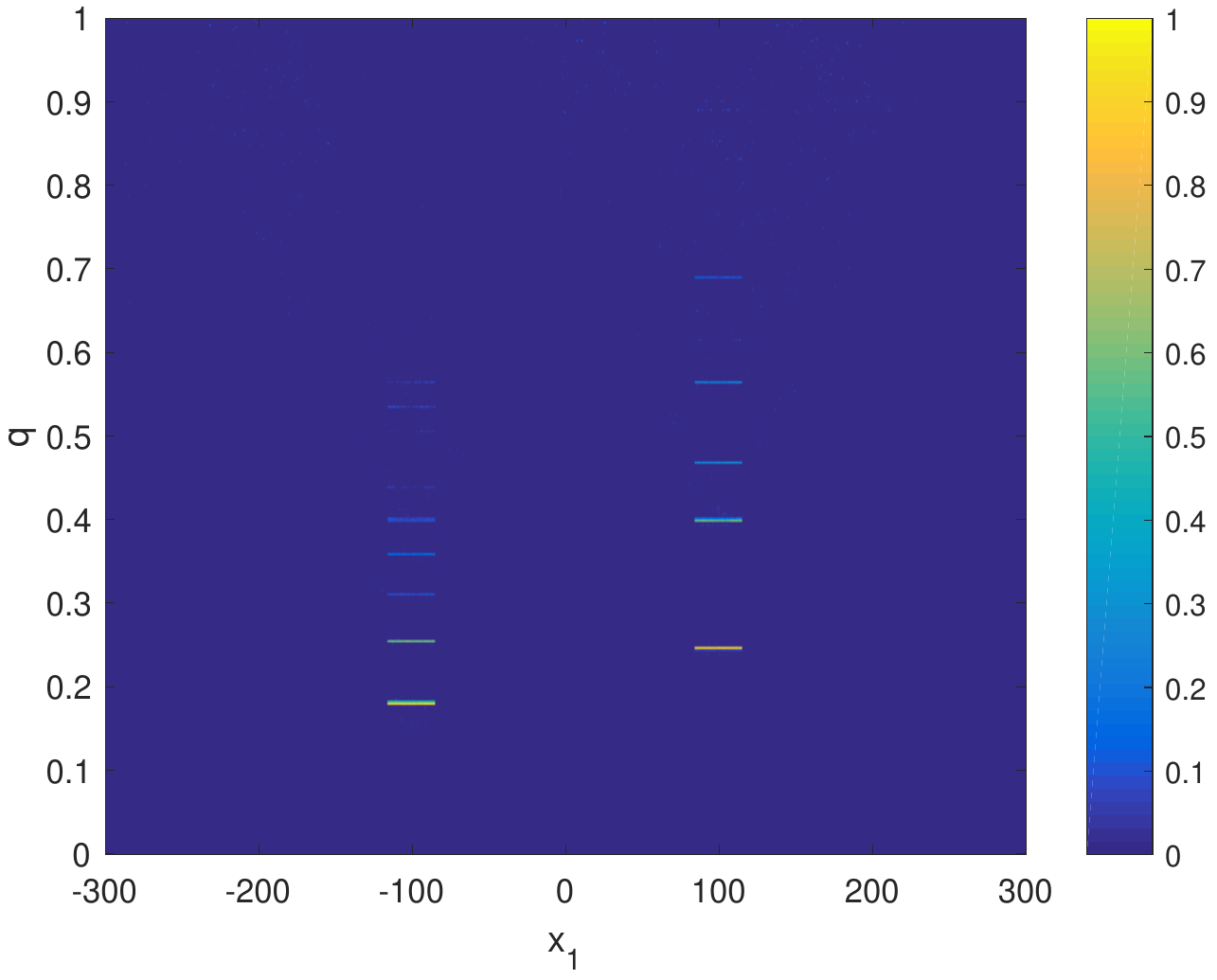}
\end{tabular}
\caption{2-spherical phantom reconstructions for all $(x_2,c_{\text{avg}})\in\{205,410,615\}\times\{1,10\}$, and the Ground Truth (GT) on the top row. To clarify, the GT does not vary with $c_{\text{avg}}$ and is included for comparison.}
\label{F1}
\end{figure}

The noisy data used in our simulations is generated as a single draw. Here $c_{\text{avg}}$ (the average photon counts per detector) controls the level of noise, i.e., larger $c_{\text{avg}}$ implies less noise and vise-versa. We consider two average count levels $c_{\text{avg}}=10$ and $c_{\text{avg}} =1$. To see how $c_{\text{avg}}$ equates to relative noise metrics, let $\epsilon_{\text{ls}}=\frac{\|\vb-\vb_{\textbf{e}}\|_2}{\|\vb\|_2}$ be the relative least square error. Then the values of  $\epsilon_{\text{ls}}$ for each experiment conducted, corresponding to each $(x_2,c_{\text{avg}})\in\{205,410,615\}\times\{1,10\}$ considered, are displayed in table \ref{T0}.
\begin{table}[!h]
\begin{subtable}{.49\linewidth}\centering
{
\begin{tabular}{| c | c | c | c | c | c | c | c |}
\hline
$\epsilon_{\text{ls}}$ & $c_{\text{avg}}=10$ & $c_{\text{avg}}=1$     \\ \hline
$x_2=0$     & $0.18$ & $0.58$  \\ 
$x_2=205$     & $0.18$  & $0.56$ \\  
$x_2=615$     & $0.18$ & $0.57$ \\  \hline
\end{tabular}}
\caption{2-spherical phantom}
\end{subtable}
\begin{subtable}{.49\linewidth}\centering
{
\begin{tabular}{| c | c | c | c | c | c | c | c |}
\hline
$\epsilon_{\text{ls}}$ & $c_{\text{avg}}=10$ & $c_{\text{avg}}=1$     \\ \hline
$x_2=0$     & $0.23$ & $0.73$  \\ 
$x_2=205$     & $0.22$  & $0.71$ \\  
$x_2=615$     & $0.22$ & $0.71$ \\  \hline
\end{tabular}}
\caption{4-spherical phantom}
\end{subtable}
\caption{Relative least square error values $\epsilon_{\text{ls}}$ for all experiments conducted. The values in the left-hand columns give the $x_2$ coordinates (in mm) of the sphere centers and the equations of the scanning line profiles (as illustrated in figure \ref{figPhan}).}
\label{T0}
\end{table}

\subsection{The reconstruction method}
\label{recon}
We aim to solve the classical linear inverse problem $A\vy =\vb_{\textbf{e}}$. Specifically we aim to find
\begin{equation}
\label{equTV}
\argmin_{\vy\in [0,\infty)^{nm}}\sum_k\left[\paren{A\vy}_k -\paren{\vb_{\textbf{e}}}_k\log\paren{A\vy}_k\right]+\lambda\text{TV}_{\beta}(\vy),
\end{equation}
where $\paren{\vb_{\textbf{e}}}_k$ and $\paren{A\vy}_k$ are the $k^{\text{th}}$ entries of $\vb_{\textbf{e}}$ and $A\vy$, respectively, and
\begin{equation}
\text{TV}_{\beta}(\vy)=\sum_{j=1}^m\sum_{i=1}^n\paren{\paren{\nabla Y}_{ij}^2+\beta^2}^{\frac{1}{2}}
\end{equation}
is a smoothed TV norm, where $Y\in\mathbb{R}^{n\times m}$ is $\vy$ reshaped into an $n\times m$ matrix (image), and $\nabla Y$ is the gradient image. The negative Poisson log-likelihood function in the first term of \eqref{equTV} is included since we expect the photon arrivals to follow a Poisson noise model, as is, for example, used in \cite{hassan2016snapshot,greenberg2013snapshot,maccabe2012pencil}. The hyperparameter $\beta>0$ is included so that the gradient of $\text{TV}_{\beta}$ is defined at zero, and the smoothing parameter $\lambda$ controls the level of TV regularization. TV regularization is applied with good results on Bragg scatter data in \cite{greenberg2013snapshot,maccabe2012pencil}, and thus why it is chosen as a regularization approach. To solve \ref{equTV} we implement the code ``JR\_PET\_TV" of \cite{ehrhardt2014joint} with non-negativity constraints, which is applied in that paper to Positron Emission Tomography (PET). That is, we input the noisy Bragg integral data $\vb_{\textbf{e}}$ to JR\_PET\_TV, with the constraints $\vy\in [0,\infty)^{nm}$. 

\begin{table}[!h]
\begin{subtable}{.49\linewidth}\centering
{
\begin{tabular}{| c | c | c | c | c | c | c | c |}
\hline
Mean $F_1$ score & $c_{\text{avg}}=10$ & $c_{\text{avg}}=1$     \\ \hline
$x_2=0$     & $0.90$ & $0.84$  \\ 
$x_2=205$     & $0.90$  & $0.90$ \\  
$x_2=615$     & $0.91$ & $0.90$ \\  \hline
\end{tabular}}
\caption{2-spherical phantom}
\end{subtable}
\begin{subtable}{.49\linewidth}\centering
{
\begin{tabular}{| c | c | c | c | c | c | c | c |}
\hline
Mean $F_1$ score & $c_{\text{avg}}=10$ & $c_{\text{avg}}=1$     \\ \hline
$x_2=0$     & $0.86$ & $0.67$  \\ 
$x_2=205$     & $0.89$  & $0.70$ \\  
$x_2=615$     & $0.87$ & $0.73$ \\  \hline
\end{tabular}}
\caption{4-spherical phantom}
\end{subtable}
\caption{Mean $F_1$ score results. The values in the left-hand columns give the $x_2$ coordinates (in mm) of the sphere centers and the equations of the scanning line profiles (as illustrated in figure \ref{figPhan}).}
\label{T1}
\end{table}

\subsection{Results}
The image reconstructions presented in this section are obtained in the following way. For  each line profile $x_2\in\{205,410,615\}$ and count level $c_{\text{avg}}\in\{1,10\}$ considered, we perform a set of image reconstructions for all $\lambda\in\{\frac{j}{10} : 1\leq j\leq9\}\cup\{j : 1\leq j\leq 10\}$ and $\beta\in\{0.001,0.01,0.1\}$ ($19\times 3=57$ reconstructions for each $(x_2,c_{\text{avg}})$) using the reconstruction method detailed in section \ref{recon}, and present the results with $F_1$ score closest to the mean $F_1$ score over all $\lambda,\beta$. The mean $F_1$ score can be considered as a lower bound for the performance as it is the expected result after picking $\lambda$ and $\beta$ at random (i.e. with no method of hyperparameter selection). $\lambda\in [0.1,10]$ and $\beta\in[0.001,0.1]$ were chosen by experimentation to be reasonable ranges for hyperparameter selection. 

See table \ref{T1} for our results in terms of mean $F_1$ score on the 2-spherical and 4-spherical phantom, and see figures \ref{F1} and \ref{F2} for the corresponding image reconstructions in all six (6) imaging scenarios considered (i.e. for $(x_2,c_{\text{avg}})\in\{205,410,615\}\times\{1,10\}$).
\begin{figure}
\setlength{\tabcolsep}{5pt}
\begin{tabular}{ c|cc }  
  &  $c_{\text{avg}}=10$ & $c_{\text{avg}}=1$ \\ \hline \\[-0.4cm] 
\rotatebox{90}{\hspace{2.5cm} GT}   & 
   \includegraphics[ width=0.4\linewidth, height=0.4\linewidth, keepaspectratio]{4sph} &
   \includegraphics[ width=0.4\linewidth, height=0.4\linewidth, keepaspectratio]{4sph} \\
 \rotatebox{90}{\hspace{1.75cm} $x_2=0$mm}  & 
   \includegraphics[ width=0.4\linewidth, height=0.4\linewidth, keepaspectratio]{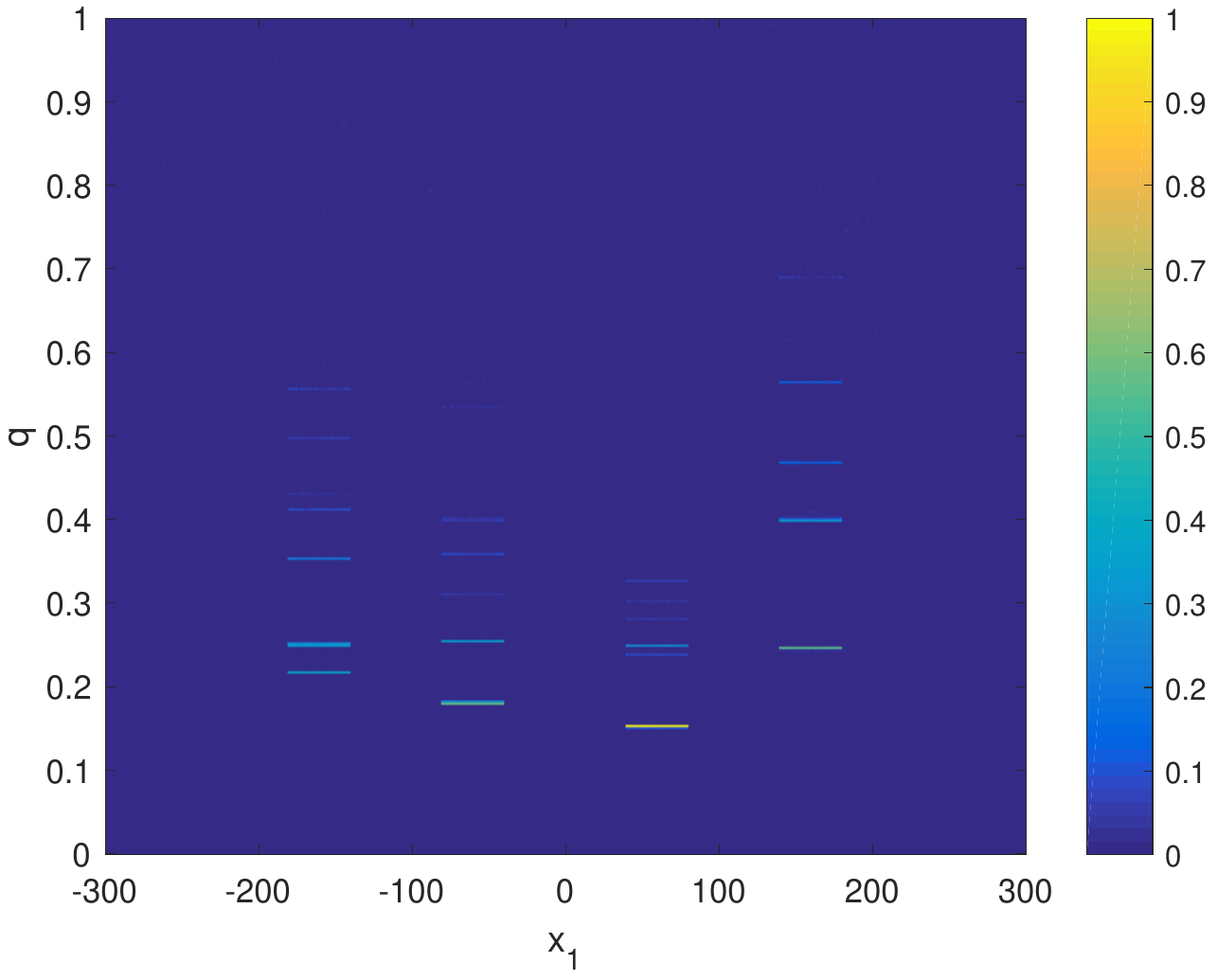} &
   \includegraphics[ width=0.4\linewidth, height=0.4\linewidth, keepaspectratio]{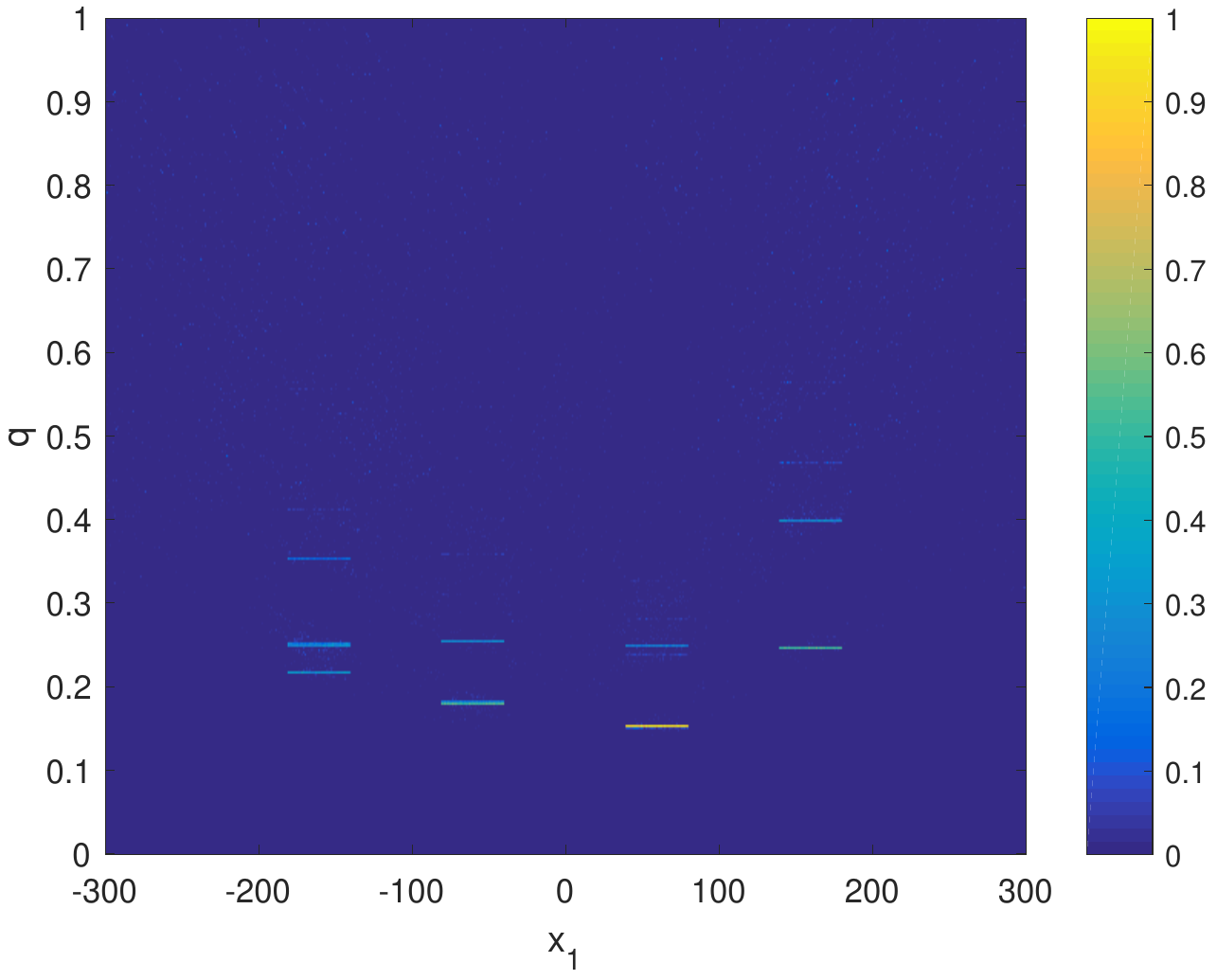} \\
 \rotatebox{90}{\hspace{1.5cm} $x_2=615$mm} &
  \includegraphics[ width=0.4\linewidth, height=0.4\linewidth, keepaspectratio]{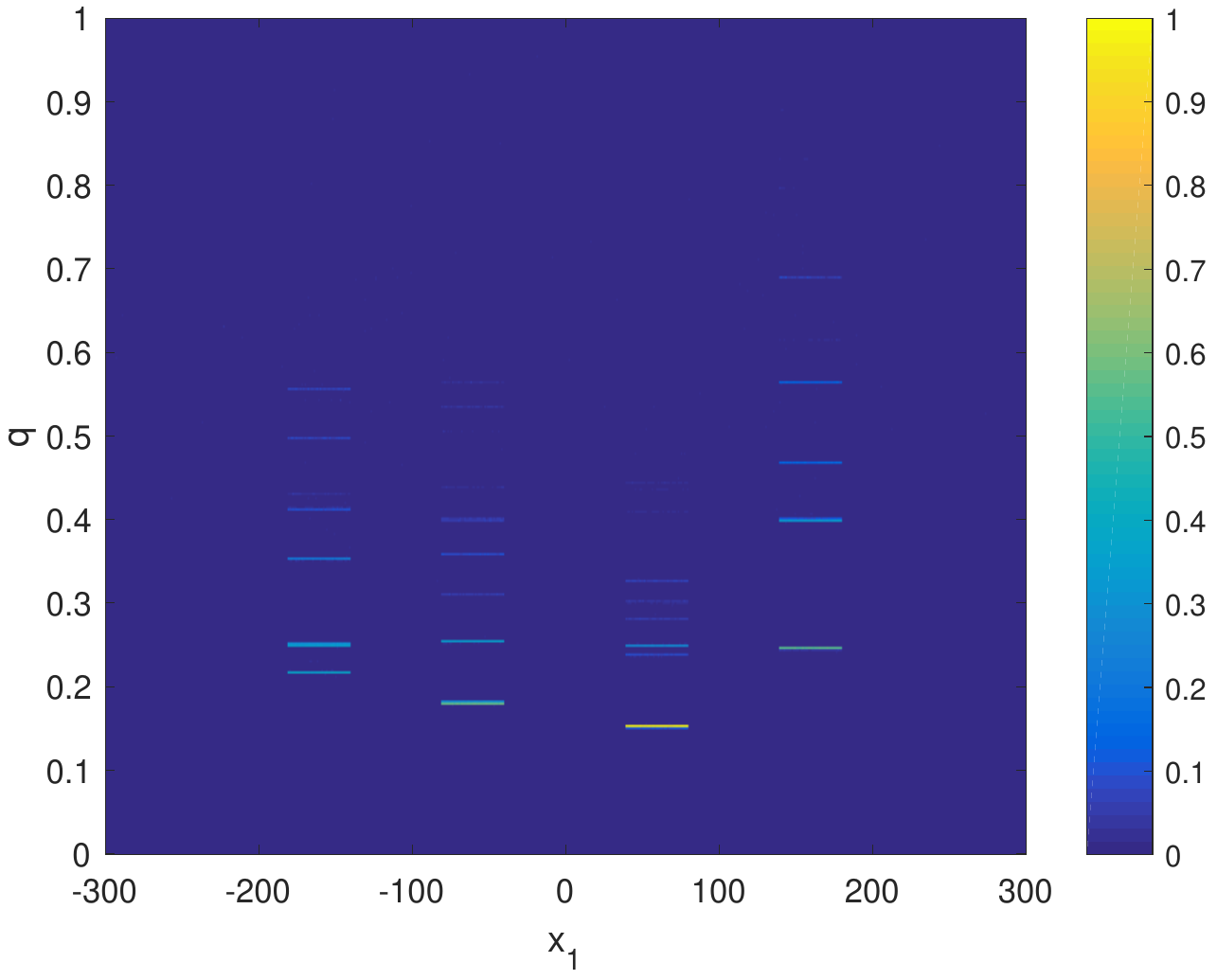} &
   \includegraphics[ width=0.4\linewidth, height=0.4\linewidth, keepaspectratio]{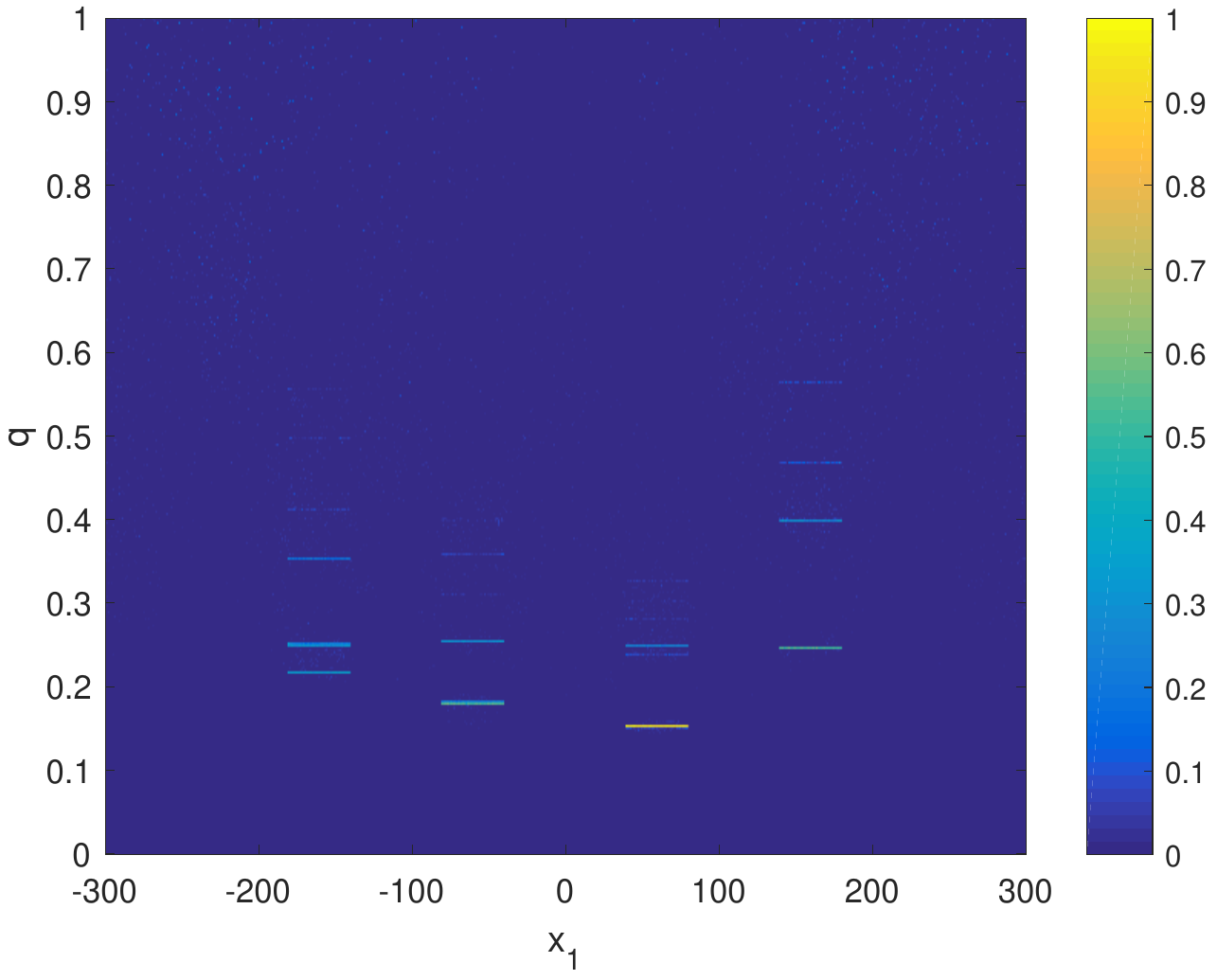} \\
 \rotatebox{90}{\hspace{1.5cm} $x_2=205$mm} &
   \includegraphics[ width=0.4\linewidth, height=0.4\linewidth, keepaspectratio]{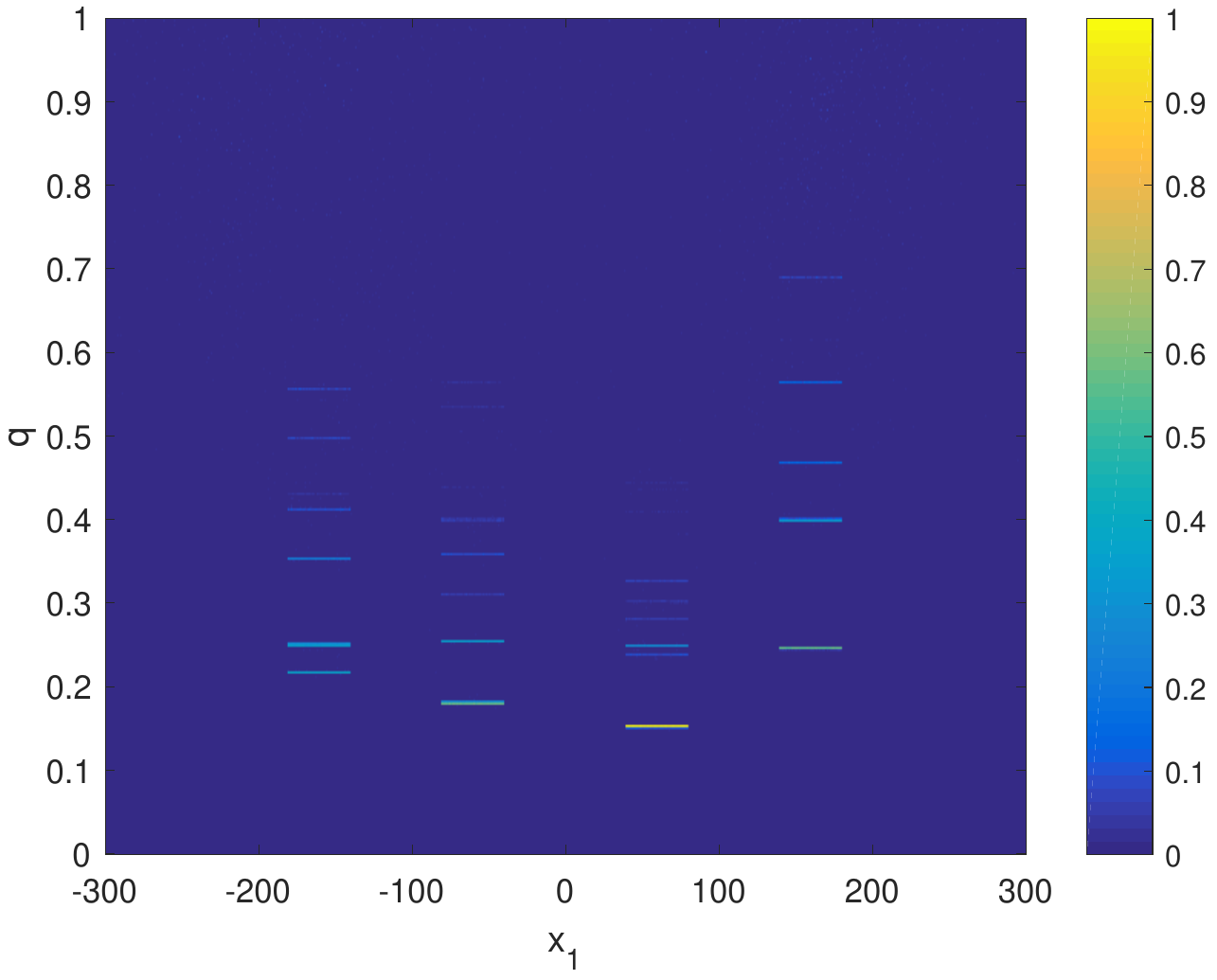} &
   \includegraphics[ width=0.4\linewidth, height=0.4\linewidth, keepaspectratio]{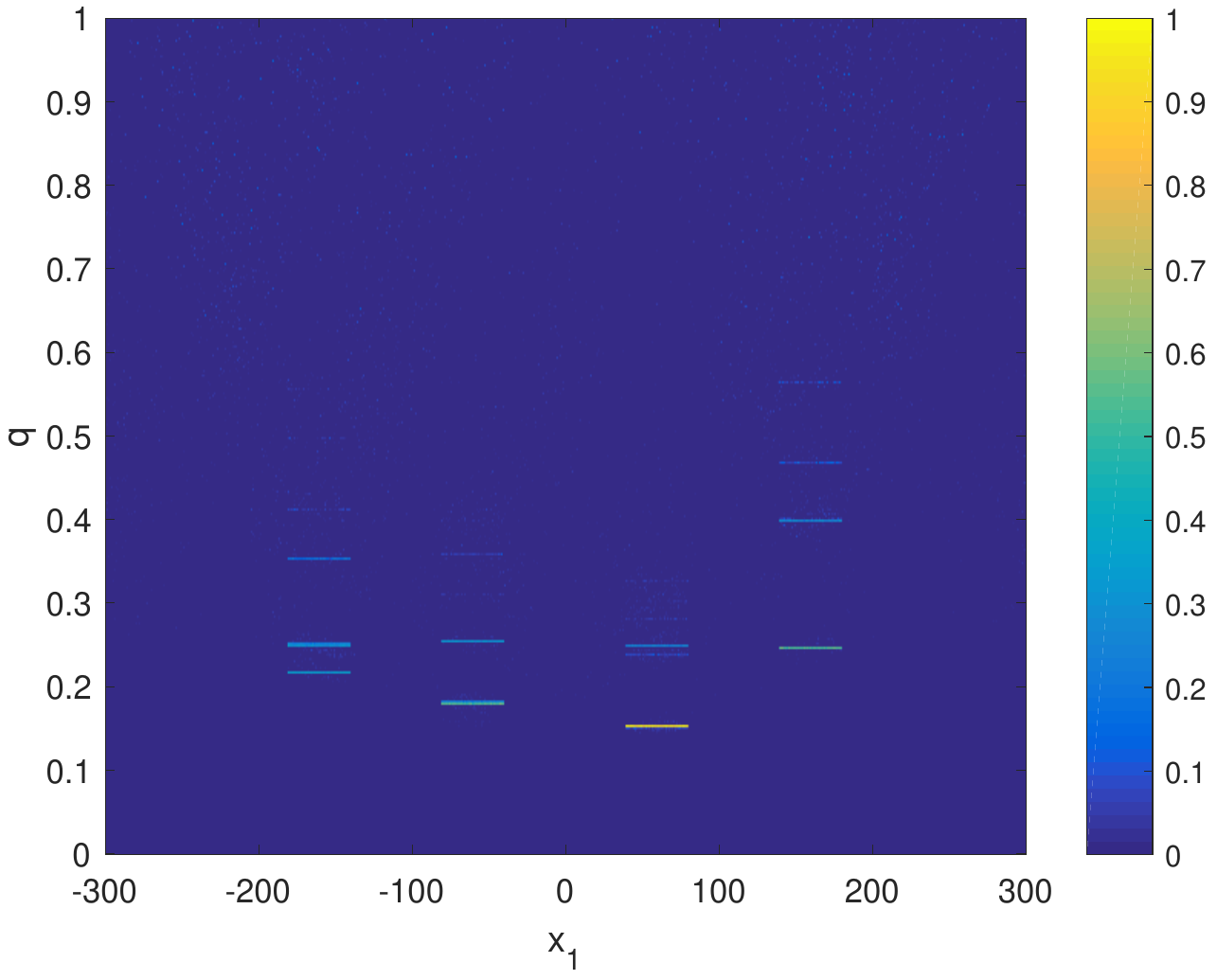}
\end{tabular}
\caption{4-spherical phantom reconstructions for all $(x_2,c_{\text{avg}})\in\{205,410,615\}\times\{1,10\}$, and the Ground Truth (GT) on the top row. To clarify, the GT does not vary with $c_{\text{avg}}$ and is included for comparison.}
\label{F2}
\end{figure}

The 2-spherical phantom results are good, with high image quality and mean $F_1$ score, over all $x_2$ and $c_{\text{avg}}$, and we see little variation in the quality of the results as $x_2$ varies. 
The analysis in Remark \ref{rem1} would suggest a difference in stability with varying $x_2$ but we do not observe this here in the simulations conducted. We notice that the image degradation is greater in the top half of the images (i.e. for $q>0.5$). This is because resolution of the image edges, and hence the problem stability, decreases with increasing $q$ \cite[Theorem 3.10]{webber2020microlocal} (see also figure 2 of that paper). Also the Bragg peaks in the $q>0.5$ space are notably smaller in magnitude than those in $q<0.5$ space, and thus contribute less to the overall variance in the data. This makes the smaller Bragg peaks harder to detect and recover. We see similar effects in conventional X-ray CT where a metal or high density object can dominate the reconstruction, and the more subtle features (e.g. soft tissue) are harder to identify \cite[figure 1(c)]{ehrhardt2014joint}. 

The image quality and $F_1$ scores in the 4-spherical phantom reconstructions are lower than that of the 2-spherical phantom, and the level of image degradation is greater overall, particularly at higher noise levels (i.e. when $c_{\text{avg}}=1$). In the $c_{\text{avg}}=1$ case we notice a high level of distortion in the $q>0.5$ half space and the smaller Bragg peaks are not recovered accurately. The more significant Bragg peaks in the $q<0.5$ space are recovered accurately however in all cases. The numerical results presented here are thus in line with the microlocal analysis of \cite{webber2020microlocal}, and verify the main theorems presented here, as we do not see the presence of any artifacts due to null space. Furthermore, the reconstruction method and TV regularization scheme presented appears to be effective in suppressing the boundary-type artifacts observed in \cite{webber2020microlocal} in Landweber and Filtered Back-Projection (FBP) reconstructions.

}

\section{Conclusions and further work}
We have presented new {injectivity results} for generalized Radon transforms in $\mathbb{R}^{n+1}$, which describe the integrals of $L^2_0$ functions over the $n$-dimensional surfaces of revolution of a $C^2$ curve class. We first introduced the Bragg transform $\mathfrak{B}$ in the $n=1$ case in section \ref{bgsec1}, which has a motivating application to a BST problem in airport baggage screening. Here we showed that the Bragg intensity could be modelled by $\mathfrak{B}f$, where $f$ is the differential cross section. We then went on to prove the injectivity of $\mathfrak{B}$ in Theorem \ref{main1}. Here we found evidence that the source width $w$ affected the stability of the inversion, in the sense that smaller $w$ allowed for a greater sampling rate in the frequency space and hence a more stable inversion (the converse effect being true for larger $w$). This was discussed in Remark \ref{rem1}. {A more rigorous stability analysis of the Bragg operators is left for future work, where we aim to build upon the theory of \cite{webber2020microlocal}, and derive estimates for the noise amplification using Sobolev space theory, as is, for example, done for the classical hyperplane Radon transform in \cite[pages 42 and 94]{natterer}.}
In section \ref{ndim} we introduced $\mathfrak{B}_n$, which generalizes the Bragg integration to $n>1$. {Injectivity results} for $\mathfrak{B}_n$ were provided in Theorem \ref{main2}. While these are promising results, which show the invertibility of the BST problem in the linearized case, moving forward we aim to consider the inversion properties of the non-linear models (e.g. including the effects due to attenuation).

We gave example machine parameters for the portal design in section \ref{MD}, chosen so that the injectivity conditions \eqref{md} of Theorem \ref{main4} were satisfied. Here we discovered that, as $w$ increased, the upper bounds on $\epsilon$ increased allowing for more freedom in the machine design. In future work we aim to show how one can choose the machine parameters so as to optimize the stability of the Bragg inversion process, under the constraints of invertible design (i.e. condition \eqref{md}).

{In section \ref{results} we presented image reconstruction from $\mathfrak{B}_{a}f$ data using a discrete (algebraic) approach and TV regularization. Overall our method performed well, and was shown to offer a high level of image quality and $F_1$ score in the presence of significant noise ($c_{\text{avg}}=10$ or $18\%$ minimum least squares error). Even at higher noise levels ($c_{\text{avg}}=1$ or $56\%$ minimum least squares error), while the recovery in $q>0.5\AA^{-1}$ space suffered due to lower stability, the more significant Bragg peaks in $q<0.5\AA^{-1}$ space were well recovered, and thus we can expect a good level of materials characterization using our approach.}
\vspace{-0.05cm}
\section*{Acknowledgements}
We would like to thank Professor Eric Todd Quinto for his helpful suggestions and insight towards the article. In particular towards Theorem \ref{main2} and the proof of analyticity of $J$ in \eqref{J(u)}.
This material is based upon work supported by the U.S.\ Department of
Homeland Security, Science and Technology Directorate, Office of
University Programs, under Grant Award 2013-ST-061-ED0001. The views
and conclusions contained in this document are those of the authors
and should not be interpreted as necessarily representing the official
policies, either expressed or implied, of the U.S. Department of
Homeland Security.

\appendix
\section{Calculation of derivatives}
\label{appA}
Here we show additional steps towards the calculation of the derivatives of $g$ and $q_1$ of equations \eqref{gBST} and \eqref{q13} respectively. 
\vspace{-0.05cm}
\subsection{The derivative of $g$}
\label{appA1}
Let $f_1(z)=\frac{1}{4z^2(1-z^2)}-x_2$ and let $f_2(z)=\frac{1-2z^2}{2z\sqrt{1-z^2}}$. Then $g(z)=\sqrt{f_1(z)}-f_2(z)$. Using the product rule, the derivative of $f_2$ is
\begin{equation}
\begin{split}
f'_2(z)&=\frac{-4z}{2z\sqrt{1-z^2}}+(1-2z^2)\paren{-\frac{1}{2z^2\sqrt{1-z^2}}+\frac{z}{2z(1-z^2)^{\frac{3}{2}}}}\\
&=\frac{-4z^2(1-z^2)-(1-2z^2)^2}{2z^2(1-z^2)^{\frac{3}{2}}}=-\frac{1}{2z^2(1-z^2)^{\frac{3}{2}}}.
\end{split}
\end{equation}
We have $f'_1(z)=\frac{2z^2-1}{2z^3(1-z^2)^2}$. Let $f_3(z)=\sqrt{1-4x_2^2z^2(1-z^2)}$. Then
\begin{equation}
\begin{split}
g'(z)&=\frac{2z^2-1}{4z^3(1-z^2)^2\sqrt{f_2(z)}}+\frac{1}{2z^2(1-z^2)^{\frac{3}{2}}}\\
&=\frac{z(2z^2-1)}{2z^3(1-z^2)^{\frac{3}{2}}f_3(z)}+\frac{zf_3(z)}{2z^3(1-z^2)^{\frac{3}{2}}f_3(z)}\\
&=\frac{f_3(z)-(1-2z^2)}{2z^2(1-z^2)^{\frac{3}{2}}f_3(z)}=g(z)h(z),
\end{split}
\end{equation}
where $h$ is as in equation \eqref{g'=gh}.
\vspace{-0.05cm}
\subsection{The derivative of $q_1$}
\label{appA2}
Let $f_1(x_1)=\frac{\sqrt{1+x_1}}{\sqrt{2}}$, let $f_2(x_1)=x_1^2-(1-x_2^2)$ and let
$f_3(x_1)=\frac{1}{h_1(x_1)}$, where $h_1$ is as in equation \eqref{hBST1}. Then $q_1=f_1\circ (f_2f_3)$ and
\begin{equation}
\begin{split}
q'_1&=\paren{f'_2f_3+f'_3f_2}\left[f'_1\circ (f_2f_3)\right]=\frac{\paren{f'_2f_3+f'_3f_2}}{4q_1}.
\end{split}
\end{equation}
The derivatives of $f_2$ and $f_3$ are $f'_2(x_1)=2x_1$ and
$$f'_3(x_1)=-\frac{x_1\paren{\Phi^2(x_2)+2(x_1^2+x_2^2+1)}}{h^3_1(x_1)}.$$
It follows that
\begin{equation}
\label{Papp}
\begin{split}
q'_1(x_1)=\frac{P_1(x_1)}{4q_1(x_1)h^3_1(x_1)},
\end{split}
\end{equation}
{where, after substituting $\epsilon=\Phi(x_2)$
\begin{equation}
\label{P1exp}
\begin{split}
\frac{P_1(x_1)}{x_1}&=2h^2_1(x_1)-\paren{\epsilon^2+2(x_1^2+x_2^2+1)}f_2(x_1)\\
&=\underbrace{2\paren{x_1^2+(1+x_2)^2}\paren{x_1^2+(1-x_2)^2+\epsilon^2}}_{2h^2_1(x_1)}\\
&\hspace{3cm}-\paren{\epsilon^2+2(x_1^2+1+x_2^2)}\underbrace{\paren{x_1^2-(1-x_2^2)}}_{f_2(x_1)}\\
&=2x_1^4+2x^2_1\paren{(1+x_2)^2+(1-x_2)^2}+2\epsilon^2x^2_1+2\epsilon^2(1+x_2)^2\\
&\hspace{3cm}+2(1+x_2)^2(1-x_2)^2 -\epsilon^2x^2_1+\epsilon^2(1-x^2_2)-2x^4_1\\
&\hspace{4cm}-2x^2_1(1+x^2_2)+2x^2_1(1-x^2_2)+2(1+x^2_2)(1-x^2_2)\\
&=\underbrace{4x^2_1(1+x^2_2)}_{2x^2_1\paren{(1+x_2)^2+(1-x_2)^2}}+\epsilon^2x^2_1+2\epsilon^2(1+x_2)^2+2(1+x_2)^2(1-x_2)^2\\
&\hspace{3cm} +\epsilon^2(1-x^2_2)\underbrace{-4x^2_1x^2_2}_{-2x^2_1(1+x^2_2)+2x^2_1(1-x^2_2)}+2(1+x^2_2)(1-x^2_2)\\
&=P_{\epsilon^2}(x_1)+P_{\sim \epsilon^2}(x_1),
\end{split}
\end{equation}
where $P_{\epsilon^2}$ denotes the terms in the expansion with $\epsilon^2$ coefficients and $P_{\sim \epsilon^2}$ denotes the remaining terms. Some of the simplifications and expansions in \eqref{P1exp} are highlighted using underbraces. We have
\begin{equation}
\begin{split}
P_{\epsilon^2}(x_1)&=\epsilon^2\paren{x^2_1+2(1+x_2)^2+(1-x^2_2)}\\
&=\epsilon^2\paren{x^2_1+3+4x_2+x^2_2}\\
&=\epsilon^2\paren{x^2_1+(x_2+1)(x_2+3)},
\end{split}
\end{equation}
and
\begin{equation}
\begin{split}
P_{\sim \epsilon^2}(x_1)&=4x^2_1(1+x^2_2)+2(1+x_2)^2(1-x_2)^2-4x^2_1x^2_2+2(x^2_2+1)(1-x^2_2)\\
&=4x^2_1+4x^2_1x^2_2+2-4x^2_2+2x^4_2-4x^2_1x^2_2-2x^4_2+2\\
&=4(1-x^2_2+x^2_1).
\end{split}
\end{equation}
Thus it follows that $P_1$ of equation \eqref{P1exp} reduces to the polynomial in equation \eqref{PBST1}.} This completes the derivation of $q'_1$.

\medskip
Received xxxx 20xx; revised xxxx 20xx.
\medskip

\medskip
\textit{E-mail address}: \href{mailto:james.webber@tufts.edu}{james.webber@tufts.edu}

\textit{E-mail address}: \href{mailto:eric.miller@tufts.edu}{eric.miller@tufts.edu}
\medskip

\end{document}